\documentclass[10pt]{article}

\usepackage{a4wide}
\usepackage{amssymb}
\usepackage{amsfonts}
\usepackage{amsmath}
\input xy
\xyoption{arrow} \xyoption{matrix}

\date{}

\newtheorem{proposition}{Proposition}[section]
\newtheorem{theorem}[proposition]{Theorem}
\newtheorem{lemma}[proposition]{Lemma}

\newtheorem{corollary}[proposition]{Corollary}

\def\Hom{{\rm Hom}}
\def\der{\partial }

\def\nFM0{{\nu }_{F,M_0}}
\def\nFN0{{\nu }_{F,N_0}}
\def\nGN0{{\nu }_{G,N_0}}

\def\N0{ {\bf N}_0 }

\def\t{\otimes}

\def\v{\varphi}
\def\ra{\rightarrow}

\def\lra{\leftrightarrow}
\def\Xpm{X^{\pm }}

\def\s{\sigma}
\def\Z{\mathbb{Z}}

\def\l1{{\lambda}_1}

\def\a{\alpha}
\def\a0{ {\alpha }_0}
\def\a1{ {\alpha }_1}

\def\l{\lambda}

\def\nFGM0{{\nu }_{F,G,M_0}}

\def\nFN0{{\nu}_{F,N_0}}

\def\sm{{\sigma}^m}

\def\sm1{{\sigma}^{-1}}

\def\smtp1{{\sigma}^{-t+1}}

\def\S1{S^{-1}}

\def\Xpm1{X^{\pm 1}_1}

\def\sPM1{{\sigma }^{\pm 1}}
\def\sMP1{{\sigma }^{\mp 1 }}

\def\d{\delta}

\def\di{{\rm d.ind}}

\def\L{\Lambda}

\def\OO{{\cal O}}

\def\CD{{\cal D}}

\def\Ytm1{Y^{t-1}}
\def\Yim1{Y^{i-1}}

\def\CK{{\cal K}}
\def\CL{{\cal L}}

\def\CF{{\cal F}}
\def\CG{{\cal G}}
\def\CH{{\cal H}}
\def\ass{{\rm ass}}

\def\supp{{\rm supp }}

\def\Aut{{\rm Aut}}
\def\bK{\overline{K}}

\def\bA{\overline{A}}

\def\ad{{\rm ad }}
\def\dim{{\rm dim }}

\def\Iso{{\rm Iso }}

\def\ker{ {\rm ker } }

\def\SL2Z{ {\rm SL}_2({\bf Z}) }

\def\CR{ {\cal R}}

\def\th{ \theta }

\def\CL{{\cal L}}

\def\Gp1{ G^{1 , 1 } }
\def\P11{ P^{-1 , 1 } }
\def\Pp1{ P^{1 , 1 } }

\def\th{\theta}
\def\CE{{\cal E}}
\def\CV{{\cal V}}

\def\nCLsr{{}^\nu\kern-2pt {\cal L}^{\sigma , \rho  }}
\def\nP{{}^\nu \kern-2pt P}
\def\nL{{}^\nu\kern-2pt L}
\def\nLL{{}^\nu\kern-2pt \Lambda}
\def\nPsr{{}^\nu\kern-2pt P^{\sigma , \rho  }}
\def\nLsr{{}^\nu\kern-2pt L^{\sigma , \rho  }}
\def\nuCL{{}^\nu\kern-2pt  {\cal L}}
\def\nCLsr{{}^\nu\kern-2pt {\cal L}^{\sigma , \rho  }}
\def\nCL1m{{}^\nu\kern-2pt {\cal L}^{-1 , 1  }}
\def\x1nu{x^\frac{1}{\nu}}
\def\xm1nu{x^{-\frac{1}{\nu}}}

\def\rad{{\rm rad}}

\def\CR{ {\cal R}}

\def\ra{\rightarrow }

\def\CB{{\cal B}}

\def\coker{{\rm coker}}
\def\CT{{\cal T}}

\def\CC{ {\cal C}}
\def\CE{ {\cal E} }
\def\CH{ {\cal H}}
\def\CP{ {\cal P}}

\def\nAM0{{\nu }_{{\cal A},M_0}}
\def\nAN0{{\nu }_{{\cal A},N_0}}

\def\End{ {\rm End }}

\def\CR{ {\cal R }}
\def\CP{ {\cal P }}
\def\det{ {\rm det }}
\def\ad{ {\rm ad }}

\def\bI{\overline{I}}

\def\bR{\overline{R}}

\def\tr{{\rm tr}}

\def\ga{\mathfrak{a}}
\def\gb{\mathfrak{b}}
\def\gc{\mathfrak{c}}

\def\gp{\mathfrak{p}}
\def\gq{\mathfrak{q}}

\def\bJ{\overline{J}}

\def\GL{{\rm GL}}
\def\SL{{\rm SL}}

\def\Hom{{\rm Hom}}

\def\di!{\frac{\der^i}{i!}}
\def\dik!{\frac{\der^k_i}{k!}}
\def\hA{\widehat{A}}

\def\Max{{\rm Max}}

\def\N{\mathbb{N}}

\def\0{\overline{0}}
\def\1{\overline{1}}

\def\Ln1{\L_{n,\overline{1}}}

\def\oa{\overline{a}}

\def\a1{a_{\overline{1}}}

\def\S{\Sigma}

\def\sign{{\rm sign}}
\def\vn1{\overrightarrow{n-1}}

\def\im{{\rm im}}

\def\soc{{\rm soc}}

\def\mJ{\mathbb{J}}
\def\mI{\mathbb{I}}

\def\ann{{\rm ann}}
\def\lann{{\rm l.ann}}
\def\rann{{\rm r.ann}}
\def\Cen{{\rm Cen}}

\def\mT{\mathbb{T}}
\def\ind{{\rm ind}}

\def\Frac{{\rm Frac}}

\def\K1{{\rm K}_1}

\def\lind{{\rm l.ind}}
\def\rind{{\rm r.ind}}
\def\hmI1{\widehat{\mI}_1}
\def\tmI1{\widetilde{\mI}_1}
\def\tmJ1{\widetilde{\mJ}_1}
\def\hB1{\widehat{B}_1}
\def\hCB1{\widehat{\CB}_1}
\def\codim{{\rm codim}}
\def\bCF{\overline{\CF}}
\def\ffrac{{\rm frac}}
\def\bS{\overline{S}}

\def\linv{{\rm l.inv}}
\def\rinv{{\rm r.inv}}

\begin{document}

\author{V. V. \  Bavula 
}

\title{The algebra of integro-differential operators
on an affine line and its modules}

\maketitle

\begin{abstract}

For the algebra $\mI_1= K\langle x, \frac{d}{dx}, \int \rangle$ of
polynomial integro-differential operators over a field $K$ of
characteristic zero, a classification of  simple modules  is
given. It is proved that  $\mI_1$ is a left and  right coherent
algebra. The {\em Strong Compact-Fredholm Alternative} is proved
for $\mI_1$.
The endomorphism algebra of each simple $\mI_1$-module is a {\em
finite dimensional} skew field. In contrast to the first Weyl
algebra, the centralizer of a non-scalar  integro-differential
operator can be a noncommutative, non-Noetherian, non-finitely
generated algebra which is not a domain.
 It is proved that neither left nor
right quotient ring of $\mI_1$ exists but there exists the {\em
largest left quotient ring } and the {\em largest right quotient
ring} of $\mI_1$,
 they are not
 $\mI_1$-isomorphic but $\mI_1$-{\em anti-isomorphic}. Moreover,
 the factor ring of the largest right quotient ring modulo its
 only proper ideal is isomorphic to the quotient ring of the first
 Weyl algebra. An analogue of the Theorem of Stafford (for the Weyl
 algebras) is proved for $\mI_1$: each {\em finitely generated}
 one-sided ideal of $\mI_1$ is $2$-generated.


 {\em Key Words: the algebra of polynomial integro-differential
 operators, the Strong Compact-Fredholm Alternative, coherent algebra,
  the Weyl algebra, simple module, compact
 integro-differential operators,  Fredholm operators, centralizer,
 the largest left/right quotient ring.
}

 {\em Mathematics subject classification
 2000:   16D60, 16D70, 16P50, 16U20.}

\end{abstract}


\section{Introduction}

Throughout, ring means an associative ring with $1$; module means
a left module;
 $\N :=\{0, 1, \ldots \}$ is the set of natural numbers; $K$ is a
field of characteristic zero and  $K^*$ is its group of units;
$P_1:= K[x]$ is a polynomial algebra in one variable $x$ over $K$;
$\der:=\frac{d}{d x}$; $\End_K(P_1)$ is the algebra of all
$K$-linear maps from $P_1$ to $P_1$,  and $\Aut_K(P_1)$ is its
group of units (i.e.,  the
 group of all the invertible linear maps from $P_1$ to $P_1$); the
subalgebras  $A_1:= K \langle x , \der \rangle$ and
 $\mI_1:=K\langle x,\der ,  \int\rangle $ of $\End_K(P_1)$
  are called the (first) {\em Weyl algebra} and the {\em algebra  of polynomial
integro-differential operators} respectively  where $\int: P_1\ra
P_1$, $ p\mapsto \int p \, dx$, is   {\em integration},  i.e.,
$\int : x^n \mapsto \frac{x^{n+1}}{n+1}$ for all $n\in \N$.  The
algebra $\mI_1$ is neither left nor right Noetherian and not a
domain. Moreover, it contains infinite direct sums of nonzero left
and right ideals, \cite{algintdif}.

$\noindent $

In Section \ref{CSI1MOD}, a classification of simple
$\mI_1$-modules is given (Theorem \ref{31May10}), it is similar to
the classification of simple $A_1$-modules from
\cite{Bav-AlgAnaliz92} and to Block's classification \cite{Bl 1},
\cite{Bl 3} (Block was the first who did it for $A_1$ in his
seminal, long paper \cite{Bl 3}. Later  advancing Block's ideas
and using an approach of generalized Weyl algebras an alternative
classification was given in \cite{Bav-AlgAnaliz92} and the proof
comprises only several pages). One might expect (which is not
obvious from the outset)  a close connection between the simple
modules of the algebras $A_1$ and $\mI_1$. The reality is even
more surprising: the algebra $\mI_1$ has exactly `one less' simple
module than the algebra $A_1$ (Theorem \ref{7Jun10}). So,
surprisingly,  inverting on the right the derivation $\der$
(i.e.,  adding the integration $\int$ which is
 a {\em one-sided} inverse of $\der$, $\der \int =1$,  but not two-sided) `kills' only
a {\em single} simple $A_1$-module.

$\noindent $

In Section \ref{FDIMKERCO},  for modules and algebras two new
concepts are introduced:  the  Compact-Fredholm Alternative and
the Strong Compact-Fredholm Alternative. The Strong
Compact-Fredholm Alternative is proved for the algebra $\mI_1$
(Theorem \ref{B30May10}) which says that on {\em all} simple
$\mI_1$-modules the action of each element of $\mI_1$ is either
{\em compact} (i.e.,  the image is finite dimensional) or,
otherwise, {\em Fredholm} (i.e.,  the kernel and the cokernel are
finite dimensional). 
 So,
the algebra
$$\mI_1=\CC_{\mI_1}\coprod\CF_{\mI_1}$$ is the disjoint
union of the sets  of  compact  and Fredholm  operators respectively. The set of compact operators
$\CC_{\mI_1}$ coincides with the only proper (two-sided) ideal
$$F=\bigoplus_{i,j\in \N} K(\int^i\der^j-\int^{i+1}\der^{j+1})$$ of
the algebra $\mI_1$ (Theorem \ref{B30May10}).
The endomorphism algebra of each
simple $\mI_1$-module is a finite dimensional division algebra
(Theorem \ref{6Jun10}).
\begin{itemize}
\item (Theorem \ref{A30May10}) {\em Let $a\in \mI_1$, $\cdot
a:\mI_1\ra \mI_1$, $b\mapsto ba$, and  $l_{\mI_1}$ be the
length function on the set of left $\mI_1$-modules. Then}
\begin{enumerate}
\item $l_{\mI_1}(\ker (\cdot a))<\infty $ {\em iff}
$l_{\mI_1}(\coker (\cdot a))<\infty$  {\em iff} $a\not\in F$.
\item $l_{\mI_1}(\im (\cdot a))<\infty $  {\em iff }  $a\in F$.
\end{enumerate}
\end{itemize}
 Various
indices are introduced for non-compact integro-differential
operators, the $M$-{\em index} and the {\em length index} and it is
proved that they are invariant under addition of  compact operator
(Lemma \ref{D30May10} and Lemma \ref{C30May10}).

$\noindent $

In Section \ref{COHERENT}, it is proved that the algebra $\mI_1$
is a left and right coherent algebra (Theorem \ref{a19Sep10}).
 In particular, it is
proved that the intersection of finitely many  finitely
generated left (or right) ideals of $\mI_1$ is again a finitely
generated left (or right) ideal of $\mI_1$ (Theorem
\ref{19Sep10}).

\begin{itemize}
\item  (Theorem \ref{21Sep10}) {\em Let $a\in \mI_1$. Then}
\begin{enumerate}
\item $\ker_{\mI_1}(\cdot a)$ {\em and $\coker_{\mI_1}(\cdot a)$
are finitely generated left $\mI_1$-modules.} \item $\ker_{\mI_1}(
a\cdot )$ {\em  and $\coker_{\mI_1}( a\cdot )$ are finitely
generated right $\mI_1$-modules.}
\end{enumerate}
\end{itemize}
The Theorem of  Stafford  says that {\em each left or right ideal
of the Weyl algebra $A_n$ over a field of characteristic zero is
generated by two elements} \cite{Staford-1978-2gen}.  The same result holds for a larger class of algebras that includes the, so-called, generalized Weyl algebras \cite{Bav-Const-96}.  The simplicity of the Weyl algebra $A_n$ plays a crucial role in the
proof. An analogue of this result is proved for the
(non-Noetherian, non-simple) algebra $\mI_1$.
\begin{itemize}
\item (Theorem \ref{A19Sep10}) {\em Each finitely generated left
(or right) ideal of the algebra $\mI_1$ is generated by two
elements.}
\end{itemize}

Section \ref{CENTRALIZERS}: The centralizer of each non-scalar
element of the Weyl algebra $A_1$ is a finitely generated
commutative (hence Noetherian) domain (see Burchnall and Chaundy
\cite{Burchnall-Chaundy-1923} and Amitsur \cite{Amitsur-58}). In
contrast to the Weyl algebra $A_1$, the centralizer of a
non-scalar element of the algebra $\mI_1$ can be a non-finitely
generated, noncommutative, non-Noetherian algebra  which is not a
domain  (Proposition \ref{a1Jun10}.(1)). Theorem \ref{C11Oct10}
describes in great detail the structure of  centralizers of
non-scalar elements of the algebra $\mI_1$. Theorem \ref{C11Oct10} is too technical to be explained
in the Introduction, it is a generalization of the Theorem of
Amitsur on the centralizer for the first Weyl algebra $A_1$. We
state here only two corollaries of Theorem \ref{C11Oct10}.

\begin{itemize}
\item {\rm (Corollary \ref{xa10Oct10})} {\em  Let $a\in
\mI_1\backslash K$ and $H:=\der x$.  Then the following statements
are equivalent.}
\begin{enumerate}
\item $a\not\in K[H]+F$. \item  $\Cen_{\mI_1}(a)$ {\em is a
finitely generated $K[a]$-module. \item $\Cen_{\mI_1}(a)$ is a
left Noetherian algebra. \item $\Cen_{\mI_1}(a)$ is a right
Noetherian algebra. \item $\Cen_{\mI_1}(a)$ is a finitely
generated and  Noetherian algebra.}
\end{enumerate}

\item (Corollary \ref{a30Oct10}) {\em  Let $a\in \mI_1\backslash
K$. Then}
\begin{enumerate}
\item $\Cen_{\mI_1}(a)$ {\em is a finitely generated algebra iff
$a\not\in (K[H]+F)\backslash (K+F)$.
 \item  $\Cen_{\mI_1}(a)$ is a finitely
generated, not Noetherian/not left Noetherian/not right Noetherian
algebra iff} $a\in (K+F)\backslash K$. \item  The algebra
$\Cen_{\mI_1}(a)$ is not  finitely generated, not Noetherian/not
left Noetherian/not right Noetherian
 iff $a\in (K[H]+F)\backslash (K+F)$.
\end{enumerate}
\end{itemize}

In Section \ref{SIKXM}, for all elements $a\in \mI_1\backslash F$, an explicit formula for the index $\ind
(a_{K[x]})$ is found
(Proposition \ref{a12Jun10}.(1)) where $a_{K[x]}:K[x]\ra K[x]$,
$p\mapsto ap$. Classifications  are given of elements $a\in \mI_1$ that satisfy the following properties:  the map $a_{K[x]}$ is a bijection (Theorem
\ref{12Jun10}),   a surjection (Theorem \ref{15Jun10}),  an
injection (Theorem \ref{16Jun10}). In case when the map $a_{K[x]}$
is a bijection an explicit inversion formula is found (Theorem
\ref{12Jun10}.(4)). The kernel and the cokernel of the linear map
$a_{K[x]}$ are found in the cases when the map $a_{K[x]}$ is either  surjective
or injective (Theorem \ref{15Jun10} and Theorem \ref{16Jun10}).
\begin{itemize}
\item (Proposition \ref{a19Jun10})
\begin{enumerate}
\item {\em For each element $a\in \mI_1\backslash F$ with $n:=
\dim_K(\coker (a_{K[x]}))$, there exists an element $\der^n+f$ for
some $f\in F$ (resp. $s\in (1+F)^*$)   such that the map
$(\der^n+f)a_{K[x]}$ (resp. $s\der^n s^{-1}a_{K[x]}$) is a
surjection. In this case, $\ker((\der^n+f)a_{K[x]}) = \ker
(a_{K[x]})$ (resp. $\ker (s\der^ns^{-1}a_{K[x]})= \ker
(a_{K[x]}))$. } \item {\em For each element $a\in \mI_1\backslash
F$ with
 $n:= \dim_K(\ker (a_{K[x]}))$, there exists an element $\int^n +g$
  for
some $g\in F$ (resp. $s\in (1+F)^*$) such that the map $a(\int^n
+g)_{K[x]}$ (resp. $s\int^n s^{-1}a_{K[x]}$)  is an injection. In
this case, $\im (a(\int^n +g)_{K[x]})=\im ( a_{K[x]})$ (resp.
$\;\;\;$  $\im (s\int^ns^{-1}a_{K[x]})= \im (a_{K[x]}))$.} \item
{\em For each element $a\in \mI_1\backslash F$ with $m:=
\dim_K(\ker (a_{K[x]}))$ and $n :=\dim_K(\coker (a_{K[x]}))$, there
exist elements $\der^n+f$ and $\int^m+g$ for some $f,g\in F$
(resp. $s,t\in (1+F)^*$) such that the map $(\der^n+f) a(\int^m
+g)_{K[x]}$ (resp. $s\der^ns^{-1} at\int^mt^{-1}_{K[x]})$ is a
bijection.}
\end{enumerate}
\end{itemize}
Proposition \ref{a6Jun10} gives necessary and sufficient
conditions for each of the vector spaces $\ker_{\ker_{\mI_1}(\cdot
b)}(a\cdot )$, $\coker_{\ker_{\mI_1}(\cdot b)}(a\cdot )$,
 $\ker_{\coker_{\mI_1}(\cdot b)}(a\cdot )$ and   $ \coker_{\coker_{\mI_1}(\cdot b)}(a\cdot
 )$   to be finite
dimensional.

The algebra $K+F=K+M_\infty (K)$ admits the (usual) determinant
map $\det $ (see (\ref{detaKF})),  and the group of units
$\mI_1^*$ of the algebra $\mI_1$ is equal to $\{ u\in K+F\, | \,
\det (u)\neq 0\}$, \cite{algintdif}.

For an element $u\in \mI_1$, let $\linv (u) :=\{ v\in \mI_1\, | \,
vu=1\}$ and $\rinv (u) :=\{ v\in \mI_1\, | \, uv=1\}$, the {\em
sets of left and right inverses} for the element $u$. In 1942,
Baer \cite{Baer-inverses-1942} and,  in 1950,  Jacobson
\cite{Jacobson-onesidedinv-1950} begun to study one sided
inverses.  The next theorem describes all the left and right
inverses of elements in $\mI_1$.
\begin{itemize}
\item (Corollary \ref{c25Sep10})
\begin{enumerate}
\item An element $a\in \mI_1$ admits a left inverse iff $a=a'
\int^n$ for some natural number $n\geq 0$ and an element $a'\in K^*+F$
such that  $a_{K[x]}$ is an injection (necessarily,
$n=\dim_K(\coker(a_{K[x]})))$. In this case, $\linv (a) = \{ b\in
\der^n b'\, | \, b'\in K^*+F, \int^n \der^n b'a'\int^n
\der^n=\int^n \der^n \}$. $\CL (\mI_1) = \{ a\in (K^*+F)\int^n\, |
\, n\in \N, a_{K[x]}$ is an injection$\}$.
 \item An element $b\in \mI_1$ admits a right inverse iff  $b\in
\der^n b'$ for some natural number $n\geq 0$ and an element $b'\in K^*+F$
such that  $b_{K[x]}$ is a surjection (necessarily,
$n=\dim_K(\ker(b_{K[x]})))$. In this case, $\rinv (b) = \{
a'\int^n \, | \, a'\in K^*+F, \int^n \der^n b'a'\int^n
\der^n=\int^n \der^n \}$. $\CR (\mI_1) =\{ b\in \der^n (K^*+F)\, |
\, n\in \N, b_{K[x]}$ is a surjection$\}$.
\end{enumerate}
\end{itemize}

In Section \ref{COSI}, it is proved that the monoid $\CL
(\mI_1):=\{ a\in \mI_1\, | \, ba=1$ for some $b\in \mI_1\}$ of
left invertible elements of the algebra $\mI_1$ is generated by
the group of units $\mI_1^*$ of the algebra $\mI_1$ and the element $\int$.
Moreover, $\CL (\mI_1) = \bigsqcup_{i\in \N} \mI_1^*\int^n$, the
disjoint union (Theorem \ref{3Oct10}.(1,2)). Similarly, the monoid
$\CR (\mI_1):=\{ b\in \mI_1\, | \, ba=1$ for some $a\in \mI_1\}$
of right invertible elements of the algebra $\mI_1$ is generated
by the group
 $\mI_1^*$ and the element $\der$. Moreover, $\CR
(\mI_1) = \bigsqcup_{i\in \N} \der^n \mI_1^*$, the disjoint union
(Theorem \ref{3Oct10}.(3,4)). For each left invertible element
$a\in \CL (\mI_1)$ (respectively, right invertible element $b\in
\CR (\mI_1)$) the set of its left (respectively, right) inverses
is  found (Theorem \ref{A3Oct10}).

In Section \ref{TAI1J1}, we introduce two algebras $\tmI1$ and
$\tmJ1$ that turn out to be the {\em largest right quotient ring}
$\Frac_r(\mI_1)$  and the {\em largest left quotient ring}
$\Frac_l(\mI_1)$ of the algebra $\mI_1$ (Theorem \ref{18Jun10}).
The algebra $\tmI1$ is the subalgebra of $\End_K(K[x])$ generated
by the algebra $\mI_1$ and the (large) set
$$\mI_1^0:= \mI_1\cap
\Aut_K(K[x]),$$ the set of all the elements of $\mI_1$ that are
{\em invertible} linear maps in $K[x]$; the group of units
$\mI_1^*$ of the algebra $\mI_1$ is a small part of the monoid
$\mI_1^0$. In particular, the set $\mI_1^0$ is the {\em largest}
(w.r.t. inclusion) {\em regular right
 Ore set} in $\mI_1$  but it is  not a left Ore set of $\mI_1$,
and $\tmI1 = \mI_1{\mI_1^0}^{-1}$ (Theorem \ref{18Jun10}.(4)). The
algebras $\tmI1$ and $\tmJ1$ contain the only proper ideal, $\CC
(\tmI1 )$ and $\CC (\tmJ1 )$,  respectively (precisely the ideal of
compact operators in $\tmI1$ and $\tmJ1$ respectively). The factor
algebras $\tmI1 /\CC (\tmI1 )$ and $\tmJ1 / \CC (\tmJ1 )$ are
canonically isomorphic to the skew field of fractions $\Frac
(A_1)$ of the first Weyl algebra $A_1$ and its opposite skew field
$\Frac (A_1)^{op}$ respectively (Theorem \ref{14Jun10}.(4) and
Corollary  \ref{a14Jun10}.(4)).
 The simple (left and right) modules for the algebras $\tmI1$ and
 $\tmJ1$ are classified, the groups of units of the algebras
 $\tmI1$ and $\tmJ1$ are found (Proposition \ref{A21Jun10}).

 In Section \ref{LLRQR}, it is proved that the rings $\tmI1$ and
 $\tmJ1$ are not $\mI_1$-isomorphic but $\mI_1$-{\em
 anti-isomorphic}. The sets of right regular, left regular and
 regular elements of the algebra $\mI_1$ are described  (Lemma
\ref{b18Jun10}.(1), Corollary \ref{c18Jun10}.(1) and Corollary
\ref{b2Jul10}).

$\noindent $

{\bf The Conjecture/Problem of Dixmier} (1968) [still open]: {\em
is an algebra endomorphism of the Weyl algebra $A_1$ an
automorphism?} It turns out that an analogue of the
Conjecture/Problem of Dixmier is true for the algebra $\mI_1$:

\begin{itemize}
\item {\rm (Theorem 1.1, \cite{Bav-cdixintdif})} {\em Each algebra
endomorphism of the algebra $\mI_1$ is an automorphism.}
\end{itemize}

The present paper is instrumental in proving this result.

The algebras $\mI_n:=\mI_1^{\t n}$ ($n\geq 1$) of polynomial integro-differential operators are studied in detail in \cite{algintdif}, their groups of automorphisms  are found in \cite{intdifaut}.


\section{Classification of simple $\mI_1$-modules}\label{CSI1MOD}

In this section, a classification of simple $\mI_1$-modules is
given (Theorem \ref{31May10}) and it is compared to a similar
classification of simple modules over the Weyl algebra $A_1$
(Theorem \ref{7Jun10}).

The algebra $\mI_1$  is generated by the elements $\der $, $H:=
\der x$ and $\int$ (since $x=\int H$) that satisfy the defining
relations (Proposition 2.2, \cite{algintdif}): $$\der \int = 1,
\;\; [H, \int ] = \int, \;\; [H, \der ] =-\der , \;\; H(1-\int\der
) =(1-\int\der ) H = 1-\int\der .$$ The elements of the algebra
$\mI_1$,  
\begin{equation}\label{eijdef}
e_{ij}:=\int^i\der^j-\int^{i+1}\der^{j+1}, \;\; i,j\in \N ,
\end{equation}
satisfy the relations $e_{ij}e_{kl}=\d_{jk}e_{il}$ where $\d_{jk}$
is the Kronecker delta function. Notice that
$e_{ij}=\int^ie_{00}\der^j$. The matrices of the linear maps
$e_{ij}\in \End_K(K[x])$ with respect to the basis $\{ x^{[s]}:=
\frac{x^s}{s!}\}_{s\in \N}$ of the polynomial algebra $K[x]$  are
the elementary matrices, i.e.,
$$ e_{ij}*x^{[s]}=\begin{cases}
x^{[i]}& \text{if }j=s,\\
0& \text{if }j\neq s.\\
\end{cases}$$
Let $E_{ij}\in \End_K(K[x])$ be the usual matrix units, i.e.,
$E_{ij}*x^s= \d_{js}x^i$ for all $i,j,s\in \N$. Then
\begin{equation}\label{eijEij}
e_{ij}=\frac{j!}{i!}E_{ij},
\end{equation}
 $Ke_{ij}=KE_{ij}$, and
$F:=\bigoplus_{i,j\geq 0}Ke_{ij}= \bigoplus_{i,j\geq
0}KE_{ij}\simeq M_\infty (K)$, the algebra (without 1) of infinite
dimensional matrices. Using induction on $i$ and the fact that
$\int^je_{kk}\der^j=e_{k+j, k+j}$, we can easily  prove that
\begin{equation}\label{Iidi}
\int^i\der^i = 1-e_{00}-e_{11}-\cdots - e_{i-1,
i-1}=1-E_{00}-E_{11}-\cdots -E_{i-1,i-1}, \;\; i\geq 1.
\end{equation}


{\bf A $\Z$-grading on the algebra $\mI_1$ and the canonical form of
an integro-differential operator, \cite{algintdif}}. The algebra
$\mI_1=\bigoplus_{i\in \Z} \mI_{1, i}$ is a $\Z$-graded algebra
($\mI_{1, i} \mI_{1, j}\subseteq \mI_{1, i+j}$ for all $i,j\in
\Z$) where 
\begin{equation}\label{I1iZ}
\mI_{1, i} =\begin{cases}
D_1\int^i=\int^iD_1& \text{if } i>0,\\
D_1& \text{if }i=0,\\
\der^{|i|}D_1=D_1\der^{|i|}& \text{if }i<0,\\
\end{cases}
\end{equation}
 the algebra $D_1:= K[H]\bigoplus \bigoplus_{i\in \N} Ke_{ii}$ is
a commutative non-Noetherian subalgebra of $\mI_1$, $ He_{ii} =
e_{ii}H= (i+1)e_{ii}$  for $i\in \N $ (notice that
$\bigoplus_{i\in \N} Ke_{ii}$ is the direct  sum of non-zero
ideals of $D_1$); $(\int^iD_1)_{D_1}\simeq D_1$, $\int^id\mapsto
d$; ${}_{D_1}(D_1\der^i) \simeq D_1$, $d\der^i\mapsto d$,   for
all $i\geq 0$ since $\der^i\int^i=1$.
 Notice that the maps $\cdot\int^i : D_1\ra D_1\int^i$, $d\mapsto
d\int^i$,  and $\der^i \cdot : D_1\ra \der^iD_1$, $d\mapsto
\der^id$, have the same kernel $\bigoplus_{j=0}^{i-1}Ke_{jj}$.

Each element $a$ of the algebra $\mI_1$ is the unique finite sum
\begin{equation}\label{acan}
a=\sum_{i>0} a_{-i}\der^i+a_0+\sum_{i>0}\int^ia_i +\sum_{i,j\in
\N} \l_{ij} e_{ij}
\end{equation}
where $a_k\in K[H]$ and $\l_{ij}\in K$. This is the {\em canonical
form} of the polynomial integro-differential operator
\cite{algintdif}.

$\noindent $

{\it Definition}. Let $a\in \mI_1$ be as in (\ref{acan}) and let
$a_F:=\sum \l_{ij}e_{ij}$. Suppose that $a_F\neq 0$ then
\begin{equation}\label{degFa}
\deg_F(a) :=\min \{ n\in \N \, | \, a_F\in \bigoplus_{i,j=0}^n
Ke_{ij}\}
\end{equation}
is called the $F$-{\em degree} of the element $a$;
$\deg_F(0):=-1$.

$\noindent $

Let $v_i:=\begin{cases}
\int^i& \text{if }i>0,\\
1& \text{if }i=0,\\
\der^{|i|}& \text{if }i<0.\\
\end{cases}$

Then $\mI_{1,i}=D_1v_i= v_iD_1$ and an element $a\in \mI_1$ is the
unique  finite  sum 
\begin{equation}\label{acan1}
a=\sum_{i\in \Z} b_iv_i +\sum_{i,j\in \N} \l_{ij} e_{ij}
\end{equation}
where $b_i\in K[H]$ and $\l_{ij}\in K$. So, the set $\{ H^j\der^i,
H^j, \int^iH^j, e_{st}\, | \, i\geq 1; j,s,t\geq 0\}$ is a
$K$-basis for the algebra $\mI_1$. The multiplication in the
algebra $\mI_1$ is given by the rule:
$$ \int H = (H-1) \int , \;\; H\der = \der (H-1), \;\; \int e_{ij}
= e_{i+1, j}, \;\; e_{ij}\int= e_{i,j-1}, \;\; \der e_{ij}=
e_{i-1, j}\;\; e_{ij} \der = \der e_{i, j+1}.$$
$$ He_{ii} = e_{ii}H= (i+1)e_{ii}, \;\; i\in \N, $$
where $e_{-1,j}:=0$ and $e_{i,-1}:=0$.  $\noindent $

$\noindent $

{\bf The ideal $F$ of compact operators of $\mI_1$}.
 Let $V$ be an infinite dimensional vector space over a field $K$.
A linear map $\v \in \End_K(V)$ is called a {\em compact} linear
map/operator if $\dim_K(\im (\v ))<\infty$. The set $\CC =\CC (V)$
of all compact operators is a (two-sided) ideal of the algebra
$\End_K(V)$.
 The algebra
$\mI_1$ has the only proper ideal
$$F=\bigoplus_{i,j\in \N}Ke_{ij}
\simeq M_\infty (K),$$ the ideal of compact operators in $\mI_1$,
$F=\mI_1\cap \CC (K[x])$ (Corollary \ref{aB30May10}),
 $F^2= F$. The factor algebra $\mI_1/F$ is canonically isomorphic to
the skew Laurent polynomial algebra $B_1:= K[H][\der, \der^{-1} ;
\tau ]$, $\tau (H) = H+1$, via $\der \mapsto \der$, $ \int\mapsto
\der^{-1}$, $H\mapsto H$ (where $\der^{\pm 1}\alpha = \tau^{\pm
1}(\alpha ) \der^{\pm 1}$ for all elements $\alpha \in K[H]$). The
algebra $B_1$ is canonically isomorphic to the (left and right)
localization $A_{1, \der }$ of the Weyl algebra $A_1$ at the
powers of the element $\der$ (notice that $x= \der^{-1} H$).
Therefore, they have the common skew field of fractions, $\Frac
(A_1) = \Frac (B_1)$, the {\em first Weyl skew field}. The algebra
$B_1$ is a subalgebra of the skew Laurent polynomial algebra
$\CB_1:= K(H) [ \der, \der^{-1} ; \tau ]$ where $K(H)$ is the
field of rational functions over the field $K$ in $H$. The algebra
$\CB_1 = S^{-1} B_1$ is the left and right localization of the
algebra $B_1$ at the multiplicative set $S=K[H]\backslash \{ 0\}$.
The algebra $\CB_1$ is a {\em noncommutative Euclidean domain},
i.e.,  the left and right division algorithms with remainder hold
with respect to the length function $l$ on $B_1$:
$$l(\alpha_m
\der^m+\alpha_{m+1}\der^{m+1}+\cdots +\alpha_n\der^n)=n-m$$ where
$\alpha_i\in K(H)$, $\alpha_m\neq 0$, $\alpha_n\neq 0$, and
$m<\cdots <n$. In particular, the algebra $\CB_1$ is a principal
left and right ideal domain. A $\CB_1$-module $M$ is simple iff
$M\simeq \CB_1/\CB_1b$ for some irreducible element $b\in \CB_1$,
and $\CB_1/\CB_1b\simeq \CB_1/\CB_1c$ iff the elements $b$ and $c$
are {\em similar} (that is,  there exists an element $d\in \CB_1$
such that $1$ is the greatest common right divisor of $c$ and $d$,
and $bd$ is the least common left multiple of $c$ and $d$).

$\noindent $

{\bf The involution $*$ on the algebra $\mI_1$}. The algebra
$\mI_1$ admits the involution $*$ over the field $K$:
$$\der^* =
\int,\;\;\;  \int^* = \der \;\;\; {\rm  and }\;\;\; H^* = H,$$ i.e.,  it is a $K$-algebra
{\em anti-isomorphism} ($(ab)^* = b^* a^*$) such that $a^{**} =a$.
Therefore, the algebra $\mI_1$ is {\em self-dual}, i.e.,  it is
isomorphic to its opposite algebra $\mI_1^{op}$. As a result, the
left and right properties of the algebra $\mI_1$ are the same.
Clearly, $e_{ij}^* = e_{ji}$ for all $i,j\in \N$, and so $F^* =
F$.

$\noindent $

{\bf A classification of  simple $\mI_1$-modules}. Since the field
$K$ has characteristic zero, the group $\langle \tau \rangle
\simeq \Z$ acts freely on the set $\Max (K[H])$ of maximal ideals
of the polynomial algebra $K[H]$. That is, for each maximal ideal
$\gp \in \Max (K[H])$, its orbit $\OO (\gp ) :=\{ \tau^i (\gp ) \,
| \, i\in \Z \}$ contains infinitely many elements. For two
elements $\tau^i (\gp )$ and $\tau^j(\gp )$ of the orbit $\OO (\gp
)$, we write $\tau^i (\gp ) <\tau^j (\gp )$ if $i<j$. Let $\Max
(K[H])/\langle \tau \rangle$ be the set of all $\langle \tau
\rangle$-orbits in $\Max (K[H])$. If $K$ is an algebraically
closed field then $\gp = (H-\l )$, for some scalar $\l \in K$,
$\Max (K[H])\simeq K$, and $\Max(K[H])/\langle \tau \rangle\simeq
K/ \Z$.

For elements $\alpha , \beta \in K[H]$, we write $\alpha <\beta$
if for all maximal ideals $\gp , \gq \in\Max (K[H])$ that belong
to the same orbit and such that $\alpha \in \gp $ and $\beta \in
\gq$ we have $\gp <\gq$. If $K$ is an algebraically closed field
then $p <q$ iff $\l -\mu>0$ for all roots $\l$ and $\mu$ of the
polynomials $p$ and $q$ respectively (if they exist)  such that
$\l -\mu \in \Z$.

$\noindent $

{\it Definition}. An element $b= \der^{-m}\beta_{-m}+\cdots
+\beta_0\in B_1$, where all $\beta_i\in K[H]$, $\beta_{-m}\neq 0$,
$\beta_0\neq 0$,  and $m>0$, is called a {\em normal} element if
$\beta_0<\beta_{-m}$.

$\noindent $

Let the element $b= \der^{-m}\beta_{-m}+\cdots +\beta_0\in B_1$ be such
as in the Definition above but not necessarily normal. Then {\em
there exist polynomials $\alpha , \beta \in K[H]$ such that the
element $\beta b \alpha^{-1}\in B_1$ is normal}. For example,
$$ \alpha = \prod \{ \tau^i (\beta_0)\, | \, -s\leq i \leq 0\},
\;\;\;\;\; \beta = \prod \{ \tau ^j (\beta_0)\, | \, -m-s\leq j
\leq 1\}$$ where $s\in \N$ such that
$\tau^{_s}(\beta_0)<\beta_{-m}$ and $\tau^{-s}(\beta_0)<\beta_0$
(Proposition 16, \cite{Bav-GWA-Ker-TensoSimple-92}).

For an algebra $A$, let  $\hA$ be the set of isomorphism
classes of simple $A$-modules and, for  a simple $A$-module $M$, let
$[M]\in \hA$ be  its isomorphism class. For   an
isomorphism invariant property $\CP$ of simple $A$-modules, let
$$\hA (\CP ) := \{ [M]\in \hA\, | \, M\;\; {\rm  has \; property}\;\; \CP\}.$$ The
{\em socle} $\soc (M)$ of a module $M$ is the sum of all the
simple submodules if they exist and zero otherwise. Since the
algebra $\mI_1$ contains the Weyl algebra $A_1$,  which is a
simple infinite dimensional algebra,  each nonzero  $\mI_1$-module
is necessarily an {\em infinite dimensional} module.

The algebra $B_1=K[H](\tau, 1)$, $\tau (H) = H+1$,  is an example
of a generalized Weyl algebra (GWA) $A= D(\s , a)$ over a {\em
Dedekind} domain $D$ which is in the case of the algebra $B_1$ is
the polynomial algebra $K[H]$. Many popular algebras of small
Gelfand-Kirillov dimension are examples of $A$ (eg, the Weyl
algebra $A_1$, the {\em quantum plane}  $\L = K\langle x,y\, | \,
xy = \l yx\rangle$, $\l \in K^*$; the {\em quantum Weyl algebra}
$A_1(q) = K\langle x,\der \, | \, \der x-qx\der =1\rangle$, $q \in
K\backslash \{ 0,1\}$, all prime infinite dimensional factor
algebras of $Usl(2)$, etc). A classification of simple $A$-modules
is obtained in \cite{Bav-AlgAnaliz92, Bav-UkrMathJ-92,
Bav-GWA-Ker-TensoSimple-92, Bav-Voyst-JA-1997}.

Since the algebra $B_1$ is a factor algebra of $\mI_1$, there is
the tautological embedding
$$ \hB1\ra \hmI1, \;\; [ M]\mapsto [M].$$
 Therefore, $\hB1\subseteq \hmI1$.

\begin{theorem}\label{31May10}
{\rm (Classification of simple $\mI_1$-modules)}
\begin{enumerate}
\item $\hmI1 = \{ [K[x]]\}\coprod \hB1$.  \item The set $\hB1 =
\hB1 (K[H]-{\rm torsion})\coprod  \hB1 (K[H]-{\rm torsion \;
free})$ is the disjoint union where $\hB1 (K[H]-{\rm torsion}):=\{
[M]\in \hB1\, | \, S^{-1} M=0\}$, $\hB1 (K[H]-{\rm torsion \;
free}):=\{ [M]\in \hB1\, | \, S^{-1}M\neq 0\}$, and
$S:=K[H]\backslash \{ 0\}$. \item The map $$\Max (K[H])/ \langle
\tau \rangle\ra \hB1 (K[H]-{\rm torsion}), \;\; [\gp ] \mapsto [
B_1/B_1\gp ], $$
 is a bijection with the inverse map $[N]\mapsto \supp (N) := \{
 \gp \in \Max (K[H])\, | \, \gp \cdot n_\gp =0 $ for some $0\neq n_\gp
 \in N\}$.
 \item {\rm (Classification of simple $K[H]$-torsion free $
\mI_1$-modules)}
\begin{enumerate}
\item The map $$\hB1 (K[H]-{\rm torsion \; free})\ra \hCB1, \; \;
[M]\mapsto [ S^{-1} M], $$ is a bijection with the inverse map
$\soc_{B_1}:[ N]\mapsto [\soc_{B_1}(N)]$ where $ \soc_{B_1}(N)$ is
the socle of the $\CB_1$-module  $N$. \item Let $b=
\der^{-m}\beta_{-m}+\cdots +\beta_0\in B_1$  be a normal element
of the algebra $B_1$ which is  an irreducible element of the algebra $\CB_1$ where all
$\beta_i\in K[H]$, $\beta_{-m}\neq 0$, $\beta_0\neq 0$,  and
$m>0$. Then $B_1/ B_1\cap \CB_1 b$ is a $K[H]$-torsion free simple
$B_1$-module. Up to isomorphism, every $K[H]$-torsion free simple
$B_1$-module arises in this way, and from $b$ which is unique up
to similarity. \item An $\mI_1$-module $M$ is simple
$K[H]$-torsion free iff $M\simeq (B_1\alpha +\CB_1b)/ \CB_1b\simeq
B_1/ B_1\cap \CB_1b\alpha^{-1}$ for some irreducible in $\CB_1$
element $b\in B_1$ such that the element $\beta b\alpha^{-1} \in
B_1$ is normal where  $\alpha , \beta \in K[H]$.
\end{enumerate}
\end{enumerate}
\end{theorem}

{\it Proof}. 1.Statement 1 follows from the facts that $K[x]$ is the
only (up to isomorphism) faithful simple $\mI_1$-module
(Proposition 6.1, \cite{algintdif}), $F$ is the only proper ideal
of the algebra $\mI_1$ and $B_1= \mI_1/F$.

2. Statement 2 is obvious.

3. Statement 3 is  a particular case of Theorem 11,
\cite{Bav-GWA-Ker-TensoSimple-92}.

4. Statement 4 is  a particular case of Theorem 17,
\cite{Bav-GWA-Ker-TensoSimple-92}. $\Box $

$\noindent $

Let $\gp \in \Max (K[H])$. Then the simple $\mI_1$-module $B_1/
B_1\gp \simeq \bigoplus_{i\in \Z}\der^i K[H]/\gp$ is a semi-simple
$K[H]$-module since ${}_{K[H]}(\der^i K[H]/\gp)\simeq
K[H]/\tau^i(\gp)$.

A classification of simple $A_1$-modules was given by Block,
\cite{Bl 1} and \cite{Bl 3}.

$\noindent $

{\bf Relationships between the sets of  simple $\mI_1$-modules and $A_1$-modules}. We will
see that the algebras $\mI_1$ has ``one less'' simple module than
the Weyl algebra $A_1$ (Theorem \ref{7Jun10}.(3)). The algebra
$B_1$ is a left (but not right) localization of the algebra
$\mI_1$ at the powers of the element $\der$ \cite{algintdif},
$B_1=S_\der^{-1}\mI_1$ where $S_\der := \{ \der^i\, | \, i\in
\N\}$. The algebra $B_1=S_\der^{-1}A_1$ is the  left and right
localization of the Weyl algebra $A_1$ at $S_\der$. It follows at
once from the fact that the $\mI_1$-module $\mI_1/ \mI_1\der=
K[x]$ is  simple and $\der$-torsion (i.e.,  $S_\der^{-1} K[x]=0$)
that
$$ \hmI1 (\der-{\rm torsion}) = \{ [K[x]] \}.$$
Similarly, the $A_1$-module $A_1/A_1\der = K[x]$ is simple and
$\der$-torsion. Then
$$ \hA_1 (\der-{\rm torsion}) = \{ [K[x]] \}.$$
 The algebra $\CB_1= K(H)[\der, \der^{-1}; \tau ]$ is a left and
 right localization of the algebras $\mI_1$ and $A_1$ at $S:=
 K[H]\backslash \{ 0\}$,
 $$\CB_1= S^{-1}\mI_1= S^{-1}A_1.$$
 Therefore, the sets of simple $\mI_1$- and $A_1$-modules can be
 represented as the disjoint unions of $K[H]$-torsion and
 $K[H]$-torsion free modules:
\begin{eqnarray*}
 \hmI1 &=& \hmI1 (K[H]-{\rm torsion})\coprod  \hmI1 (K[H]-{\rm torsion \;
free}),\\
\hA_1 &=& \hA_1 (K[H]-{\rm torsion})\coprod  \hA_1 (K[H]-{\rm
torsion \; free}),\\
\end{eqnarray*}
where
\begin{eqnarray*}
  \hmI1 (K[H]-{\rm torsion})&=& \{ [K[x]]\} \coprod  \hB1 (K[H]-{\rm
  torsion}),\\
  \hmI1  (K[H]-{\rm torsion \; free})&=& \hB1 (K[H]-{\rm torsion \; free}),\\
\hA_1 (K[H]-{\rm torsion})&=& \{ [A_1/A_1x]\} \coprod  \{ [K[x]]\}
\coprod  \hA_1' (K[H]-{\rm  torsion}),
\end{eqnarray*}
where the set $\hA_1' (K[H]-{\rm  torsion})$ is obtained from the
set $\hA_1 (K[H]-{\rm torsion})$ by deleting its two elements
$[A_1/A_1x]$ and $[K[x]]$.

 Recall that the Weyl algebra $A_1$ is a GWA
$K[H](\s , H)$, $\s (H) = H-1$, over the Dedekind domain $K[H]$.
There is a similar to Theorem \ref{31May10} classification of
simple $A_1$-modules \cite{Bav-AlgAnaliz92}, \cite{Bav-UkrMathJ-92},
\cite{Bav-GWA-Ker-TensoSimple-92}. Comparing it to Theorem
\ref{31May10}, the next theorem follows.

\begin{theorem}\label{7Jun10}
\begin{enumerate}
\item The map $\hA_1' (K[H]-{\rm  torsion})\ra \hB1 (K[H]-{\rm
  torsion})$, $[M]\mapsto [ S_\der^{-1}M]$, is a bijection with the
  inverse map $[N]\mapsto [{}_{A_1}N]$.
\item The map $\hA_1 (K[H]-{\rm  torsion\; free})\ra \hB1
(K[H]-{\rm
  torsion \; free})$, $[M]\mapsto [ S_\der^{-1}M]$, is a bijection with the
  inverse map $[N]\mapsto [{\rm soc}_{A_1}N]$.
\item Combining statements 1 and 2, the map
$$ \hA_1\backslash \{ [ A_1/ A_1x]\} \ra \hmI1 = \{ [K[x]]\}
\coprod  \hB1, \;\;  [K[x]]\mapsto  [K[x]], \;\; [M]\mapsto [
S_\der^{-1}M],$$ is a bijection with the inverse map
$[K[x]]\mapsto [K[x]]$ and $[N]\mapsto [{\soc}_{A_1}\, N]$ where
$[N]\in \hB1$.
\end{enumerate}
\end{theorem}


\section{The Strong Compact-Fredholm Alternative}\label{FDIMKERCO}

In this section, the Strong Compact-Fredholm Alternative is proved
for the algebra $\mI_1$ (Theorem \ref{B30May10}), it is shown that
the endomorphism algebra of each simple $\mI_1$-module is a finite
dimensional  division algebra (Theorem \ref{6Jun10}). Various
indices are introduced for non-compact integro-differential
operators, the $M$-index and the length index, and we show  that
they are invariant under addition of  compact operator (Lemma
\ref{D30May10} and Lemma \ref{C30May10}).

$\noindent $

{\bf The (Strong) Compact-Fredholm Alternative}. Let $\CF = \CF
[K]$ be the family of all $K$-linear maps with finite dimensional
kernel and cokernel, such maps are called the {\em Fredholm linear
maps/operatorts}. So, $\CF$ is the family of {\em Fredholm} linear
maps/operators. For vector spaces $V$ and $U$, let $\CF (V,U)$ be
the set of all the linear maps from $V$ to $U$ with finite
dimensional kernel and cokernel and $\CF (V) :=\CF (V,V)$. So,
$\CF =\coprod_{V,U}\CF (V,U)$ is the disjoint union.

$\noindent $

{\it Definition}. For a linear map $\v \in \CF$, the integer
$$ \ind (\v ) := \dim \, \ker (\v ) - \dim \, \coker (\v )$$
is called the {\em index} of the map $\v$.

$\noindent $

{\it Example}. Note that $\der, \int \in \mI_1\subset
\End_K(K[x])$. Then  
\begin{equation}\label{indId}
\ind (\der^i)= i\;\; {\rm and }\;\; \ind (\int^i)= -i, \;\; i\geq
1.
\end{equation}

It is well-known that for any two linear maps $M\stackrel{a}{\ra}
N \stackrel{b}{\ra}L$ there is the long  exact sequence of vector
spaces 
\begin{equation}\label{ables}
0\ra\ker (a)\ra \ker (ba)\ra\ker (b) \ra \coker (a)\ra \coker
(ba)\ra\coker (b)\ra 0.
\end{equation}

Let $\CC$ be the family of $K$-linear maps with finite dimensional
image, such maps are called the {\em compact} linear
maps/operators. For vector spaces $V$ and $U$, let $\CC (V,U)$ be
the set of all the compact linear maps from $V$ to $U$. So, $\CC
=\coprod_{V,U}\CC (V,U)$ is the disjoint union. If $V=U$ we write
$\CC (V) := \CC (V,V)$.

$\noindent $

{\it Definitions}. Let $A$ be an algebra and $M$ be its module. We
say that, for the $A$-module $M$, the {\em Compact-Fredholm
Alternative} holds if, for each element $a\in A$,   the linear map
$a_M:M\ra M$, $m\mapsto am$,   is either
 compact  or Fredholm. We say that for the algebra $A$ the {\em Compact-Fredholm
Alternative} holds if for each {\em simple} $A$-module it does;
and the {\em Strong Compact-Fredholm Alternative} holds if, in
addition, for each element $a\in A$ either the linear maps $a_M$
are compact for all simple $A$-modules $M$ or, otherwise,  they
are all Fredholm. In the first and the second case the element $a$
is called a {\em compact} and {\em Fredholm} element respectively,
and the sets of all compact and Fredholm elements of the algebra
$A$ are denoted by $\CC_A$ and $\CF_A$ respectively; then
$$ A=\CC_A\coprod \CF_A$$
is the disjoint union (for the algebra $A$ with Strong
Compact-Fredholm Alternative). Clearly, $\CC_A$ is an ideal of the
algebra $A$ and $\CC_A \supseteq \rad (A)$, the Jacobson radical
of $A$. By (\ref{ables}), $\CF_A$ is a multiplicative monoid that
contains the group of units of the algebra $A$.

$\noindent $

The next theorem shows that each non-compact integro-differential
operator has only finitely many linearly independent solutions in
each (left or right) $\mI_1$-module of finite length. An analogous
result is known for the Weyl algebra $A_1$ \cite{McCon-Rob-JA-73}
(for simple $A_1$-modules). Notice that each nonzero element of
the Weyl algebra $A_1$ is a non-compact operators, i.e.,  $A_1\cap
F=0$.

\begin{theorem}\label{B30May10}
{\rm (Strong Compact-Fredholm Alternative)} For the algebra
$\mI_1$ the Strong Compact-Fredholm Alternative holds (for left
and right simple $\mI_1$-modules); $\CC_{\mI_1}=F$ and
$\CF_{\mI_1}= \mI_1\backslash F$. Moreover,
\begin{enumerate}
\item Let $a\in \mI_1$, ${}_{\mI_1}M$ be a nonzero $\mI_1$-module
of finite length and $a_M:M\ra M$, $m\mapsto am$. Then
\begin{enumerate}
\item $\dim_K(\ker (a_M))<\infty$ iff $\dim_K(\coker
(a_M))<\infty$ iff $a\not\in F$. \item If, in addition,
${}_{\mI_1}M$ is simple then $\dim_K(\im (a_M))<\infty$ iff $a\in
F$.
\end{enumerate}
\item Let $a\in \mI_1$, $N_{\mI_1}$ be a nonzero right
$\mI_1$-module of finite length and $a_N:N\ra N$,  $n\mapsto na$.
Then
\begin{enumerate}
\item  $\dim_K(\ker (a_N))<\infty$ iff $\dim_K(\coker
(a_N))<\infty$ iff $a\not\in F$. \item If, in addition,
$N_{\mI_1}$ is simple then $\dim_K(\im (a_N))<\infty$ iff $a\in
F$.
\end{enumerate}
\end{enumerate}
\end{theorem}

{\it Proof}. 1(a). Using induction on the length of the module $M$
and the Snake Lemma we see that statement 1(a) holds iff it does
for each simple $\mI_1$-module. If $a\in F$ then, by Theorem
\ref{31May10}, for each simple $\mI_1$-module $M$, the kernel and
the cokernel of the map $a_M$ are infinite dimensional spaces. It
remains to show that if $a\not\in F$ then the kernel and the
cokernel of the map $a_M$ are finite dimensional spaces. First
we consider the case when $M=K[x]$.  The algebra $\mI_1$ is a
$\Z$-graded algebra. The usual $\N$-grading of the polynomial
algebra $K[x]=\bigoplus_{i\geq 0}Kx_i$ is also a $\Z$-grading of
the $\mI_1$-module $K[x]$. It determines the filtration
$K[x]=\bigcup_{i\geq 0} K[x]_{\leq i}$ where $K[x]_{\leq
i}:=\bigoplus_{j=0}^iKx^j$. Clearly, when $a=\der^i$ or
$a=\int^i$, the map $a_{K[x]}$  has  finite dimensional  kernel and
cokernel. Using (\ref{ables}) and (\ref{acan}), we may assume that
all the polynomials $a_i$, $i>0$, in the decomposition
(\ref{acan}) are equal to zero and $a_0\neq 0$. Fix a natural
number $m$ such that
 all $\l_{ij}=0$ for all $i,j\geq m$  and
$a_0(H)*x^l=a(l+1)x^l\neq 0$  for all $l\geq m$. Then
$$a K[x]_{\leq l}\subseteq K[x]_{\leq l}\;\; {\rm   for\;  all}\;\; l\geq m,$$ and the
element $a$ acts as a bijection on the factor space $K[x]/
K[x]_{\leq m}$. This implies that the map $a_{K[x]}$ has finite
dimensional kernel and cokernel.

Suppose that $M\not\simeq K[x]$, i.e.,  $M$ is a simple $B_1$-module
(Theorem \ref{31May10}.(1)). The algebra $B_1$ is an example of
the generalized Weyl algebra (GWA) $A= K[H](\tau , u)$ over the
Dedekind domain $K[H]$ where $\tau (H) = H+1$ and $u=1$. Statement
2 was established for the GWA $A$ over an {\em algebraically
closed} field and for simple $A$-modules (Theorem 4,
\cite{Bav-AlgAnaliz92}). Using the classification of simple
$\mI_1$-modules (Theorem \ref{31May10}) and the classification of
simple $A$-modules in \cite{Bav-AlgAnaliz92} (or
\cite{Bav-GWA-Ker-TensoSimple-92}) we see that, for each simple
$B_1$-module $M$, the $B_1(\bK )$-module $B_1(\bK
)\t_{B_1}M=\bK\t_KM$ has finite length (where $\bK$ is the
algebraic closure of the field $K$), and the result follows from
Theorem 4, \cite{Bav-AlgAnaliz92}.

There is a more elementary way to show that the $B_1(\bK )$-module
$\bK \t_K M$ has finite length. For, notice that over an arbitrary
field (i) the algebra $B_1$ is a simple Noetherian domain of
Gelfand-Kirillov dimension 2 which is an almost commutative
algebra with respect to the standard filtration $\CV$ associated
 with the generators $\der$ , $\der^{-1}$ and $H$ (i.e.,  the
associated graded algebra is a commutative finitely generated
algebra);  (ii) the Gelfand-Kirillov dimension of each simple
$B_1$-module is 1; (iii) every finitely generated $B_1$-module of
Gelfand-Kirillov dimension 1 has finite length which  does not
exceed the multiplicity of the module (with respect to the
filtration $\CV$); (iv) since the Gelfand-Kirillov dimension and
the multiplicity are invariant under  field extensions, by (ii)
and (iii), the $B_1(\bK )$-module $\bK \t_KM$ has finite length.

1(b). $(\Rightarrow )$ Each simple $\mI_1$-module $M$ is infinite
dimensional. So, if $\dim_K(\im (a_M))<\infty$ then necessarily
$a\in F$, by statement 1(a).

$(\Leftarrow )$ If $a\in F$ then, clearly,  $\dim_K(\im
(a_M))<\infty$ for all simple $\mI_1$-modules $M$: this is obvious
when $M\simeq {}_{\mI_1}K[x]$ but if $M\not\simeq {}_{\mI_1}K[x]$
then $a_M=0$ (Theorem \ref{31May10}.(1)).

2. In view of the involution $*$ on the algebra $\mI_1$, statement
2 follows from statement 1. $\Box $


\begin{corollary}\label{9Sep10}
$(\CC_{\mI_n})^*=\CC_{\mI_n}$ and $ (\CF_{\mI_1})^*= \CF_{\mI_1}$,
i.e.,  the involution $*$ preserves the compact and Fredholm
elements of the algebra $\mI_1$.
\end{corollary}

{\it Proof}. This follows from $F^*=F$, $\CC_{\mI_1} = F$ and
$\CF_{\mI_1} = \mI_1\backslash F$.  $\Box $

\begin{corollary}\label{aB30May10}
$\mI_n\cap \CC (K[x])=\CC_{\mI_1}=F$.
\end{corollary}

{\it Proof}. By the very definition, $\mI_1\subseteq
\End_K(K[x])$. Since $F\subseteq \CC (K[x])$ and $(\mI_1\backslash
F)\bigcap  \CC (K[x])=\emptyset$ (by Theorem \ref{B30May10}), the first
equality follows. The second equality is obvious (Theorem \ref{B30May10}).  $\Box $

$\noindent $

{\bf The endomorphism algebras of simple $\mI_1$-modules are finite dimensional}. The endomorphism algebra $\End_A(M)$ of a simple module over an algebra $A$ is a division algebra.

\begin{theorem}\label{6Jun10}
\begin{enumerate}
\item $\dim_K(\End_{\mI_1}(M))<\infty$ for all simple
$\mI_1$-modules $M$. \item  $\End_{\mI_1}(K[x])\simeq K$.\item
$\End_{\mI_1}(B_1/B_1\gp ) \simeq K[H]/\gp$ for all $\gp \in \Max
(K[H])$.
\end{enumerate}
\end{theorem}

{\it Proof}. 2. $\End_{\mI_1}(K[x])\simeq  \End_{\mI_1}(\mI_1/
\mI_1\der )\simeq \ker_{K[x]}(\der ) =K$.

3. Recall that the $K[H]$-module $B_1/ B_1\gp = \bigoplus_{i\in
\Z} \der^i K[H]/ \gp$ is the direct sum of simple pairwise
non-isomorphic $K[H]$-modules ${}_{K[H]}( \der^i K[H]/ \gp )\simeq
K[H]/ \tau^i (\gp )$. Therefore, $\End_{\mI_1}(B_1/B_1\gp ) \simeq
\ann_{B_1/B_1\gp }(\gp ) =  K[H]/\gp$.

1. By Theorem  \ref{31May10} and statements  2 and 3, it remains
to show that the endomorphism algebra of the $\mI_1$-module $M$ is
finite dimensional  in the case when $M$ is $K[H]$-torsion free,
i.e.,   $M\simeq B_1/ B_1\cap \CB_1 b$ where the element $b\in B_1$
is as in Theorem \ref{31May10}.(4).  Since $\End_{\mI_1}(M) \simeq
\ann_M(B_1\cap \CB_1b) \subseteq \ker_M(b)$ and $b\not\in F$, the
algebra $\End_{\mI_1}(M)$ is finite dimensional by Theorem
\ref{B30May10}.(1a).  $\Box $

$\noindent $

{\bf The $M$-index and its properties}.

$\noindent $

{\it Definition}. For each element $a\in \mI_1\backslash F$ and a
(left or right) $\mI_1$-module $M$ of finite length, we define the
$M$-{\em index} of the element $a$,
$$ \ind_M(a) :=\dim_K(\ker (a_M))-\dim_K(\coker (a_M)).$$
The set $ \mI_1\backslash F$ is a multiplicative monoid. For all
elements $a,b\in  \mI_1\backslash F$, 
\begin{equation}\label{Mind}
\ind_M(ab) = \ind_M(a) +\ind_M(b).
\end{equation}
This follows from (\ref{ables}). Therefore, the $M$-index $\ind_M
:\mI_1\backslash F \ra \Z$ is a monoid homomorphism. The next lemma shows that the $M$-index is invariant under addition of compact operator.

\begin{lemma}\label{D30May10}
Let $a\in \mI_1\backslash F$ and $f\in F$. Then
$\ind_M(a+f)=\ind_M(a)$ for all left or right  $\mI_1$-modules $M$
of finite length. i.e.,  $\CF_{\mI_1}+\CC_{\mI_1}\subseteq
\CF_{\mI_1}$.
\end{lemma}

{\it Proof}. Note that $\der^if=0$ for all $i\gg 0$, and $\der^i,
\der^ia\not\in F$. Then, by (\ref{Mind}),
$$ \ind_M(\der^i) +\ind_M(a) = \ind_M(\der^i a) = \ind_M(\der^i
(a+f))= \ind_M(\der^i) + \ind_M(a+f),$$ hence
$\ind_M(a)=\ind_M(a+f)$.  $\Box $

$\noindent $

Let $\s$ be an automorphism of an algebra $A$ and $M$ be an
$A$-module. The {\em twisted} $A$-module ${}^\s M$  by the
automorphism $\s$  coincides with $M$ as a vector space but the
$A$-module structure on ${}^\s M$ is given by the rule $a\cdot m=
\s (a) m$ for all elements $a\in A$ and $m\in M$. For all elements
$a\in \mI_1\backslash F$, an automorphism $\s \in \Aut_{K-{\rm
alg}}(\mI_1)$ and an $\mI_1$-module $M$ of finite length (Theorem
\ref{B30May10}), 
\begin{equation}\label{indMsa}
\ind_M(\s (a)) = \ind_{{}^\s M}(a).
\end{equation}
If, in addition, ${}^\s_{\mI_1}M\simeq M$ for some automorphism
$\s$, then
\begin{equation}\label{1indMsa}
\ind_M(\s (a)) = \ind_M(a)\;\; {\rm for \; all}\;\; a\in
\mI_1\backslash F.
\end{equation}
If ${}^\s_{\mI_1}M\not\simeq M$ then, in general, the equality
(\ref{1indMsa}) {\em does not hold and this observation is a very
effective tool in proving that two modules are non-isomorphic}.

$\noindent $

{\em Example}. For each $\l \in K^*$, the $\mI_1$-module $V_\l :=
B_1/B_1(\der -\l )$ is simple: the $K[H]$-module homomorphism
$K[H]\ra V_\l$, $\alpha \mapsto \alpha +B_1 (\der - \l )$, is
obviously an epimorphism; therefore it is a $K[H]$-module
isomorphism as each nonzero $B_1$-module is infinite dimensional
and each proper factor module of the $K[H]$-module $K[H]$ is
finite dimensional; then $V_\l$ is a simple $B_1$-module such that
${}_{K[H]}V_\l \simeq K[H]$. For each $\l \in K^*$, there is an
automorphism $t_\l $ of the algebra $\mI_1$ given by the rule
$$t_\l (\int ) =\l \int, \;\;\; t_\l (\der ) = \l^{-1} \der, \;\;\; t_\l
(H)=H.$$ Since $t_\l (e_{ij})=\l^{i-j}e_{ij}$, $t_\l (F) = F$,
 $t_\l $ induces the automorphism $t_\l$ of the factor
algebra $B_1= \mI_1 / F$  by the rule $t_\l (\der ) = \l^{-1}
\der$ and $t_\l (H) = H$. Since $t_\l t_\mu = t_{\l \mu }$ for all
$\l , \mu \in K^*$, the {\em algebraic torus} $\mT^1:=\{ t_\l \, |
\, \l \in K^*\}\simeq K^*$ is a subgroup of the groups
$\Aut_{K-{\rm alg}}(\mI_1)$ and $\Aut_{K-{\rm alg}}(B_1)$.
Clearly, $V_\l \simeq {}^{t_{\l^{-1}}}V_1$ for all $\l \in K^*$
(since $t_{\l^{-1}}(\der -\l ) = \l (\der -1)$). The simple
$\mI_1$-modules $\{ V_\l \, | \, \l \in K^*\}$ is the
$\mT^1$-orbit of the module $V_1$. We aim to show that
\begin{equation}\label{VlVm}
{}_{\mI_1}V_\l \simeq {}_{\mI_1}V_\mu \;\; {\rm iff}\;\; \l = \mu.
\end{equation}
Without loss of generality we may assume that $\mu =1$. Then the
result follows at once from the fact that 
\begin{equation}\label{1VlVm}
\ind_{V_1}(\der -\l ) =\begin{cases}
1& \text{if }\l =1,\\
0& \text{otherwise}.\\
\end{cases}
\end{equation}
{\it Proof of (\ref{1VlVm})}.  When we identify the $\mI_1$-module $V_1$ with $K[H]$ (as we did
above),  the action of the element $\der -\l$ on $V_1$ is
identified with the action of the linear map $\tau - \l$ on $K[H]$
where $\tau (H) = H+1$ since $\der \cdot 1 = 1$. The  map $\tau
-1$ is surjective with kernel $K$, then $\dim_K(\ker_{V_1}(\der
-1)) =1$. The map $\tau - \l$ (where $\l \neq 1$) is an
isomorphism, therefore $\ind_{V_1}(\der -\l ) =0$. $\Box$

\begin{theorem}\label{A30May10}
Let $a\in \mI_1$, $\cdot a:\mI_1\ra \mI_1$, $b\mapsto ba$, $a\cdot
: \mI_1\ra \mI_1$, $b\mapsto ab$, and  $l_{\mI_1}$ be the
length function on the set of (left or right) $\mI_1$-modules. Then
\begin{enumerate}
\item
\begin{enumerate}
\item $l_{\mI_1}(\ker (\cdot a))<\infty $ iff $l_{\mI_1}(\coker
(\cdot a))<\infty$ iff $a\not\in F$. \item $l_{\mI_1}(\im (\cdot
a))<\infty $ iff iff $a\in F$. 
\end{enumerate}
\item
\begin{enumerate}
 \item  $l_{\mI_1}(\ker (
a\cdot))<\infty $ iff $l_{\mI_1}(\coker ( a\cdot))<\infty$ iff
$a\not\in F$. \item $l_{\mI_1}(\im (a\cdot ))<\infty $ iff iff
$a\in F$. 
\end{enumerate}
\end{enumerate}
\end{theorem}

{\it Proof}. 1(a). Applying the Snake Lemma to the commutative
diagram of $\mI_1$-modules
$$
\xymatrix{0\ar[r] & F\ar[r]\ar[d]^{\cdot a}  & \mI_1 \ar[r]\ar[d]^{\cdot a} & B_1 \ar[r]\ar[d]^{\cdot \overline{a}} & 0 \\
0\ar[r] & F\ar[r]  & \mI_1 \ar[r] & B_1 \ar[r] & 0 }
$$
(where $\oa:=a+F\in B_1$) yields the long exact sequence of
$\mI_1$-modules 
\begin{equation}\label{les1}
0\ra\ker_F (\cdot a)\ra \ker_{\mI_1} (\cdot a)\ra\ker_{B_1} (\cdot
\oa) \ra \coker_F (\cdot a)\ra \coker_{\mI_1} (\cdot
a)\ra\coker_{B_1} (\cdot \oa)\ra 0.
\end{equation}
If $a\in F$ then $l_{\mI_1}(\ker_F(\cdot a))=\infty$ and
$l_{\mI_1}(\coker_{B_1}(\cdot \oa))=l_{\mI_1}(B_1)=\infty$ since
$\oa =0$. Therefore, $l_{\mI_1}(\ker_{\mI_1}(\cdot a))=\infty$ and
 $l_{\mI_1}(\coker_{\mI_1}(\cdot a))=\infty$, by (\ref{les1}). To
 finish the proof it suffices  to show that
 $l_{\mI_1}(\ker_{\mI_1}(\cdot a))<\infty$ and
 $l_{\mI_1}(\coker_{\mI_1}(\cdot a))<\infty$ provided $a\not\in
 F$. By (\ref{les1}), it suffices to show that the $\mI_1$-modules
 $\ker_F(\cdot a)$,  $\coker_F(\cdot a)$,  $\ker_{B_1}(\cdot \oa)$ and
   $\coker_{B_1}(\cdot \oa)$ have finite length. Since $\oa \neq 0$,
   $\ker_{B_1}(\cdot \oa ) =0$ and $l_{\mI_1}(\coker_{B_1}(\cdot
   \oa )) = l_{B_1}(\coker_{B_1}(\cdot
   \oa ))<\infty$ since the algebra $B_1$ is a localization of the
   Weyl algebra $A_1$ for which the analogous property holds (i.e.,
   $l_{A_1}(A_1/ A_1u)<\infty$ for all nonzero elements $u\in
   A_1$). Notice that
   $$F= \mI_1E_{00}\mI_1= \mI_1E_{00}E_{00}\mI_1=E_{\N , 0}E_{0,
   \N}=E_{\N , 0}\t E_{0, \N}$$
   where ${}_{\mI_1}E_{\N , 0}= \mI_1E_{00} = \bigoplus_{i\in
   \N} KE_{i0} \simeq {}_{\mI_1}K[x]$, $E_{00}\mapsto 1$, and
 $(E_{0, \N })_{\mI_1}= E_{00}\mI_1 = \bigoplus_{i\in
   \N}KE_{0i} \simeq K[\der]_{\mI_1}$, $E_{00}\mapsto 1$ where $K[\der ]_{\mI_1}\simeq \mI_1/ \int \mI_1$ is a simple right $\mI_1$-module. By
   Theorem \ref{B30May10}.(2a), the linear map $a_{K[\der ]}: K[\der ] \ra K[\der ]$, $ v\mapsto av$,  has
   finite dimensional kernel and cokernel since the module
    $K[\der ]_{\mI_1}$ is simple. Since
\begin{equation}\label{kcoE}
\ker_F(\cdot a) = E_{\N , 0}\t \ker (a_{K[\der ]}), \;\;
\coker_F(\cdot a) \simeq E_{\N , 0}\t \coker (a_{K[\der ]}),
\end{equation}
and the $\mI_1$-module $E_{\N , 0}\simeq K[x]$ is simple, we see
that  
\begin{equation}\label{kcoE1}
l_{\mI_1}(\ker_F(\cdot a)) = \dim_K(\ker (a_{K[\der ]}))<\infty ,
\;\; l_{\mI_1}(\coker_F(\cdot a)) =\dim_K(\coker (a_{K[\der
]}))<\infty .
\end{equation}
The proof of statement 1(a) is complete.

1(b). $(\Rightarrow )$ Note that $l({}_{\mI_1}\mI_1)=\infty$. So,
if $l_{\mI_1}(\im (\cdot a))<\infty$ then necessarily $a\in F$, by
statement 1(a).

$(\Leftarrow )$ If $a\in F$ then, clearly
$l{}_{\mI_1} (\im (\cdot a))<\infty$.


2. Statement 2 follows from statement 1 by applying the involution
$*$ of the algebra $\mI_1$ and using the equality $F^*=F$. $\Box $

$\noindent $

{\bf The left and right length indices}.

$\noindent $

{\it Definition}. For each element $a\in \mI_1\backslash F$, we
define its left and right {\em length index} respectively as
follows (by Theorem \ref{A30May10})
$$ \lind (a) = l_{\mI_1}(\ker (\cdot a))-l_{\mI_1}(\coker (\cdot a)), \;\;
\rind (a) = l_{\mI_1}(\ker ( a\cdot))-l_{\mI_1}(\coker ( a\cdot)).
$$

The set $ \mI_1\backslash F$ is a multiplicative monoid. For all
elements $a,b\in  \mI_1\backslash F$, 
\begin{equation}\label{xcoE2}
\lind (ab) = \lind (a) +\lind (b),\;\; \rind (ab) = \rind (a)
+\rind (b).
\end{equation}
These equalities  follow from (\ref{ables}). Therefore, the indices  $\lind ,
\rind  :\mI_1\backslash F \ra \Z$ are  monoid homomorphisms. The next lemma show that the left and right indices are invariant under addition of compact operator.

\begin{lemma}\label{C30May10}
Let $a\in \mI_1\backslash F$ and $f\in F$. Then $\lind (a+f)=\lind
(a)$  and $\rind (a+f)=\rind (a)$.
\end{lemma}

{\it Proof}. Note that $\der^if=0$ for all $i\gg 0$, and $\der^i,
\der^ia\not\in F$. Then, by (\ref{xcoE2}),
$$ \lind (\der^i) +\lind (a) = \lind (\der^i a) = \lind (\der^i
(a+f))= \lind (\der^i) + \lind (a+f),$$ hence $\lind (a)=\lind
(a+f)$. Replacing $\lind$ by $\rind$ in the argument above, the
second equality follows.  $\Box $



\section{The algebra $\mI_1$ is a coherent algebra}\label{COHERENT}

In this section, we prove that the algebra $\mI_1$ is a left and
right coherent algebra (Theorem \ref{a19Sep10}) and that every
finitely generated left (or right) ideal of the algebra $\mI_1$ is
generated by two elements (Theorem \ref{A19Sep10}).

Let $V$ be a vector space. A linear map $\v : V\ra V$ is called a
{\em locally nilpotent} map if $V = \bigcup_{i\geq 1}\ker (\v^i)$,
i.e.,  for each element $v\in V$, $\v^i v =0$ for some $i=i(v)$.

\begin{lemma}\label{a7Jun10}
Let $a,b\in \mI_1$. Let  $a\cdot$ and $\cdot b$ be  the left and
right multiplication maps in $\mI_1$ by the elements  $a$ and $b$
respectively.
\begin{enumerate}
\item There are short exact sequences of left $\mI_1$-modules:
\begin{enumerate}
\item  $0\ra \ker_F(\cdot b)\ra \ker_{\mI_1}(\cdot b)\ra
\ker_{B_1}(\cdot b) \ra 0$, \item $0\ra \coker_F(\cdot b)\ra
\coker_{\mI_1}(\cdot b)\ra \coker_{B_1}(\cdot b) \ra 0$.
\end{enumerate}
\item There are long exact sequences of vector spaces:
\begin{enumerate}
\item  $0\ra \ker_{\ker_F(\cdot b)}(a\cdot )\ra
\ker_{\ker_{\mI_1}(\cdot b)}(a\cdot )\ra \ker_{\ker_{B_1}(\cdot
b)}(a\cdot ) \ra \coker_{\ker_F(\cdot b)}(a\cdot )\ra
$

$\coker_{\ker_{\mI_1}(\cdot b)}(a\cdot )\ra
\coker_{\ker_{B_1}(\cdot b)}(a\cdot ) \ra 0$,

\item  $0\ra \ker_{\coker_F(\cdot b)}(a\cdot )\ra
\ker_{\coker_{\mI_1}(\cdot b)}(a\cdot )\ra
\ker_{\coker_{B_1}(\cdot b)}(a\cdot ) \ra \coker_{\coker_F(\cdot
b)}(a\cdot )\ra \coker_{\coker_{\mI_1}(\cdot b)}(a\cdot )\ra
\coker_{\coker_{B_1}(\cdot b)}(a\cdot ) \ra 0$.
\end{enumerate}
\end{enumerate}
\end{lemma}

{\it Proof}. 1. By (\ref{les1}), we have the long exact sequence
of  left $\mI_1$-modules:
$$ 0\ra \ker_F(\cdot b)\ra \ker_{\mI_1}(\cdot b)\ra \ker_{B_1}(\cdot
b) \stackrel{\d}{\ra} \coker_F(\cdot b)\ra \coker_{\mI_1}(\cdot
b)\ra \coker_{B_1}(\cdot b) \ra 0.$$ The map $\d$ is equal to zero
since the element $\der$ acts as an invertible linear map on
$\ker_{B_1}(\cdot b)$ but its action on $\coker_F(\cdot b)$ is a
locally nilpotent map. So, the long exact sequence breaks down
into two short exact sequences (as above).

2. Applying the Snake Lemma to the commutative diagrams  of short
exact sequences of vector spaces
$$
\xymatrix{0\ar[r] & \ker_F(\cdot b)\ar[r]\ar[d]^{a\cdot }  & \ker_{\mI_1}(\cdot b) \ar[r]\ar[d]^{a\cdot } & \ker_{B_1}(\cdot b)  \ar[r]\ar[d]^{a\cdot } & 0 \\
0\ar[r] & \ker_F(\cdot b) \ar[r]  & \ker_{\mI_1}(\cdot b)  \ar[r]
&  \ker_{B_1}(\cdot b) \ar[r] & 0, }$$
$$
\xymatrix{0\ar[r] & \coker_F(\cdot b)\ar[r]\ar[d]^{a\cdot }  & \coker_{\mI_1}(\cdot b) \ar[r]\ar[d]^{a\cdot } & \coker_{B_1}(\cdot b)  \ar[r]\ar[d]^{a\cdot } & 0 \\
0\ar[r] & \coker_F(\cdot b) \ar[r]  & \coker_{\mI_1}(\cdot b)
\ar[r] &  \coker_{B_1}(\cdot b) \ar[r] & 0. }
$$
yields the long exact sequences of statement 2. $\Box $

\begin{theorem}\label{21Sep10}
Let $a\in \mI_1$. Then
\begin{enumerate}
\item $\ker_{\mI_1}(\cdot a)$ and $\coker_{\mI_1}(\cdot a)$ are
finitely generated left $\mI_1$-modules. \item $\ker_{\mI_1}(
a\cdot )$ and $\coker_{\mI_1}( a\cdot )$ are finitely generated
right $\mI_1$-modules.
\end{enumerate}
\end{theorem}

{\it Proof}.  An algebra $\mI_1$ is self-dual, so it suffices to
prove only statement 1. If $a\not\in F$ then statement 1 follows
from Theorem \ref{A30May10}.(1a). We may assume that $a\in F$ and
$a\neq 0$. By Lemma \ref{a7Jun10}.(1), there are short exact
sequences of left $\mI_1$-modules
$$ 0\ra \ker_F(\cdot a) \ra \ker_{\mI_1}( \cdot a) \ra B_1\ra 0
\;\; {\rm and}\;\;  0\ra \coker_F(\cdot a) \ra \coker_{\mI_1}( \cdot a) \ra
B_1\ra 0.$$ Fix elements $u\in \ker_{\mI_1}(\cdot a) $ and $v\in
\coker_{\mI_1}(\cdot a)$ that are mapped to $1\in B_1$ by the
second maps in the short exact sequences. Then $ u = 1+f$ for some
element $f\in F$ and $v= 1+g+\mI_1a$ for some element
$g\in F$. Since
$\mI_1u\subseteq \ker_{\mI_1}(\cdot a)$ and $l_{\mI_1}(\mI_1/
\mI_1u) = l_{\mI_1}(\coker_{\mI_1}(\cdot u)) <\infty$ (Theorem
\ref{A30May10}.(2a) as $u\not\in F$), we see that
$l_{\mI_1}(\ker_{\mI_1}(\cdot a)/ \mI_1 u )<\infty$, this means
that the left $\mI_1$-module $\ker_{\mI_1}(\cdot a)$ is finitely
generated. Similarly, since $\mI_1v\subseteq \coker_{\mI_1}(\cdot
a)$ and $l_{\mI_1}(\mI_1/ \mI_1(1+g))=
l_{\mI_1}(\coker_{\mI_1}(\cdot (1+g)))<\infty $ (Theorem
\ref{A30May10}.(2a) as $1+g\not\in F$),  we see that
$$ l_{\mI_1}(\coker_{\mI_1} (\cdot a)/\mI_1 v) = l_{\mI_1} (\mI_1
/ (\mI_1(1+g)+\mI_1a))\leq l_{\mI_1}(\mI_1/ \mI_1(1+g))<\infty.$$
This means that the left $\mI_1$-module $\coker_{\mI_1}(\cdot a)$
is finitely generated. $\Box $


\begin{theorem}\label{19Sep10}
The  intersection of finitely many  finitely generated left
(resp. right) ideals of the algebra $\mI_1$ is again  a finitely
generated left (resp. right) ideal of $\mI_1$.
\end{theorem}

{\it Proof}. Since the algebra $\mI_1$ is self-dual it suffices to
prove the statement for, say, left ideals, and only for two of
them. Let $I$ and $J$ be finitely generated left ideals of the
algebra $\mI_1$. If one of them, say $I$, belongs to the ideal $F$
then necessarily the left ideal $I$ is a finitely generated
semi-simple left $\mI_1$-module, hence so is the intersection
$I\cap J$. In particular, $I\cap J$ is  finitely generated.

We assume that neither $I$ nor $J$ belongs to $F$. Then their
images $\bI$ and $\bJ$ under the ring epimorphism $\mI_1\ra
\mI_1/F=B_1$ are nonzero left ideals of the ring $B_1= A_{1,x}$
which is the localization of the Weyl algebra $A_1$ at the powers
of the element $x$. Then $\bI \cap \bJ\neq 0$ (since $\bI\cap A_1
\neq 0$,  $\bJ \cap A_1\neq 0$, and the intersection of two
nonzero left ideals in $A_1$ is a nonzero left ideal). Take an
element $a\in \mI_1$ such that $0\neq a+F\in \bI \cap \bJ$. Then
$a\not\in F$ and $\mI_1 a\subseteq I\cap J$. Since
$l_{\mI_1}(\mI_1 / \mI_1 a)<\infty$ (by Theorem
\ref{A30May10}.1(a)), the ideal $I\cap J\neq 0$ is finitely
generated. $\Box $

$\noindent $

A finitely generated module is a {\em coherent} module if every
finitely generated submodule is finitely presented. A ring $R$  is
a {\em left} (resp. {\em right}) {\em coherent ring} if the module
${}_RR$ (resp. $R_R$) is coherent. {\em A ring $R$ is a left
coherent ring iff, for each element $r\in R$, $\ker_R(\cdot r)$ is
a finitely generated left $R$-module and the intersection of two
finitely generated left ideals is finitely generated}, Proposition
13.3, \cite{Stenstrom-RingQuot}. Each left Noetherian ring is left
coherent but not vice versa.

\begin{theorem}\label{a19Sep10}
The algebra $\mI_1$ is a  left and  right coherent algebra.
\end{theorem}

{\it Proof}. The theorem follows from Theorem \ref{21Sep10},
Theorem \ref{19Sep10}   and Proposition 13.3,
\cite{Stenstrom-RingQuot}. $\Box $

\begin{theorem}\label{A19Sep10}
\begin{enumerate}
\item Every finitely generated left (resp. right) ideal of the
algebra $\mI_1$  is generated by two elements. \item Let $I$ be a
left (resp. right) ideal of $\mI_1$. Then
\begin{enumerate}
\item If $I\not\subseteq F$ then the left (resp. right) ideal $I$
 is generated by two elements. \item If $I\subseteq F$ and $I$ is a
finitely generated left (resp. right) ideal then $I$ is generated
by a single element.
\end{enumerate}
\end{enumerate}
\end{theorem}

{\it Proof}. The algebra $\mI_1$ is self-dual, so it suffices to
prove the statements for left ideals.

 1. Statement 1 follows from statement 2.

2(a). Since $I\not\subseteq F$, we can fix an element $a\in
I\backslash F$. By Theorem \ref{A30May10}.1(a), the factor module
$I/\mI_1a$ has finite length, and so is cyclic, by Proposition
\ref{a21Sep10}. Therefore, ${}_{\mI_1}I$ is generated by two
elements.

2(b). ${}_{\mI_1}F=\bigoplus_{i\in \N}\mI_1e_{0i}$  is the direct
sum of simple isomorphic $\mI_1$-modules $\mI_1e_{0i}=
\bigoplus_{j\in\N}Ke_{ji}$. Then $I\simeq \bigoplus_{i=0}^s \mI_1
e_{0i}=:I'$ for some $s$. Clearly, $I'= \mI_1\th$ where $\th =
e_{00}+e_{11}+\cdots +e_{ss}$ since $\der^s\th = e_{0s}$ and so
$\mI_1 e_{0s}\in \mI_1 \th$ and $\mI_1 \th =
\mI_1(e_{00}+e_{11}+\cdots +e_{s-1,s-1})+\mI_1 e_{0s}$. Using a
similar argument we see that
$$ \mI_1\th =\mI_1(e_{00}+e_{11}+\cdots +e_{s-2,s-2})+\mI_1 e_{0,s-1}+\mI_1
e_{0s}=\cdots = \sum_{i=1}^s \mI_1e_{0j}=I'. \;\;\Box $$


\begin{proposition}\label{a21Sep10}
Every left or right $\mI_1$-module of finite length is cyclic
(i.e.,  generated by a single element).
\end{proposition}

{\it Proof}. The algebra $\mI_1$ is self-dual, so it suffices to
show that every left $\mI_1$-module of finite length $M$ is
cyclic. We use induction on the length $l=l_{\mI_1}(M)$ of the
module $M$. The case $l=1$ is trivial. So, let $l>1$, and we
assume that the statement holds for all $l'<l$. Fix a simple
submodule $U=\mI_1 u \simeq \mI_1 / \ga$ of $M$  where $u\in M$
and  $\ga = \ann_U(u)$. Let $V:= M/ U$. Then $l_{\mI_1}(V)=l-1$,
and, by induction, the $\mI_1$-module $V= \mI_1 \overline{v}
\simeq \mI_1 / \gb$ is cyclic where $\overline{v} \in V$ and $\gb
= \ann_V(\overline{v})$. Fix an element $v\in M$ such that
$\overline{v}=v+U$. Let $\gc = \ga \cap \gb$. Then
$l({}_{\mI_1}(\mI_1/\gc ))<\infty$ since $\mI_1/ \gc$ can be seen
as a submodule of the finite length $\mI_1$-module $\mI_1/ \ga
\bigoplus \mI_1/\gb$. Then $\gc \neq 0$ since
$l{}_{\mI_1}(\mI_1)=\infty$. Moreover, $\gc\not\subseteq F$ since
$l_{\mI_1}(\mI_1/F)=l_{B_1}(B_1)=\infty$. We claim that there
exists an element $a\in \mI_1$ such that $\gc au\neq 0$. Suppose
not, i.e.,   $\gc \mI_1 U=0$, we seek a contradiction. Then $\gc
\mI_1 = F$ or $\mI_1$ since $F$ and $\mI_1$ are the only nonzero
ideals of the algebra $\mI_1$, \cite{algintdif}. The first case is
not possible since $\gc\not\subseteq F$. The second case is not
possible since $u\neq 0$. This finishes the proof of the claim.
 Fix  $c\in \gc$ and $a\in
\mI_1$ such that $cau\neq 0$. Then the  element $w=au+v$ is a
generator for the $\mI_1$-module $M$. Indeed, $cw= cau\neq 0$.
Then, $U=\mI_1 cau\subseteq \mI_1 w$, and so $v\in \mI_1 w$.
Therefore, $M= \mI_1 w$. $\Box $

$\noindent $


\section{Centralizers}\label{CENTRALIZERS}

The centralizers of non-scalar elements of the Weyl algebra $A_1$
share many pleasant properties. Recall Amitsur's well-known
theorem on the centralizer \cite{Amitsur-58} which states that the
 centralizer $\Cen_{A_1}(a)$ of any non-scalar element $a$ of the
 Weyl algebra $A_1$ is a commutative algebra and a free
 $K[a]$-module of finite rank (see also Burchnall and Chaundy \cite{Burchnall-Chaundy-1923}).
  In particular, the centralizer
 $\Cen_{A_1}(a)$ is a commutative finitely generated (hence
 Noetherian)  algebra. It turns out that this result also holds
 for certain generalized Weyl algebra \cite{Bav-AlgAnaliz92},
 \cite{Bav-Ferrara-05},
 \cite{Bav-DixPr5-05}, \cite{Bav-DixPr6-CA-06} and some
 (quantum) algebras see  \cite{Burchnall-Chaundy-1923},
 \cite{Dix}, \cite{Goodearl-Cen-83},
 \cite{Berest-Wilson-MAD-2008}, \cite{Artamonov-Cohn-1999},
  \cite{Mazorchuk-2001ComCentr}, \cite{Hellstrom-Silvestrov-2007}.
  Proposition \ref{a1Jun10} shows that the situation is completely
  different for the algebra $\mI_1$. Theorem \ref{C11Oct10}
  presents in great detail the structure of the centralizers of
  the non-scalar
  elements of the algebra  $\mI_1$, and Corollary \ref{xa10Oct10}
  answers the questions of when the centralizer $\Cen_{\mI_1} (a)$
  is a finitely generated $K[a]$-module where $a\in \mI_1$, or  is
  a finitely generated Noetherian algebra. Corollary
  \ref{a30Oct10} classifies the non-scalar elements of the
  algebra $\mI_1$ such that their centralizers are finitely generated
  algebras. The next proposition will be used in the proof of
  Theorem \ref{C11Oct10}.

For an element $a\in \mI_1$, let $\Cen_{\mI_1}(a):=\{ b\in \mI_1\,
| \, ab=ba\}$ be its  {\em centralizer} in $\mI_1$ and
$\Cen_F(a):=F\cap \Cen_{\mI_1}(a)$.

\begin{proposition}\label{a1Jun10}
\begin{enumerate}
\item Let $\alpha \in K[H]\backslash K$. Then $\Cen_{\mI_1}(\alpha
) = D_1\bigoplus C_\alpha \bigoplus C_\alpha^*$ where the vector
space $C_\alpha :=$ $  \bigoplus_{i\geq 1, j\geq 0} \{ Ke_{i+j,
j}\, | \, j+1$ is a root of the polynomial $\tau^i (\alpha ) -
\alpha\}$,  and $\Cen_{\mI_1}(\alpha )^*=\Cen_{\mI_1}(\alpha )$. In
general, the centralizer $\Cen_{\mI_1}(\alpha ) $ is a
noncommutative, not left and not right Noetherian, not finitely
generated algebra which is not a domain as the following example
shows: $ \Cen_{\mI_1}((H-3/2)^2) = D_1\bigoplus Ke_{ 10}\bigoplus
Ke_{01}$; but $\Cen_{\mI_1}(H^k)=D_1$ is a {\em commutative}, not
Noetherian, not finitely generated algebra which is not a domain
for all $k\geq 1$;  and  $ \Cen_{\mI_1}(H-3/2) = D_1\neq \Cen_{\mI_1}((H-3/2)^2)$.

\item Let $a\in \mI_1$. Then $\dim_K
(\Cen_F(a))<\infty$ iff $a\not\in K[H]+F$. \item
$\Cen_{\mI_1}(\der^i ) = K[\der ]$ and $\Cen_{\mI_1}(\int^i ) =
K[\int ]$ for all $i\geq 1$. \item $\Cen_{\mI_1}(x^i) = K[x]$ for
all $i\geq 1$.
\end{enumerate}
\end{proposition}

{\it Proof}. 1. Recall that $B_1= \mI_1 / F$. Notice that
 $\Cen_{B_1}(\alpha )= K[H]$. This follows from the fact that the
 algebra of $\tau^i$-invariants
 $K[H]^{\tau^i}:=\{ \beta \in K[H]\, |, \tau^i (\beta ) = \beta \}$  is equal to $ K$ for all integers $0\neq i\in \Z$. Since
 $\alpha^* = \alpha$, $\Cen_{\mI_1}(\alpha
)^* = \Cen_{\mI_1}(\alpha )$. The element $\alpha \in D_1$ is a
homogeneous element of the $\Z$-graded algebra $\mI_1$. Therefore,
its centralizer is a homogeneous subalgebra of $\mI_1$. Since
$K[H]\subseteq D_1\subseteq \Cen_{\mI_1}(\alpha )$ and
$\Cen_{B_1}(\alpha )= K[H]$, we see that $\Cen_{\mI_1}(\alpha )
=D_1\bigoplus \bigoplus_{i\geq 1}(C_i\bigoplus C_i^*)$ where
$C_i:= F_{1,i}\cap \Cen_{\mI_1}(\alpha ) $ and $F_{1,i} :=
\bigoplus_{j\geq 0} Ke_{i+j, j}$. Each direct summand $Ke_{i+j,j}$
of the vector space $F_{1, i}$ is a $D_1$-bimodule and $\alpha \in
D_1$, hence $C_i=\sum \{ Ke_{i+j, j}\, | \, \alpha e_{i+j,
j}=e_{i+j, j}\alpha \}$. The equality in the brackets is
equivalent to the equality $\alpha (i+j+1) = \alpha (j+1)$, i.e.,
$(\tau^i (\alpha )-\alpha )(j+1) =0$, i.e.,  $j+1$ is a root of the
polynomial $\tau^i (\alpha ) -\alpha $. Therefore,
$\Cen_{\mI_1}(\alpha ) = D_1\bigoplus C_\alpha \bigoplus
C_\alpha^*$. If $\alpha = H^k$ then $(\tau^i (H^k) - H^k) (j+1) =
(i+j+1)^k-(j+1)^k >0$, and so $\Cen_{\mI_1} (H^k)= D_1$ is a
commutative, not Noetherian, not finitely generated algebra which
is not a domain. If $\alpha = (H-\frac{3}{2})^2$ then $0= (\tau^i
(\alpha ) - \alpha ) (j+1) = i(2j+i-1)$, and so $C_\alpha =
Ke_{10}$. The centralizer $\Cen_{\mI_1} (\alpha )$ is a
noncommutative algebra (since $He_{10}= 2e_{10} \neq e_{10} =
e_{10}H$) and is not a domain (since $e_{10}^2=0$). The factor
algebra $\CD := D_1/ (Ke_{00}+Ke_{11})$ is a commutative, not
Noetherian, not finitely generated algebra. Since the algebra
$\CD$ is the factor algebra algebra of the algebra
$\Cen_{\mI_1}(\alpha )$ modulo the ideal
$\bigoplus_{i,j=0}^1Ke_{ij}$, the algebra $\Cen_{\mI_1} (\alpha )$
is not a finitely generated algebra which is neither left nor right
Noetherian.  If $\alpha = H-3/2$ then the equation $0=(\tau^i(\alpha ) - \alpha ) (j+1)=i$  has no solution since $i\geq 1$. Therefore,  $ \Cen_{\mI_1}(H-3/2) = D_1\neq \Cen_{\mI_1}((H-3/2)^2)$. Notice that, for the Weyl algebra $A_1$, $\Cen_{A_1}(p(a))= \Cen_{A_1}(a)$ for all elements $a\in A_1$ and $p(t) \in K[t]$, \cite{Dix}.

2. Suppose that $a\in K[H]+F$, i.e.,  $a= \alpha
+\sum_{i,j=0}^n\l_{ij}e_{ij}$ for some $n$ where $\alpha \in K[H]$
and $\l_{ij}\in K$. Then $\dim_K(\Cen_F(a))=\infty$ since
$\bigoplus_{i=n+1}^\infty Ke_{ii}\subseteq \Cen_F(a)$. It remains
to show that if $a\not\in K[H]+F$ then $\dim_K(\Cen_F(a))<\infty$.
By (\ref{acan1}), $a=\sum_{i=m}^na_iv_i+a_F$ where
$a_F:=\sum_{i,j=0}^l\l_{ij}e_{ij}$ for some elements $a_i\in
K[H]$, $a_n\neq 0$, $a_m\neq 0$, and $\l_{ij}\in K$. Fix a natural
number, say $N$, such that $N>l$, $Ke_{s+n, t}\ni a_n(H)v_n\cdot
e_{st} = a_n(s+n+1)e_{s+n, t}\neq 0$ and $Ke_{t,s-m}\ni
e_{ts}\cdot a_m(H)v_m= a_m(s+1) e_{t, s-m}\neq 0$ for all $s\geq
N$ and $t\in \N$.

{\em Claim 1. If $m<0$ then $\Cen_F(a) \subseteq
\bigoplus_{j=0}^{N-1}E_{\N , j}$ where $E_{\N ,
j}:=\bigoplus_{i\in \N}Ke_{ij}$.}

Suppose that this is not the case then there exists an element
$c=\sum c_{ij} e_{ij}$ where $c_{ij} \in K$ such that  $t:= \max
\{ j  \, | \, c_{ij}\neq 0$ for some $i\in \N \} \geq N$, we seek
a contradiction. Fix an element $s\in \N$ such that $c_{st}\neq
0$. Then $ ce_{t+|m|, t+|m|}=0$  and $a_Fe_{t+|m|, t+|m|}=0$, and
so
\begin{eqnarray*}
 0&=& e_{ss} \cdot 0\cdot e_{t+|m|, t+|m|}= e_{ss} [c,a]e_{t+|m|, t+|m|}=
 e_{ss} cae_{t+|m|, t+|m|}-e_{ss} ace_{t+|m|, t+|m|} \\
 &=& (\sum c_{sj} e_{sj}) \cdot ( \sum_{i=m}^n a_iv_i )\cdot e_{t+|m|,
 t+|m|}=c_{st} (e_{st} a_mv_m) e_{t+|m|, t+|m|}= c_{st} a_m(t+1) e_{s,
 t+|m|}\\
 &\neq & 0,
\end{eqnarray*}
since $e_{st} a_mv_m= a_m(t+1) e_{s,  t+|m|}\neq 0$, by the choice
of $N$,  where $a_m(t+1)$ is the value of the polynomial $a_m(H)$
at $H=t+1$. This contradiction, $0\neq 0$, finishes the proof of
Claim 1.
 By applying the involution $*$ of the algebra $\mI_1$ to Claim 1
 we obtain the following statement.

{\em Claim 2. If $n>0$ then $\Cen_F(a) \subseteq
\bigoplus_{i=0}^{N-1}E_{i, \N }$ where $E_{i, \N}:=\bigoplus_{j\in
\N}Ke_{ij}$.}

If $m<0$ and $n>0$ then, by Claims 1 and 2,
$$ \Cen_F(a) \subseteq (\bigoplus_{j=0}^{N-1}E_{\N , j})\bigcap (\bigoplus_{i=0}^{N-1}E_{i, \N })
= \bigoplus_{i,j=0}^{N-1} Ke_{ij},$$ and so
$\dim_K(\Cen_F(a))<\infty$.

In view of the involution $*$ of the algebra $\mI_1$, to finish
the proof of statement 2 it suffices to consider the case where
$m<0$ and $m\leq n\leq 0$. By Claim 1, $\Cen_F(a)\subseteq
V:=\bigoplus_{j=0}^{N-1}E_{\N , j}$. Suppose that
$\dim_K(\Cen_F(a))=\infty$  (we seek a contradiction). Fix a
natural number $M$ such that $M>N+|m|$. The subspace $U:=
\bigoplus_{j=0}^{N-1}\bigoplus _{i\geq M}Ke_{ij}$ of $V$ has
finite codimension, i.e.,  $\dim_K(V/U)<\infty$. Therefore, $I:=
\Cen_F(a) \cap U\neq 0$ since $\dim_K(\Cen_F(a)) = \infty$. Choose
a nonzero element, say $u=\sum u_{ij} e_{ij}$ (where $u_{ij}\in
K$), of the intersection $I$. Let $p:= \min \{ i\in \N \, | \,
u_{ij}\neq 0$ for some $j\in \N \}$. Then $p\geq M$, by the choice
of the element $u$. Fix an element $q\in \N$  such that
$u_{pq}\neq 0$. Then $e_{p-|m|, p-|m|}a_F=0$ (since $p-|m|\geq
M-|m|>N+|m|- |m|=N>l$) and $e_{p-|m|, p-|m|}u=0$ (by the choice of
$p$ and since $m<0$), and so
\begin{eqnarray*}
 0&=& e_{p-|m|, p-|m|} \cdot 0\cdot e_{qq}= e_{p-|m|, p-|m|} [a,u]e_{qq}=
 e_{p-|m|, p-|m|}aue_{qq}-e_{p-|m|, p-|m|}uae_{qq} \\
 &=& e_{p-|m|, p-|m|}\cdot  ( \sum_{j=m}^n a_jv_j )\cdot (\sum_{i\geq p} u_{iq}
 e_{iq})= (e_{p-|m|, p-|m|}a_mv_m) \cdot u_{pq}e_{pq}
 \\
 &=&a_m(p-|m|+1) u_{pq}
 e_{p-|m|,q}
 \neq  0,
\end{eqnarray*}
since $e_{p-|m|, p-|m|}a_mv_m\neq 0$, by the choice of $N$,  and
$p-|m|\geq M-|m|>N+|m|-|m|=N$. The contradiction, $0\neq 0$,
proves statement 2.

3. Let us prove the first equality then the second one becomes
obvious:
$$
\Cen_{\mI_1}(\int^i)=\Cen_{\mI_1}((\der^i)^*)=(\Cen_{\mI_1}(\der^i))^*=K[\der
] ^*= K[\int ].$$ Notice that $\Cen_F(\der^i)=0$, $K[\der ]
\subseteq C:=\Cen_{\mI_1}(\der^i)$ and $\Cen_{B_1}(\der^i) =
K[\der, \der^{-1}]$, and therefore $C\cap (F+K[\der ] )=K[\der ]$
and $C\subseteq F+K[\der ] +K[\int ]$. To finish the proof it
suffices to show that $C\subseteq F+K[\der ]$. If $C\not\subseteq
F+K[\der]$ then there exists an element $c=\int +f$ for some
element $f\in F$ (this follows from the fact that $\der^i\int^i=1$
and $K[\der ] \subseteq C$), we seek a contradiction. Necessarily,
$f\neq 0$ as $\int \not\in C$.  Then $ \der c= 1+\der f\in C\cap
(F+K[\der ] )= K[\der ] $, and so $\der f=0$, i.e.,  $f=\sum_{j\geq
0} \l_je_{0j}$ where $\l_j\in K$ and not all $\l_j$ are equal to
zero. Similarly, $c\der = 1-e_{00}+f\der\in C\cap (F+K[\der])=
K[\der]$, and so $e_{00}= f\der = \sum_{j\geq 0} \l_j e_{0, j+1}$,
a contradiction. Then, $C=K[\der]$.

4. Notice that $\Cen_F(x^i) =0$, $K[x]\subseteq C:= \Cen_{\mI_1}
(x^i)$ and $\Cen_{B_1}(x^i) ] =K[x]$. The last equality implies
that $C\subseteq F+K[x]$, then the first two give the result:
$C=\Cen_F(x^i) +K[x]= K[x]$. $\Box $

$\noindent $

The following trivial lemma is a reason why the centralizers of
elements may share  `exotic' properties.

\begin{lemma}\label{a3Jun10}
Let $r$ be an element of a ring $R$ and $\ga := \lann_R(r) \cap
\rann_R(r)$. Then $\ga R\ga \subseteq \Cen_R(r)$.
\end{lemma}


The ideal $F=\bigoplus_{i,j\in \N} Ke_{ij} =\bigoplus_{i,j\in \N}
KE_{ij}\simeq M_\infty (K)$ of the algebra $\mI_1$  admits the
{\em trace} (linear) map
$$\tr : F\ra K, \;\;  \sum \l_{ij} e_{ij}=\sum \mu_{ij}E_{ij} \mapsto
\sum \l_{ii}=\sum \mu_{ii}$$ since $E_{ij} =\frac{i!}{j!}e_{ij}$.
Clearly, for all elements $\alpha \in K[H]$, 
\begin{equation}\label{trdi}
\tr ([\alpha \der^i, F])= \tr ([\alpha \int^i , F]) = \tr ([F,F])
=0, \;\; i\geq 0,
\end{equation}
since $[\alpha (H) v_i, e_{st}]= \alpha(i+s+1) e_{i+s, t}-\alpha
(t+1) e_{s,t-i}$ where $[a,F]:= \{ [a,f]:= af -fa \, | \, f\in
F\}$ for an element $a\in \mI_1$, and $[F,F]$ is the linear
subspace of $F$ generated by all the commutators $[f,g]$ where
$f,g\in F$. Therefore, by (\ref{acan1}) and (\ref{trdi}),
\begin{equation}\label{trI1F}
\tr ([\mI_1, F])=0,
\end{equation}
where $[\mI_1, F]$ is the subspace of $F$ generated by all the
commutators $[a,f]$ where $a\in \mI_1$ and $ f\in F$.

\begin{lemma}\label{a11Oct10}
For all positive integers $i$ and $j$, and for all elements
$f,g\in F$, $[\der^i +f, \int^j+g]\neq 0$.
\end{lemma}

{\it Proof}. Suppose that  $[\der^i +f, \int^j+g]=0$, we seek a
contradiction. Then $[(\der^i +f)^j, (\int^j+g)^i]=0$, so we may
assume that $i=j$ since $(\der^i+f)^j= \der^{ij}+f'$ and
$(\int^j+g)^i= \int^{ij}+g'$ for some elements $f',g'\in F$. Then
the equality $[\der^i +f, \int^i+g]= 0$ implies the equalities (see
(\ref{Iidi}))
$$ e_{00}+e_{11}+\cdots + e_{i-1, i-1}= [\der^i , \int^i]=
-[\der^i , g]-[f,\int^i]-[f,g].$$ Applying the trace map and using
(\ref{trdi}) we get a contradiction, $i=0$.  $\Box $

$\noindent $

The algebra $B_1$ is a subalgebra of the algebra $\CB_1= K(H) [
\der ,\der^{-1} ; \tau ] $ which is the (two-sided) localization
of the algebra $B_1$ at the (left and right) denominator set
$K[H]\backslash \{ 0\}$. A polynomial
$f=\l_nH^n+\l_{n-1}H^{n-1}+\cdots +\l_0\in K[H]$ of degree $n$ is
called a {\em  monic}  polynomial  if the leading coefficient
$\l_n$ of $f$ is $1$. A rational function $h\in K(H)$ is called a
{\em monic}  rational function  if $h=f/g$ for some monic
polynomials $f,g$. A homogeneous element $u=\alpha \der^i$ of
$\CB_1$ is called {\em monic} iff
 $\alpha $ is a monic rational function. We can extend in the obvious way the
 concept
of the  degree of a polynomial to the field of rational functions
by setting,
 $\deg_H (h) = \deg_H (f) - \deg_H (g)$, for $h=f/g\in K(H)$. If $h_1, h_2 \in K(H)$
then $\deg_H (h_1h_2)=\deg_H(h_1) +\deg_H(h_2)$, and
 $\deg_H (h_1+h_2)\leq \max \{ \deg_H(h_1) , \deg_H(h_2)\}$.
  We denote by ${\rm sign} (n)$ and by $|n|$ the {\em sign} and the
{\em absolute value} of $n\in \mathbb{Z}\backslash \{ 0\}$, respectively.

\begin{proposition}\label{A16Oct10}
{\rm (Proposition 2.1, \cite{Bav-DixPr5-05})}
\begin{enumerate}
\item Let $u=\alpha \der^n$ be a monic element of $\CB_1$  with
$n\in \Z\backslash \{  0\}$.  The centralizer
$\Cen_{\CB_1}(u)=K[v,v^{-1}]$ is a Laurent polynomial ring in a
uniquely defined variable $v=\beta \der^{{\rm sign} (n) s}$ where
$s$ is the least positive divisor of $n$ for which there exists a
monic element $\beta =\beta_s \in K(H)$, (necessarily unique) such
that
\begin{eqnarray*}
 \label{bap} \beta \,\tau^s (\beta )\, \tau^{2s}(\beta ) \cdots
\tau^{(n/s -1)s}(\beta )&=&\alpha, \;\; {\rm if}\;\; n>0,\\
\label{bam} \beta \, \tau^{-s} (\beta ) \,\tau^{-2s}(\beta )
\cdots \tau^{-(|n|/s -1)s}(\beta )&=&\alpha, \;\; {\rm if}\;\;
n<0.
\end{eqnarray*}
\item Let $u\in K(H)\backslash K$. Then $\Cen_{\CB_1}(u)=K(H)$.
\end{enumerate}
\end{proposition}

The algebra $B_1= K[H][\der, \der^{-1} ; \tau ]=\bigoplus_{i\in
\Z}K[H]\der^i$ is a $\Z$-graded algebra where $K[H]\der^i$ is the
$i$'th graded component. Each nonzero element $b\in B_1$ is a
unique finite sum $b=\sum \beta_i \der^i$ where $\beta_i\in K[H]$.
Let 
\begin{equation}\label{psipm}
\pi'_+(b):=\max \{ i\in \Z \, | \, \beta_i\neq 0\} \;\; {\rm and
}\;\; \pi'_-(b):=\min \{ i\in \Z \, | \, \beta_i\neq 0\}
\end{equation}
For all elements $a,b\in B_1\backslash \{ 0\}$ and $\alpha \in
K[H]\backslash \{ 0\}$, $\pi'_{\pm}(ab)=
\pi'_{\pm}(a)+\pi'_{\pm}(b)$ and $\pi'_{\pm}(\alpha a)
=\pi'_{\pm}(a)$. The element $\beta_n\der^n$ where $n=\pi'_+(a)$
is called the {\em leading term} of the element $b$, and the
element $\beta_m\der^m$ where $m=\pi'_-(a)$ is called the {\em
least term} of $b$.

\begin{corollary}\label{a16Oct10}
\begin{enumerate}
\item Let $b\in B_1\backslash K$. Suppose that $n= \pi'_+(b)>0$
and $g_1, g_2\in \Cen_{B_1} (b)$ be such that $m = \pi_+'(g_1)=
\pi_+'(g_2)$. Let $\beta_1 \der^m$ and $\beta_2 \der^m$ be the
leading terms of the elements $g_2$ and $g_2$ respectively where
$\beta_1, \beta_2\in K[H]$. Then $K\beta_1 = K\beta_2$. \item Let
$ b = \alpha \der^n$ where $\alpha \in K[H]\backslash K$ and $n\in
\Z\backslash \{ 0\}$. Then $\Cen_{B_1} (b) \cap B_{1, -\sign (n)}
= K$ where $B_{1, -} := K[H][\der^{-1}; \tau^{-1}]\subseteq B_1$
and $B_{1, +} := K[H][\der; \tau ]\subseteq B_1$.
\end{enumerate}
\end{corollary}

{\it Proof}. 1. Let $ \beta \der^n$ be the leading term of the
element $b$. The elements $b$ and $g_1$ (respectively, $b$ and
$g_2$) commute then so do their leading terms. Then the result
follows from Proposition \ref{A16Oct10}.(1) since $\Cen_{B_1}
(\beta \der^n) \subseteq \Cen_{\CB_1} (\beta \der^n)$.

2. Without loss of generality we may assume that $\alpha$ is a
monic since $\Cen_{B_1} (\l b) = \Cen_{B_1} (b)$ for all $\l \in
K^*$. In view of the $(\pm )$-symmetry we may assume that $n>0$.
By Proposition \ref{A16Oct10}, $\Cen_{\CB_1} (b) = K[v,v^{-1}]$
for some element $v =\beta \der^s$ where $s>0$, $s|n$,  $\beta \in
K(H)$, and   $v^\frac{n}{s}=(\beta \der^s)^\frac{n}{s}=\beta
\tau^s (\beta ) \cdots \tau^{(n/s-1)s} (\beta ) \der^n = \alpha
\der^n= b$. Clearly,
$$ 1\leq \deg_H(\alpha ) = \deg_H(\prod_{j=0}^{n/s-1} \tau^{js}
(\beta )) = \sum_{j=0}^{n/s-1} \deg_H(\tau^{js} (\beta )) =
\sum_{j=0}^{n/s-1} \deg_H(\beta ) =\frac{n}{s}\deg_H(\beta ),$$
and so $\deg_H(\beta )>0$. Clearly, $\Cen_{B_1} (b) \cap B_{1,-}
\subseteq \Cen_{\CB_1} (b) \cap B_{1, -} \subseteq K[v^{-1}]\cap
B_1= \bigoplus_{i\in \N} (Kv^{-i}\cap B_1)$. For all integers
$i\geq 1$,
$$ v^{-i} = (\der^{-s}\beta^{-1})^i = \der^{-si} \beta^{-1}
\tau^{s} (\beta^{-1}) \cdots \tau^{s(i-1)}(\beta^{-1}) \not\in
B_1$$ since $\deg_H(\beta^{-1} \tau^{s} (\beta^{-1}) \cdots
\tau^{s(i-1)}(\beta^{-1}))=s\deg_H(\beta^{-1}) = -s\deg_H(\beta )
<0$. Therefore, $\Cen_{B_1} (b) \cap B_{1, -}=K$. $\Box $

$\noindent $

For each element $a\in \mI_1\backslash F$, at least one of the
elements $b_i\in K[H]$ in (\ref{acan1}) is nonzero. Let
\begin{equation}\label{pipmb}
\pi_+(a) :=\max \{ i\in \Z \, | \, b_i\neq 0\}, \;\; \pi_-(a)
:=\min \{ i\in \Z \, | \, b_i\neq 0\}.
\end{equation}
For the  element $a\in \mI_1\backslash F$, the summands $ l_+(a):=b_m v_m$
and $l_-(a):=b_nv_n$ where  $m:= \pi_+(a)$ and $n:= \pi_-(a)$ are called
the {\em largest} and {\em least terms of $a$ modulo} $F$
respectively. It is a useful  observation that if $[a,b]=0$ then
\begin{equation}\label{pipab}
[l_+(a) , l_+(b)],  [l_-(a) , l_-(b)]\in F.
\end{equation}
For an algebra $A$, the subspace $[A, A]$ of $A$ generated by all
the commutators $[a,b]:=ab-ba$ where $a, b\in A$ is called the
{\em commutant} of the algebra $A$. By (\ref{acan1}), $\mI_1=
\bigoplus_{i\in \Z}K[H]v_i \bigoplus \bigoplus_{s,t\in \N}
Ke_{st}$. Let 
\begin{equation}\label{xiII}
\xi : \mI_1 \ra \mI_1
\end{equation}
be the projection onto the direct summand $\bigoplus_{i\in \N}
Ke_{ii}$ of $\mI_1$.

We identify each linear map $b\in \mI_1 \subseteq \End_K(K[x])$
with its $\N \times \N$ matrix with respect to the $K$-basis $\{
x^{[s]}:=\frac{x^s}{s!}\, | \, s\in \N \}$ of the vector space
$K[x]$. For each natural number $d\in \N$, let $e_d:=
e_{00}+e_{11}+\cdots + e_{dd}$ and $e_d':= 1-e_d$. Then $1=
e_d+e_d'$ is the sum of two orthogonal idempotents in the algebra
$\mI_1 \subseteq \End_K(K[x])$. Therefore, $K[x]= \im_{K[x]}(e_d)
\bigoplus \im_{K[x]}(e_d')=K[x]_{\leq d}\bigoplus K[x]_{>d}$ and,
for each element $b\in \mI_1$, 
\begin{equation}\label{bbij}
b=\begin{pmatrix} b_{11} & b_{12}\\ b_{21} & b_{22}
 \end{pmatrix} , b_{11}=e_dbe_d, \;\; b_{12}:=e_dbe_d', \;\;
 b_{21}= e_d'be_d, \;\; b_{22} := e_d'be_d'.
\end{equation}
Notice that $b_{11}\in F_{\leq d} := \bigoplus_{i,j=0}^d Ke_{ij}$,
$b_{12}\in \bigoplus_{i=0}^d\bigoplus_{j>d} Ke_{ij}$, $b_{21}\in
\bigoplus_{i>d}\bigoplus_{j=0}^dKe_{ij}$ and $b_{22} \in
\bigoplus_{i,j>d}Ke_{ij}$. For each natural number $d\in \N$,
consider the algebra $e_d'\mI_1 e_d'\subseteq \mI_1$ where $e_d'$
is its identity element. For each element $b\in e_d'\mI_1 e_d'$,
the presentation (\ref{bbij}) has the form $b=\begin{pmatrix} 0 &0
\\ 0 & b_{22}
 \end{pmatrix}$. The algebras $e_d'\mI_1 e_d'$ will appear in
 Theorem \ref{C11Oct10}.(3) as a large part of the centralizer
 $\Cen_{\mI_1}(a)$ for some elements $a\in \mI_1$. The properties
 of the centralizer  $\Cen_{\mI_1}(a)$ of being a finitely
 generated algebra or a (left or right) Noetherian algebra largely
 depend on similar properties of the algebras  $e_d'\mI_1 e_d'$.
 Let us prove some  results on the algebras $e_d'\mI_1
 e_d'$ that will be used  in the proof of Theorem \ref{C11Oct10}.
 We can easily verify that
\begin{equation}\label{edi}
e_d'\der^ie_d'= e_d'\der^i, \;\; e_d'\int^ie_d'=\int^i e_d', \;\;
i\in \N .
\end{equation}
The map $\mI_1\ra e_d'\mI_1e_d'$, $a\mapsto e_d'ae_d'$, is not an
algebra homomorphism and neither an injective map as
$e_d'e_{00}e_d'=0$; but its restriction to the subalgebra $D_1$
yields the  $K$-algebra epimorphism $D_1\ra e_d'D_1e_d'$,
$a\mapsto e_d'ae_d'$,  with kernel $\bigoplus_{i=0}^dKe_{ii}$.
Notice that $e_d'D_1e_d'= e_d'D_1= D_1e_d'$ is a commutative
non-Noetherian algebra as $\bigoplus_{i>d}Ke_{ii}$ is the direct
sum of nonzero ideals $Ke_{ii}$ of the algebra $D_{1,d}:=
e_d'D_1=K[e_d'H ] \bigoplus \bigoplus_{i>d}Ke_{ii}$.
\begin{lemma}\label{a31Oct10}
Let $d\in \N$. Then
\begin{enumerate}
\item The algebra $e_d'\mI_1e_d'$ is  a finitely generated algebra
which is neither left nor right Noetherian algebra, and
$(e_d'\mI_1e_d')^*=e_d'\mI_1e_d'$. Moreover,  the algebra
$e_d'\mI_1e_d'$ is generated by the elements $e_d'\der$, $ \int
e_d'$, $e_d'H$ and $e_{d+1, d+1}$; and it contains infinite direct
sums of nonzero left and right ideals. \item The algebra
$e_d'\mI_1e_d'$ is a $\Z$-graded algebra which is a homogeneous
subring of the $\Z$-graded algebra $\mI_1$:
$$e_d'\mI_1e_d'=\bigoplus_{i\geq 1} D_{1,d}\der^i\bigoplus D_{1,d}\bigoplus \bigoplus_{i\geq 1} \int^i D_{1,d}$$
where $D_{1,d}:= e_d'D_1$. \item Let $\alpha \in K[H]\backslash
K$. Then $\Cen_{e_d'\mI_1e_d'}(e_d'\alpha ) =
e_d'\Cen_{\mI_1}(\alpha ) e_d'= D_{1, d} \bigoplus C_{\alpha , d}
\bigoplus C_{\alpha , d}^*$ where $C_{\alpha , d}:=
\bigoplus_{i\geq 1, j>d}\{ Ke_{i+j, j} \, | \, j+1$ is a root of
the polynomial $\tau^i (\alpha ) - \alpha \}$.
\end{enumerate}
\end{lemma}

{\it Proof}. 2. The element $e_d'\in D_1$ is a homogeneous element
of the $\Z$-graded algebra $\mI_1$ of graded degree 0. Therefore,
the algebra $e_d'\mI_1e_d'$ is a homogeneous subring of the
algebra $\mI_1$, and, by (\ref{I1iZ}) and (\ref{edi}),
$$ e_d'\mI_1e_d'= e_d'(\sum_{i\geq 1} D_1\der^i+D_1+\sum_{i\geq 1} \int^i D_1) e_d'=
\sum_{i\geq 1} D_{1, d}\der^i+D_{1,d}+\sum_{i\geq 1} \int^i
D_{1,d}, $$ and statement 2 follows.

 1. By (\ref{edi}), $(e_d'\der )^i = e_d'\der^i$ and $ (\int
 e_d')^i = \int^i e_d'$ for all $i\geq 1$. For all $i,j>d$,
 $$ e_{ij} = \int^{i-d-1}e_{d+1, d+1}\der^{i-d-1} =
 \int^{i-d-1}e_d'e_{d+1, d+1}e_d'\der^{i-d-1} =
 (\int e_d')^{i-d-1}e_{d+1, d+1}(e_d'\der )^{i-d-1}. $$
Therefore, by statement 2, the algebra $e_d'\mI_1e_d'$ is
generated by the elements $e_d'\der$, $ \int e_d'$, $e_d'H$ and
$e_{d+1, d+1}$. Since $(e_d')^* = e_d'$, we have
$(e_d'\mI_1e_d')^*=e_d'\mI_1e_d'$. The sum $\bigoplus_{j>d} E_{>d,
j}$ where $E_{>d, j} := \bigoplus_{i>d}Ke_{ij}$ (resp.
$\bigoplus_{i>d} E_{i, >d} $ where $E_{i, >d} :=
\bigoplus_{j>d}Ke_{ij}$) is an infinite  direct sum of nonzero
left ideals $E_{>d, j}$ (resp. right ideals $E_{i, >d}$) of the
algebra $e_d'\mI_1e_d'$. Therefore, the algebra $e_d'\mI_1e_d'$ is
neither left nor right Noetherian.

3. Since the elements $e_d'$ and $\alpha $ commute we have the
inclusion $C':=\Cen_{e_d'\mI_1e_d'}(e_d'\alpha ) \supseteq
e_d'\Cen_{\mI_1}(\alpha ) e_d'$. By Proposition \ref{a1Jun10}.(1),
 $ e_d'\Cen_{\mI_1} (\alpha ) e_d'=D_{1, \alpha } \bigoplus
 C_{\alpha , d}\bigoplus C_{\alpha , d}^*$. The element $e_d'\alpha
 \in D_{1, d}$ is a homogeneous element of the $\Z$-graded algebra
 $e_d'\mI_1e_d'$. Since $\pi (e_d'\alpha ) = \alpha$ where $\pi :
 \mI_1 \ra \mI_1/F = B_1$, $a\mapsto a+F$, we see that $\pi (
 \Cen_{e_d'\mI_1e_d'}(e_d'\alpha ))\subseteq \Cen_{B_1} (\alpha ) =
 K[H] $ (see the proof of Proposition \ref{a1Jun10}.(1)).
 Since $e_d'K[H]\subseteq D_{1, d} \subseteq C'$ and
 $\Cen_{B_1}(\alpha ) = K[H]=\pi (K[e_d' H])$, we see that $C'= D_{1,d}
 \bigoplus C_{\alpha , d} \bigoplus C_{\alpha , d}^*$, and so $C'=
 e_d'\Cen_{\mI_1} (e_d'\alpha ) e_d'$. $\Box $

$\noindent $

Let $K[t]$ be a polynomial algebra in a variable $t$ over the
field $K$ and $K(t)$ be its field of fractions. Let $M$ be a
$K[t]$-module. Then $ {\rm tor}_{K[t]}(M):= \{ m\in M \, | \,
pm=0$ for some $p\in K[t]\backslash \{ 0\} \}$ is the $K[t]$-{\em
torsion} submodule of $M$, it is the sum of all finite dimensional
$K[t]$-submodules of $M$. The {\em rank} of the $K[t]$-submodule
$M$ is $\dim_{K(t)}(K(t)\t_{K[t]}M)$. For each non-scalar element of the algebra $\mI_1$, the following theorem describes its centralizer.

\begin{theorem}\label{C11Oct10}
\begin{enumerate}
\item Let $a\in \mI_1\backslash (K[H]+F)$. Then
\begin{enumerate}
\item The centralizer $\Cen_{\mI_1}(a)$ is a finitely generated
$K[a]$-module of finite rank $\rho \geq 1$. In particular,
$\Cen_{\mI_1}(a)$ is a finitely generated, left and right
Noetherian algebra. \item Moreover, there is a $K[a]$-module
isomorphism
$$\Cen_{\mI_1}(a)\simeq K[a]^\rho \bigoplus \Cen_F(a),$$
and if $-n=\pi_-(a) <0$ (respectively, $m=\pi_+(a)>0$) then $\rho$
divides $n$ (respectively, $m$). \item  $[\Cen_{\mI_1}(a),
\Cen_{\mI_1}(a)]\subseteq \Cen_F(a)$,  $\dim_K(
\Cen_F(a))<\infty$, and  $\Cen_F(a) = {\rm
tor}_{K[a]}(\Cen_{\mI_1}(a))$.
\end{enumerate}

\item Let $a\in (K[H]+F)\backslash (K+F)$, i.e.,  $a=\alpha +f$ for
 some polynomial $\alpha \in K[H]\backslash K$ and $f\in F$,
 $d:=\deg_F(a )$. Then
\begin{enumerate}
\item The centralizer $\Cen_{\mI_1}(a)$ is  not a finitely
generated  $K[a]$-module but has finite rank $\rho :=\deg_H
(\alpha )  $ (as a
 $K[a]$-module).
 \item  $\Cen_{\mI_1}(a)$ is not a finitely generated algebra
  and neither a left nor  right Noetherian algebra. Moreover, it
  contains infinite direct sums of nonzero left and right ideals.
\item There is a $K[a]$-module isomorphism
$$\Cen_{\mI_1}(a)\simeq K[a]^\rho \bigoplus \Cen_F(a).$$
\item $\Cen_{\mI_1} (a) = e_d'K[H]\bigoplus \Cen_F(a)$ and
\begin{eqnarray*}
 \Cen_F(a) &=&\bigoplus_{j>d} Ke_{jj}
 \oplus C_{\alpha, d} \oplus C^*_{\alpha , d} \oplus
 \Cen_{F_{\leq d}}( a_{11})\oplus \\
 & &\left( \bigoplus_{j>d}\ker_{\CH_0} ((a_{11}-\alpha (j+1))\cdot ) e_{0j}\right)\oplus \bigoplus_{i>d}e_{i0} \ker_{\CV_0}
 (\cdot (a_{11}-\alpha (i+1)))
\end{eqnarray*}
is an infinite direct sum of $K[a]$-modules where $a_{11}:=
\sum_{i=0}^d \alpha (i+1)  e_{ii}+f$, $F_{\leq d}
:=\bigoplus_{i,j=0}^dKe_{ij}$, $\CH_0:= \bigoplus_{i=0}^dKe_{i0}$,
$\CV_0:=\bigoplus_{j=0}^dKe_{0j}$ and  $C_{\alpha , d}:=
\bigoplus_{i\geq 1, j> d} \{ Ke_{i+j, j}\, | \, j+1$ is a root of
the polynomial $\tau^i (\alpha ) - \alpha \}$. \item
$[\Cen_{\mI_1}(a), \Cen_{\mI_1}(a)]\subseteq \Cen_F(a)$,
$\dim_K(\Cen_F(a))=\infty$,  $\Cen_F(a)={\rm
tor}_{K[a]}(\Cen_{\mI_1}(a))$.
\end{enumerate}
\item  Let $a\in (K+F)\backslash K$, i.e.,  $a= \l +f$ for some
elements $\l \in K$ and $f = \sum_{i,j=0}^d \l_{ij}e_{ij} \in
F\backslash \{ 0 \}$ where $d:= \deg_F(f)$; let $e_d=
e_{00}+\cdots + e_{dd}$ and $e_d'= 1-e_d$. Then
\begin{enumerate}
\item The centralizer $\Cen_{\mI_1} (a)$ is a finitely generated
algebra which  is neither a left nor right Noetherian algebra.
Moreover, it contains infinite direct sums of nonzero left and
right ideals. In particular, $\Cen_{\mI_1} (a) $ is not a finitely
generated $K[a]$-module. \item $ \Cen_{\mI_1}(a) = \Cen_{F_{\leq
d}} (f) \oplus \bigoplus_{j>d} \CK e_{0j} \oplus
\bigoplus_{i>d}e_{i0} \CK'\oplus e_d'\mI_1e_d'$ where $\CK :=
\ker_{\oplus_{i=0}^dKe_{i0}} (f\cdot)$ and $\CK':= \ker
_{\oplus_{j=0}^dKe_{0j}} (\cdot f)$. \item $[\Cen_{\mI_1} (a),
\Cen_{\mI_1} (a)]\not\subseteq \Cen_F(a)$, $\dim_K(\Cen_F(a))=\infty$, ${\rm tor}_{K[a]}(\Cen_{\mI_1}(a))=\Cen_{\mI_1}(a)\neq \Cen_{F}(a)$.

\end{enumerate}
\end{enumerate}
\end{theorem}

{\it Proof}. 1. Since $a\not\in K[H]+F$ then either
$-n:=\pi_-(a)<0$ or $m:=\pi_+(a)>0$ (or both). In view of $(\pm
)$-symmetry, let us assume that $-n:=\pi_-(a)<0$, i.e.,  $a= \alpha
\der^n+\cdots$ where $\alpha \der^n$ is the least term of the
element $a$ modulo $F$ and  $\alpha \in K[H]\backslash \{ 0\}$.
Let $C:=\Cen_{\mI_1}(a)$. We claim that 
\begin{equation}\label{CIF}
C\cap \mI_{1, +}\subseteq F
\end{equation}
where $\mI_{1, +}:=\bigoplus_{i\geq 1} K[H]\int^i \bigoplus F$. If
$\alpha \in K[H]\backslash K$ this follows from (\ref{pipab}) and
Corollary \ref{a16Oct10}.(2). If $\alpha \in K\backslash \{ 0\}$
then we may assume that $\alpha =1$, by multiplying the element
$a$ by $\alpha^{-1}$. Suppose that the inclusion (\ref{CIF})
fails, and so there is an element $c\in C\cap \mI_{1, +}$ with
$l:= \pi_-(c) \geq  1$. Let $\beta \int^l$ be the least term of
the element $c$ modulo $F$ where $\beta\in K[H]\backslash \{ 0\}$. Since $\Cen_{B_1}(\pi (a) = \der^n) =
K[\der , \der^{-1}]$, we must have $\beta \in K^*$. Without loss
of generality we may assume that $\beta =1$ (by dividing $c$ by
$\beta $). According to (\ref{acan1}), the elements $a$ and $c$
can be uniquely written as the sums $a= \der^n +a'+f$ and $c =
\int^l+c'+g$ where $f,g \in F$, $a'$ is the sum $\sum b_i v_i$ in
the decomposition (\ref{acan1}) without the least term $\der^n$,
and similarly $c'$ is defined. Replacing $a$ and $c$ by $a^l$ and
$c^n$ respectively we may assume that $n=l$ in the presentations
above, i.e.,  $a= \der^n +a'+f$ and $c = \int^n +c'+g$. Notice that
$\xi ([\der^n , c'])=\xi ([a', \int^n]) = \xi ( [ a', c'])=0$ as
the elements in the brackets are sums of homogeneous elements of
the $\Z$-graded algebra $\mI_1$ of {\em positive} graded degrees.
Clearly, $ \tr \circ \xi |_F = \tr$. By (\ref{trI1F}), $\tr \circ
\xi ([F, \mI_1]) = \tr ([F, \mI_1]) =0 $. Now, applying the map
$\tr \circ \xi$ to the equality
$$ 0=[a,c]= [ \der^n , \int^n]+[a', \int^n]+[f,\int^n]+ [ \der^n ,
c']+ [ a',c']+[f,c']+ [a,g]$$ we get a contradiction:
$$ 0= \tr \circ \xi ( [ \der^n \int^n]) = \tr \circ \xi
(e_{00}+e_{11}+\cdots + e_{n-1, n-1}) = n.$$ Therefore,
(\ref{CIF}) holds. The set $C\backslash C_F$ is a multiplicative
monoid where $C_F:=\Cen_F(a)$. We denote by $\kappa$ the
composition of the two monoid homomorphisms $\pi_-:C\backslash
C_F\ra -\N$ (see (\ref{CIF})) and $\Z \ra \Z/ n \Z$, and so
$\kappa$ is a monoid homomorphism ($\kappa (uv) = \kappa (u)
+\kappa (v)$ and $\kappa (1) = 0$). Therefore, its image $G:= \im
(\kappa )$ is a cyclic subgroup of order, say $\rho$, and $\rho
|n$. Let $ G = \{ m_1=0, m_2, \ldots , m_\rho \}$. Then for each
element $m_i$ we choose an element $g_i \in C$ so that $\kappa
(g_i) = m_i$ and the number $\pi_-(g_i) \in -\N$ is the largest
possible  integer. We may choose $g_1=1$, by (\ref{CIF}). Let us
show that the $K[a]$-submodule $M:= \sum_{i=1}^\rho K[a]g_1$ of
$C$ is free, i.e.,  $M= \bigoplus_{i=1}^\rho K[a]g_i$,  and $C=
M\bigoplus C_F$.

Suppose that $\v_1g_1+\cdots +\v_\rho g_\rho\in F$ for some
elements $\v_i\in K[a]$, not all of which are equal to zero, we
seek a contradiction. Then there exist nonzero terms $\v_ig_i$ and
$\v_jg_j$ such that $\pi_- (\v_ig_i) = \pi_-(\v_jg_j)$, and so
$\kappa (g_i) = \kappa (g_j)$. This contradicts to the choice of
the elements $g_1, \ldots , g_\rho$. Therefore,
$M=\bigoplus_{i=1}^\rho K[a]g_i$ and $M\cap C_F=0$. To finish the
proof  the claim it remains to show that $C= M+C_F$. Choose an
arbitrary nonzero element $g\in C$. If $\pi_- (g)=0$ then, by
Corollary \ref{a16Oct10}.(1) and (\ref{CIF}), $g\in K+C_F =
Kg_1+C_F$. If $\pi_- (g) =k<0$ then there exists $g_i$ such that
$\kappa (g) = \kappa (g_i)$, and so $ k= \pi_- (a^s g_i)$ for some
natural number $s\in \N$ (by the choice of the elements $g_1,
\ldots , g_\rho$). By Corollary \ref{a16Oct10}.(1), there exists
$\l \in K$ such that $\pi_- (g-\l a^sg_i) >k$. Using induction on
$|k|$ or repeating the same argument several times we see that $C=
M+C_F$. Therefore, $C= M\bigoplus C_F$ where $C_F$ is a finite
dimensional ideal of the algebra $C$ (Proposition
\ref{a1Jun10}.(2)). This implies that $\rho$ is the rank of the
$K[a]$-module $C$ and $\Cen_F(a) = {\rm tor}_{K[a]}(C)$. In
particular, $C$ ia finitely generated $K[a]$-module, and so $C$ is
a finitely generated left and right Noetherian algebra.

Let us show that $uv-vu\in C_F$ for all elements $u,v\in C$. It
suffices  to show that  $uv-vu\in F$ as  $uv-vu\in C$ and $C_F=
C\cap F$. Choose an element $g\in C$ such that $\kappa (g)$ is a
generator for the group $G$. Denote by $E$ the subalgebra
(necessarily commutative) of $C$ generated by the elements $a$ and
$g$. Notice that $C_F$ is an ideal of the algebra $C$. By the
choice of $g$, the $K[a]$-module $C/ (E+C_F)$ is {\em finite
dimensional}. Therefore, there exist nonzero elements $P,Q\in
K[a]$ such that $Pu, Qv \in E+C_F$. Then
$$ PQuv\equiv (Pu) (Qv)  \equiv (Qv) (Pu)\equiv PQvu \mod C_F,$$
i.e.,  $PQ(uv-vu)\in C_F$. Since $PQ$ is a nonzero element of the
algebra $K[a]$ (which is isomorphic to a polynomial algebra over
$K$ in a single variable, since $a\not\in F$) and $F\cap K[a]=0$,
we see that $PQ\not\in F$, then $uv-vu \in C_F$ since the algebra
$B_1$ is a domain.

2. Let $C:= \Cen_{\mI_1}(a)$ and $C_F:= \Cen_F(a)$. Let $b\in
\mI_1$. By (\ref{bbij}), $b=\begin{pmatrix}
 b_{11}& b_{12}\\ b_{21} & b_{22}
 \end{pmatrix}$ and $a=\begin{pmatrix}
 a_{11}& 0\\ 0 & a_{22}
 \end{pmatrix}$ where $a_{11}:=
\sum_{i=0}^d \alpha (i+1)  e_{ii}+f$ and $\alpha (i+1)$ is the
value of the polynomial $\alpha (H)$ at $H=i+1$, $a_{22}=
e_d'ae_d'= e_d'\alpha = \alpha e_d'\in e_d'\mI_1e_d'$. Then $b\in
C$ iff $\begin{pmatrix}
 a_{11}& 0\\ 0 & a_{22}
 \end{pmatrix}\begin{pmatrix}
 b_{11}& b_{12}\\ b_{21} & b_{22}
 \end{pmatrix}=\begin{pmatrix}
 b_{11}& b_{12}\\ b_{21} & b_{22}
 \end{pmatrix}\begin{pmatrix}
 a_{11}& 0\\ 0 & a_{22}
 \end{pmatrix}$ iff $b_{11}\in \Cen_{F_{\leq d}}(a_{11})$,
 $b_{22}\in \Cen_{e_d'\mI_1e_d'} (e_d'\alpha )$, $a_{11} b_{12} =
 b_{12} a_{22}$ and $ a_{22} b_{21} = b_{21} a_{11}$. Notice that
 $\begin{pmatrix}
 b_{11}& 0\\ 0 & 0
 \end{pmatrix}, \begin{pmatrix}
 0& b_{12}\\ 0 & 0
 \end{pmatrix} , \begin{pmatrix}
 0& 0\\ b_{21} & 0
 \end{pmatrix}\in C_F$. Then, by Lemma \ref{a31Oct10}.(3),
\begin{equation}\label{CedH}
C=e_d'K[H]\bigoplus C_F
\end{equation}
and $e_d'K[H]$ is a free $K[a]$-module (i.e.,  $K[\alpha ] $-module
since $ae_d'= \alpha e_d'$) of rank $\rho = \deg_H(\alpha )$.
  As a $K[a]$-bimodule (in particular, as a $K[\ad
(a)]$-module where ${\rm ad} (a):= [a, \cdot ]$ is the inner
derivation of the algebra $\mI_1$) the ideal  $F$ is the direct sum of four
sub-bimodules:
$$ F_{\leq d},\;\;  F_{>d} :=\bigoplus_{i,j>d} Ke_{ij}, \;\; \CH
:= \bigoplus_{i=0}^d\bigoplus_{j>d} Ke_{ij}, \;\; \CV :=
\bigoplus_{i>d}\bigoplus_{j=0}^d Ke_{ij}.$$ Then
$$ \Cen_F(a) = \ker_{F_{\leq d}} (\ad (a)) \bigoplus
\ker_{F_{>d}}(\ad (a))\bigoplus \ker_{\CH }(\ad (a)) \bigoplus
\ker_{\CV}(\ad (a)).$$ $\ker_{F_{\leq d}} (\ad (a)) $ is the
centralizer of the $(d+1) \times (d+1)$ matrix $a_{11}$ in the
algebra $F_{\leq d}$ of  $(d+1) \times (d+1)$ matrices.
$\ker_{F_{>d}}(\ad (a))=\ker_{F_{>d}}(\ad (\alpha
))=\bigoplus_{j>d} Ke_{jj}
 \bigoplus C_{\alpha, d} \bigoplus C^*_{\alpha , d} $, by
 Proposition \ref{a1Jun10}.(1).  The $K[a]$-bimodule $\CH =\bigoplus_{j>d}\CH_j$ is the
 direct sum of finite dimensional $K[a]$-bimodules $\CH_j :=
 \bigoplus_{i=0}^dKe_{ij}$, and the action of the map $\ad (a)$ on
 $\CH_j$ is equal to $(a_{11}-\alpha (j+1))_{\CH_j}\cdot$.
 Therefore,
\begin{eqnarray*}
 \ker_{\CH_j} (\ad (a)) &=&\ker_{\CH_j}((a_{11}-\alpha (j+1))\cdot)
 =\ker_{\CH_0}((a_{11}-\alpha (j+1))\cdot )e_{0j},  \\
 \ker_{\CH} (\ad (a)) &=&\bigoplus_{j>d} \ker_{\CH_0} ((a_{11}-\alpha (j+1))\cdot
 )e_{0j}.
\end{eqnarray*}
Similarly, the $K[a]$-bimodule $\CV = \bigoplus_{i>d}\CV_i$ is the
direct sum of finite dimensional $K[a]$-bimodules $\CV_i =
\bigoplus_{j=0}^d Ke_{ij}$ and the action of the map $\ad (a)$ on
$\CV_i$ is equal to $\cdot (a_{11}-\alpha (i+1))$. Therefore,
\begin{eqnarray*}
 \ker_{\CV_i} (\ad (a)) &=&\ker_{\CV_i}(\cdot (a_{11}-\alpha (i+1)))
 =e_{i0}\ker_{\CV_0}(\cdot (a_{11}-\alpha (i+1))),  \\
 \ker_{\CV} (\ad (a)) &=&\bigoplus_{i>d} e_{i0}\ker_{\CV_0} (\cdot (a_{11}-\alpha (i+1))
 ).
\end{eqnarray*}
It is obvious that $\Cen_F(a)$ is a $K[a]$-torsion, infinite
dimensional, not finitely generated $K[a]$-module. Therefore, by
(\ref{CedH}), $C$ is not a finitely generated $K[a]$-module, ${\rm
tor}_{K[a]}(C) = C_F$, and the $K[a]$-module $C$ has rank $\rho =
\deg_H (\alpha )$. By (\ref{CedH}), $[C,C]\subseteq C_F$. Recall
that the left $\mI_1$-module $F$ is the direct sum
$\bigoplus_{j\in \N} E_{\N , j}$ of nonzero left ideals $E_{\N ,
j} = \bigoplus_{i\in \N} Ke_{ij}$. Similarly, the right
$\mI_1$-module $F$ is the direct sum $\bigoplus_{i\in \N} E_{i, \N
}$ of nonzero right  ideals $E_{i, \N } = \bigoplus_{j\in \N}
Ke_{ij}$. Since 
\begin{equation}\label{eijEij}
e_{jj}\in E_{\N , j} \cap  \Cen_{\mI_1} (a)  \;\; {\rm for \;
all}\; j>d;\;\;\;  e_{ii}\in E_{i, \N } \cap \Cen_{\mI_1} (a) \;\;
{\rm for \; all}\; i>d,
\end{equation}
the sums $\bigoplus_{j>d} ( E_{\N , j} \cap  \Cen_{\mI_1} (a))$
and $\bigoplus_{i>d} (E_{i, \N } \cap \Cen_{\mI_1} (a))$ are
infinite direct sums of nonzero left and right ideals of the
algebra $C$ respectively. Therefore, the algebra $C$ is neither
left nor right Noetherian. To finish the proof of statement 2 it
remains to show that the algebra $C$ is not finitely generated.
Suppose that $S$ is a finite set of algebra generators for the
algebra $C$, we seek a contradiction. We may assume that $e_d'H
\in S$. Then, by (\ref{CedH}), we may assume that $S= \{ e_d'H\}
\cup S_F$ where $S_F$ is a finite subset of $C_F$. Then we can fix
a natural number $n$ such that $S_F\subseteq F_{\leq n}$. Then $C
= K\langle S \rangle \subseteq K \langle F_{\leq n} , e_d'H\rangle
= K[e_d' H] \bigoplus F_{\leq n}$, and so $C_F \subseteq F_{\leq
n}$, a contradiction, since $e_{n+1, n+1} \in C_F\backslash
F_{\leq n }$.

3. Let $C:= \Cen_{\mI_1}(a)$ and $C_F:= \Cen_F(a)$. Since $C=
\Cen_{\mI_1}(f)$, we may assume that $a=f\in F\backslash K$. Let
$R$ be the RHS in the equality in statement 3(b). Let $b\in
\mI_1$. By (\ref{bbij}), $b=\begin{pmatrix}
 b_{11}& b_{12}\\ b_{21} & b_{22}
 \end{pmatrix}$ and $f=\begin{pmatrix}
 f& 0\\ 0 & 0
 \end{pmatrix}$.  Then $b\in C$ iff
$\begin{pmatrix}
 f& 0\\ 0 & 0
 \end{pmatrix}\begin{pmatrix}
 b_{11}& b_{12}\\ b_{21} & b_{22}
 \end{pmatrix}=\begin{pmatrix}
 b_{11}& b_{12}\\ b_{21} & b_{22}
 \end{pmatrix}\begin{pmatrix}
 f& 0\\ 0 & 0
 \end{pmatrix}$ iff $fb_{11}=b_{11}f$,
 $fb_{12}=0$ and  $b_{12}f =0$ iff $b\in R$. It follows from the
 equality $C= R$ and Lemma \ref{a31Oct10} that the centralizer $C$
 is a finitely generated algebra which is generated by the finite
 dimensional subspaces $\Cen_{F_{\leq d}}(f)$, $\CK e_{0, d+1}$,
 $e_{d+1, 0}\CK'$  and the elements $e_d'$, $e_d'\der$,
 $\int e_d'$ and  $e_{d+1, d+1}$ (since, for all $i>d$,
 $\CK e_{0i} = \CK e_{0, d+1} (e_d'\der )^{i-d-1}$ and
 $e_{i0} \CK' = (\int e_d')^{i-d-1} e_{d+1, 0} \CK ' $).
 Clearly, $[C, C]\not\subseteq C_F$ since
 $$ [ e_d' He_d', e_d' \der e_d'] \equiv [ H, \der ] \equiv - \der
 \mod F,$$
 and $\dim_K(C) \geq \dim_K(e_d' D_1)=\infty$. By (\ref{eijEij}),
 the sums $\bigoplus_{j>d} ( E_{\N , j} \cap  \Cen_{\mI_1}
(a))$ and $\bigoplus_{i>d} (E_{i, \N } \cap \Cen_{\mI_1} (a))$ are
infinite direct sums of nonzero left and right ideals of the
algebra $C$ respectively. Therefore, the algebra $C$ is neither
left nor right Noetherian. Clearly, ${\rm tor}_{K[a]}(C)=C\neq C_F$. $\Box $


\begin{corollary}\label{xa10Oct10}
Let $a\in \mI_1\backslash K$. Then the following statements are
equivalent.
\begin{enumerate}
\item $a\not\in K[H]+F$. \item  $\Cen_{\mI_1}(a)$ is a finitely
generated $K[a]$-module. \item $\Cen_{\mI_1}(a)$ is a left
Noetherian algebra. \item $\Cen_{\mI_1}(a)$ is a right Noetherian
algebra. \item $\Cen_{\mI_1}(a)$ is a finitely generated  and
Noetherian algebra.
\end{enumerate}
\end{corollary}



\begin{corollary}\label{a30Oct10}
Let $a\in \mI_1\backslash K$. Then
\begin{enumerate}
\item $\Cen_{\mI_1}(a)$ is a finitely generated algebra iff
$a\not\in (K[H]+F)\backslash (K+F)$.
 \item  $\Cen_{\mI_1}(a)$ is a finitely
generated not Noetherian/not left Noetherian/not right Noetherian
algebra iff $a\in (K+F)\backslash K$.
 \item  The algebra $\Cen_{\mI_1}(a)$ is not  finitely
generated, not Noetherian/not left Noetherian/not right Noetherian
 iff $a\in (K[H]+F)\backslash (K+F)$.
\end{enumerate}
\end{corollary}


\section{The simple $\mI_1$-module $K[x]$ and a
classification of elements of the algebra
$\mI_1$}\label{SIKXM}

In this section, a formula for the index $\ind_{K[x]}(a)$ of
elements $a\in \mI_1\backslash F$ is found (Proposition
\ref{a12Jun10}). Recall that $K[x]$ is the unique (up to
isomorphism) {\em faithful} simple $\mI_1$-module and
$\mI_1\subseteq \End_K(K[x])$. Classifications of elements $a\in
\mI_1$ are given such that the map $a_{K[x]}: K[x]\ra K[x]$, $p\mapsto
ap$, is a bijection (Theorem \ref{12Jun10}), a surjection (Theorem
\ref{15Jun10}) or an injection (Theorem \ref{16Jun10}). In case
when the map $a_{K[x]}$ is a bijection, an explicit inversion
formula  is found (Theorem \ref{12Jun10}.(4)). As a result we have
a formula for the unique polynomial solution $q=a^{-1}*p$ of the
polynomial integro-differential operator equation $a*p=q$ where
$p,q\in K[x]$ and $a\in \mI_1^0:= \mI_1\cap \Aut_K(K[x])$. The
monoid
 $\mI_1^0$ is much more massive set  than the group $\mI_1^*$ of units of
the algebra $\mI_1$ (Theorem \ref{12Jun10}.(1,2), Theorem
\ref{14Jun10}.(6)). In  case when the map $a_{K[x]}$ is surjective
or injective we found its  kernel and cokernel respectively
(Theorem \ref{15Jun10} and Theorem \ref{16Jun10}).

Each nonzero element $u$ of the skew Laurent polynomial algebra
$B_1=K[H][\der, \der^{-1}; \tau ]$ (where $\tau (H) = H+1$) is the
unique sum
$u=\alpha_s(\der^{-1})^s+\alpha_{s+1}(\der^{-1})^{s+1}+\cdots +
\alpha_d(\der^{-1})^d$ where all $\alpha_i\in K[H]$, $\alpha_d\neq
0$, and $\alpha_d(\der^{-1})^d$ is the {\em leading term} of the
element $u$. The integer $\deg_{\der^{-1}}(u)=d$ is called the
{\em degree} of the element $u$ in the noncommutative variable
$\der^{-1}$,
 $\deg_{\der^{-1}}(0):=-\infty$. For all elements $u,v\in B_1$,
 $\deg_{\der^{-1}}(uv) = \deg_{\der^{-1}}(u)+\deg_{\der^{-1}}(v)$,
  $\deg_{\der^{-1}}(u+v) \leq \max \{  \deg_{\der^{-1}}(u),
 \deg_{\der^{-1}}(v)\}$. Therefore, the minus  degree function
 $$-\deg_{\der^{-1}}: \Frac (B_1)\ra \Z \cup \{ \infty \}, \;\;
 s^{-1}u\mapsto \deg_{\der^{-1}}(s) -\deg_{\der^{-1}}(u), $$ is a
 discrete valuation on the skew  field of fractions $\Frac (B_1)$ of the algebra $B_1$
 such that $\deg_{\der^{-1}}(\alpha ) =0$ for all elements $0\neq
 \alpha \in K(H)$.

\begin{proposition}\label{a12Jun10}
Let $a\in \mI_1\backslash F$.
\begin{enumerate}
\item $\ind_{K[x]}(a) = -\deg_{\der^{-1}}(a+F)$ where $a+F\in
\mI_1/F = B_1$. \item Let $i= \ind_{K[x]}(a)$. Then
\begin{enumerate}
\item $i\geq 0$ iff $a=\der^ia'$ for some element $a'\in
\mI_1\backslash F$ with $\ind_{K[x]}(a')=0$. \item $i\leq 0$ iff
$a=a'\int^{|i|}$ for some element $a'\in \mI_1\backslash F$ with
$\ind_{K[x]}(a')=0$.
\end{enumerate}
\item Let $b\in \mI_1\backslash F$ such that $a+b\in
\mI_1\backslash F$. Then $\ind_{K[x]} (a+b) \geq \min \{
\ind_{K[x]} (a), \ind_{K[x]} (b)\}$. Therefore, the index map
$\ind_{K[x]}:\Frac (A_1) \ra \Z \cup \{ \infty \}$ is a discrete
valuation.
 \item $\ind_{K[x]}(\s (a))=\ind_{K[x]}(a)$ for all automorphisms
 $\s \in \Aut_{K-{\rm alg}}(\mI_1)$.
\end{enumerate}
\end{proposition}

{\it Remark}. In statement 2(a), $a=\der^ia'$ $(i\geq 1)$ does not
imply $a=a''\der^i$ for some element (necessarily) $a''\in
\mI_1\backslash  F$ with $\ind_{K[x]}(a'')=0$; eg,  $a=\der^i
(1+e_{i0})$ since $e_{00}\not\in \im (\cdot \der^i)$. Similarly,
in statement 2(b), $a=a'\int^{|i|}$ $(i\leq -1)$ does not imply
$a=\int^{|i|}a''$ for some element (necessarily) $a''\in
\mI_1\backslash  F$ with $\ind_{K[x]}(a'')=0$; eg,  $a=
(1+e_{0i})\int^{|i|}$ since $e_{00}\not\in \im (\int^{|i|}\cdot
)$.

$\noindent $

{\it Proof}. 1. Let $a$ be as in (\ref{acan}) and
$d=-\deg_{\der^{-1}}(a+F)$. By Lemma \ref{D30May10}, without loss
of generality we may assume that all $\l_{ij}=0$. Notice that the
maps $\int^i\cdot$, $\cdot \der^i:\mI_1\ra \mI_1$ are injections
for all natural numbers $i\geq 1$ (since $\der^i\int^i =1$) such
that $\int^i (\mI_1\backslash F) \subseteq \mI_1\backslash F$ and
 $(\mI_1\backslash F)\der^i \subseteq \mI_1\backslash F$.  If $d\geq
 0$ then $a=\der^d\int^da$ and $\ind_{K[x]}(a) =
 \ind_{K[x]}(\der^d) +\ind_{K[x]}(\int^da)=
 d+\ind_{K[x]}(\int^da)$. If $d\leq
 0$ then $a=a\der^{-d}\int^{-d}$ and $\ind_{K[x]}(a) =
 \ind_{K[x]}(\int^{-d}) +\ind_{K[x]}(a\der^{-d})=
 d+\ind_{K[x]}(a\der^{-d})$.
 Therefore, to finish the proof of statement 1 it suffices to show that if $d=0$
 then $\ind_{K[x]}(a)=0$. Since $d=0$, $a=\sum_{i\geq 1}
 a_{-i}\der^i +a_0$ with $a_0\neq 0$. The element $a$ respects the
 finite dimensional filtration $\{ K[x]_{\leq
 i}\}_{i\in \N}$ of the $\mI_1$-module
 $K[x]$, that is,  $aK[x]_{\leq i}\subseteq K[x]_{\leq i}$ for all
 $i\in \N$. Notice that $a_0(H)*x^i=a_0(i+1)x^i$ for all $i\in
 \N$. Fix a natural number, say $N$, such that $a_0(i+1)\neq 0$
 for all $i>N$. Then the linear map $a\cdot :V:=K[x]/K[x]_{\leq
 N}\ra V$ is a bijection. It follows from the short exact sequence
 of $K[a]$-modules $0\ra K[x]_{\leq N}\ra K[x]\ra V\ra 0$ that
 $\ind_{K[x]}(a) = \ind_{K[x]_{\leq N}}(a)+\ind_V(a)=0+0=0$.

2. Statement 2(a) (resp. 2(b)) follows from statement 1,
(\ref{acan}),  and the fact that $\der^iF = F$ (resp.
$F\int^{|i|}=F$ and $\der^{|i|}\int^{|i|}=1$), see the proof of
statement 1 where statement 2 was, in fact, proved.

3. Statement 3 follows from statement 1.

4. The simple $\mI_1$-module $K[x]$ is the only (up to
isomorphism) {\em faithful} simple $\mI_1$-module (Theorem
\ref{31May10}). Therefore, ${}^\s K[x]\simeq K[x]$ for all
automorphisms $\s $ of the algebra $\mI_1$. Then $\ind_{K[x]}(\s
(a)) = \ind_{{}^\s K[x]}(a) = \ind_{K[x]}(a)$, by (\ref{indMsa}).
$\Box $

$\noindent $

{\bf Classification of elements $a\in \mI_1$ such that $a_{K[x]}$
is a bijection}.  Let $\Aut_K(K[x])$ be the group of all
invertible $K$-linear maps in the vector space $K[x]$, i.e.,  it is
the group of units of the algebra $\End_K(K[x])$. The next theorem
describes the intersection $\mI_1^0 := \mI_1\cap \Aut_K(K[x])$ and
gives an inversion formula in $\Aut_K(K[x])$   for each element in
the intersection. The set $\mI_1^0$ is a multiplicative monoid the
elements of which are (left and right) regular elements of the
algebra $\mI_1$ (i.e.,  non-zero-divisors). We will see that the
multiplicatively closed set $\mI_1^0$ is the {\em largest} (with
respect to inclusion) right Ore set that consists of regular
elements  (Theorem \ref{18Jun10}.(2)). It is obvious that if the
map $a_{K[x]}$ is a surjection  (where $a\in \mI_1$) then
$a\not\in F$. Recall that $\mI_1^*=K^* (1+F)^*$, \cite{algintdif},
where $\mI_1^*$ is the group of units of the algebra $\mI_1$.

\begin{theorem}\label{12Jun10}
Let $a\in \mI_1\backslash F$ and $d:=\deg_F(a)$. Then
\begin{enumerate}
\item $a_{K[x]}\in \Aut_K(K[x])$ iff $a=\sum_{i\geq 1}
a_{-i}\der^i +a_0+\sum_{i,j\in \N}\l_{ij}e_{ij}$ (see
(\ref{acan})), $b:= a|_{K[x]_{\leq d}}\in \GL (K[x]_{\leq d})$ and
$s+1$ is not a root of the nonzero  polynomial $a_0\in K[H]$ for
all natural numbers $s>d$. \item $\mI_1^0 := \mI_1\cap
\Aut_K(K[x])\supsetneqq \mI_1^* = K^* (1+F)^*$. \item If
$\ind_{K[x]}(a) =0$ then the following statements are equivalent:
\begin{enumerate}
\item $a_{K[x]}$ is an injection, \item $a_{K[x]}$ is a
surjection, \item $a_{K[x]}$ is a bijection.
\end{enumerate}
\item {\rm (Inversion Formula)} Suppose that $a_{K[x]}\in
\Aut_K(K[x])$, $a=a_-+a_0+\sum \l_{ij}e_{ij}$ where $a_-=
\sum_{i\geq 1} a_{-i}\der^i$. Then according to the decomposition
$K[x]= K[x]_{\leq d}\bigoplus K[x]_{>d}$ where
$K[x]_{>d}:=\bigoplus_{i>d}Kx^i$ the map $a_{K[x]}$ is the matrix
 $\bigl(
\begin{smallmatrix} a_{11} & a_{12}\\ 0 & a_{22}
\end{smallmatrix}\bigr)$
 where $a_{11}:= a|_{K[x]_{\leq d}}$, $a_{22}:=a|_{K[x]/K[x]_{\leq
 d}}=(a_-+a_0)|_{K[x]/K[x]_{\leq d}}$, the matrix of the linear
 map $a_{12}: K[x]/K[x]_{\leq d}\ra K[x]_{\leq d}$ has only finite
 number of nonzero entries with respect to the monomial bases of
 the vector spaces, and
\begin{equation}\label{am1b}
\begin{pmatrix} a_{11} & a_{12}\\ 0 & a_{22}
\end{pmatrix}^{-1} =
\begin{pmatrix} a_{11}^{-1} & -a_{11}^{-1}a_{12}a_{22}^{-1}\\ 0 &
a_{22}^{-1}
\end{pmatrix}
\end{equation}
where $a_{22}^{-1} = \sum_{i\geq 0} (-1)^i
(a_0^{-1}a_-)^ia_0^{-1}$ (notice that $a_0^{-1}a_-$ is a locally
nilpotent map on the vector space $K[x]/K[x]_{\leq d}$).
\end{enumerate}
\end{theorem}

{\it Proof}. 1. Suppose that $a_{K[x]}\in \Aut_K(K[x])$. Then
$\ind_{K[x]}(a)=0$. By Proposition \ref{a12Jun10}.(1),
$a=\sum_{i\geq 1}a_{-i}\der^i+a_0+\sum \l_{ij}e_{ij}$ with
$a_0\neq 0$. Since $K[x]_{\leq d}$ is an $a$-{\em invariant finite
dimensional} subspace, we must have $b\in \GL(K[x]_{\leq d})$. The
vector space $K[x]=\bigcup_{s\geq d}K[x]_{\leq s}$ is the union of
$a$-invariant finite dimensional subspace $K[x]_{\leq s}$ and the
element $a$ acts on the one-dimensional  factor space $K[x]_{\leq
s}/K[x]_{\leq s-1}$ by scalar multiplication   $a_0(s+1)$ (for all
$s>d$) since $a_0(H)*x^s= a_0(s+1)x^s$. Now, statement 1 is
obvious.

2. Statement 2 is obvious (see statement 1).

3. For an arbitrary element $a\in \mI_1\backslash F$ with
$\ind_{K[x]}(a) =0$,  we have seen in the proof of statement 1
that  each term of the finite dimensional filtration $\{
K[x]_{\leq s}\}_{s\geq d}$ of the $\mI_1$-module $K[x]$ is an
$a$-invariant subspace. For a linear endomorphism  acting in a finite
dimensional vector space the conditions of being an injection, a
surjection, or a bijection are equivalent. Therefore, statement 3
follows.

4. Clearly, (\ref{am1b}) holds and $a_{22}^{-1} = (a_0+a_-)^{-1}=
(a_0(1+a_0^{-1}a_-))^{-1}=(1+a_0^{-1}a_-)^{-1}a_0^{-1}=\sum_{i\geq
0}(-1)^i(a_0^{-1} a_-)^ia_0^{-1}$.  $\Box $

$\noindent $

{\bf Classification of the elements $a\in \mI_1$ such that $a_{K[x]}$
is a surjection}.

\begin{theorem}\label{15Jun10}
 Let $a\in \mI_1$. Then the map
$a_{K[x]}$ is surjective iff $a=\sum_{i\geq n}a_{-i}\der^i
+a_F=a'\der^n+a_F$ where $n\geq 0$, $a_{-n}\neq 0$, all $a_{-i}\in
K[H]$, $a_F\in F$, $a':= \sum_{i\geq n}a_{-i}\der^{i-n}$, and
\begin{enumerate}
\item if $a_F=0$ then none of the natural numbers $j\geq 1$ is a
root of the polynomial $a_{-n}\in K[H]$; in this case, $a_{K[x]}'$
is a bijection and
$\ker(a_{K[x]})=\ker(\der^n_{K[x]})=\bigoplus_{i=0}^{n-1}Kx^i$;
and \item if $a_F\neq 0$ then none of the natural numbers $j\geq
d+2$ (where $d=\deg_F(a)$) is a root of the polynomial $a_{-n}\in
K[H]$ and $\im (a_{K[x]_{\leq d}})+\sum_{j=0}^{\min \{ d,n\} }
Ka'*x^{d-j}=K[x]_{\leq d}$; in this case there is  the  short
exact sequence of finite dimensional vector spaces
$$0\ra \ker
(a_{K[x]_{\leq d}})\ra \ker (a_{K[x]})\ra \ker (\d )\ra 0$$ where
$\d = a'\der^n: \bigoplus_{i=d+1}^{d+n}Kx^i\ra K[x]_{\leq
d}/aK[x]_{\leq d}$, and therefore $$\ker(a_{K[x]})=\ker
(a_{K[x]_{\leq d}})\bigoplus \bigoplus_{i=1}^s K(v_i-u_i)$$ where
$\{ v_1, \ldots , v_s\}$ is a basis for the kernel $\ker (\d )$ of
$\d$ and $u_1, \ldots , u_s\in K[x]_{\leq d}$ be any elements such
that $a'\der^n*v_i = a*u_i$.
\end{enumerate}
If the map $a_{K[x]}$ is surjective then $n:=
\ind_{K[x]}(a)=\dim_K(\ker(a_{K[x]})\geq 0$.
\end{theorem}

{\it Proof}. Suppose that the map $a_{K[x]}$ is a surjection. Then
$n:=\ind_{K[x]}(a) = \dim_K(\ker(a_{K[x]}))\geq 0$ and, by
Proposition \ref{a12Jun10}.(1), $a=\sum_{i\geq n}a_{-i}\der^i+a_F$
is the canonical form of the element $a$ where $a_{-n}\neq 0$, all
$a_j\in K[H]$, and $a_F\in F$. We can write $a=a'\der^n+a_F$ where
$a:=a_{-n}+a_{-n-1}\der+\cdots$ is an element of the set
$\mI_1\backslash F$ with $\ind_{K[x]}(a')=0$, by Proposition
\ref{a12Jun10}.(1).

Suppose that $a_F=0$, i.e.,  $a=a'\der^n$. Then $\ker_{K[x]}(a)
\supseteq \ker_{K[x]}(\der^n)= \langle 1,x,\ldots ,
x^{n-1}\rangle$, and so $\ker_{K[x]}(a) = \ker_{K[x]}(\der^n)$
since $\dim_K(\ker_{K[x]}(a))=n$. The map $a_{K[x]}'$ must be
injective: if $a'p=0$ for some polynomial $0\neq p\in K[x]$ then
$a'\der^n\int^np=0$; on the one hand $\int^n p\in \ker_{K[x]}(a)$
but on the other hand $\int^n p\not\in \langle 1,x,\ldots ,
x^{n-1}\rangle =\ker_{K[x]}(a)$ since $\deg_x(\int^n p )\geq n$, a
contradiction. By Theorem \ref{12Jun10}.(3), the map $a_{K[x]}'$
is a bijection since $\ind_{K[x]}(a')=0$. By Theorem
\ref{12Jun10}.(1), none of the natural numbers $j\geq 1$ is a root
of the polynomial $a_{-n}$. Conversely, suppose that $a=a'\der^n$
and none of the natural numbers $j\geq 1$ is a root of the
polynomial $a_{-n}$. By Theorem \ref{12Jun10}.(1), the map
$a_{K[x]}'$ is a bijection, then the map $a_{K[x]}$ is a
surjection. This finishes the proof of the theorem in the case
when $a_F=0$.

 Suppose that $a_F\neq 0$. Recall that  $d=\deg_F(a)$. Suppose, for a
moment, that the map $a_{K[x]}$ is not necessarily a surjection.
Applying the Snake Lemma to the commutative diagram of the short
exact sequence of vector spaces (where $V=K[x]_{\leq d}$ and
$U=K[x]/V$)
$$
\xymatrix{0\ar[r] & V\ar[r]\ar[d]^{ a_V}  & K[x] \ar[r]\ar[d]^{ a_{K[x]}} & U \ar[r]\ar[d]^{a_U} & 0 \\
0\ar[r] & V\ar[r]  & K[x] \ar[r] & U \ar[r] & 0 }
$$
we obtain the long exact sequence of vector spaces
\begin{equation}\label{VUa}
0\ra \ker(a_V)\ra \ker (a_{K[x]})\ra \ker (a_U)\stackrel{\d}{\ra }
\coker(a_V)\ra \coker (a_{K[x]})\ra \coker (a_U)\ra 0.
\end{equation}
Therefore, the map $a_{K[x]}$ is surjective iff the maps $\d$ and
$a_U$ are surjective. Clearly, for the element $a=a'\der^n+a_F$,
$a_U=(a'\der^n)_U$. Using the same argument as in the case (a) we
conclude that the map $a_U$ is surjective iff $a'|_U$ is an
isomorphism iff $a_{-n}|_U$ is an isomorphism iff $a_{-n}(H)*x^i =
a_{-n}(i+1) x^i \neq 0$ for $i\geq d+1$ iff none of the natural
numbers $i\geq d+2$ is a root of the polynomial $a_{-n}(H)$; and
in this case $\ker(a_U) = \ker(\der^n_U) =
\bigoplus_{i=d+1}^{d+n}Kx^i$. Let us give more details in the
proof of the only non-obvious step above: `if $a_U$ is surjective
then $a'|_U$ is an isomorphism.' Notice that the vector space $U$
is invariant under the action of the elements $a'$ and $\der^n$,
$\ker_U(\der^n) = \bigoplus_{i=d+1}^{d+n}Kx^i\subseteq
\ker_U(a'\der^n)$, and $ \der|_U$ is a surjection. Since
$$ \ker_U(a'\der^n) = \ind_U(a'\der^n) = \ind_U(a')
+\ind_U(\der^n) = 0+n=n,$$ we must have $\ker_U(a'\der^n) =
\ker_U(\der^n)$. The map $a'_U$ must be injective; if $a'u=0$ for
some element $0\neq u\in U$ then $a'\der^n \int^n u=0$; on the one
hand $\int^n u \in \ker_{K[x]}=\bigoplus_{i=d+1}^{d+n} Kx^i$ but
on the other hand $\int^n u \not\in \langle x^{d+1}, x^{d+2},
\ldots , x^{d+n}\rangle$, a contradiction. Therefore, $a'_U$ is an
injection, and so $a'_U$ is a bijection. Notice that,  for all
$i=d+1, \ldots , d+n$,
 $$\d (x^i) = a*x^i +\im
(a_V) =\begin{cases}
i(i-1)\cdots (i-n+1)a'*x^{i-n}& \text{if }i\geq n,\\
0& \text{otherwise}.\\
\end{cases}$$
 Therefore, the map $\d$ is
surjective iff $\im (a_V)+\sum_{j=0}^{\min \{ d,n\} }
Ka'*x^{d-j}=V$.  This finishes the proof of `iff' part of  case 2
of the theorem.

To finish the proof of case 2 suppose that the map $a_{K[x]}$ is
surjective. By (\ref{VUa}), there is the short exact sequence of
vector spaces $0\ra \ker (a_V)\ra \ker (a_{K[x]})\ra \ker (\d )\ra
0$ where $\d = a'\der^n: \bigoplus_{i=d+1}^{d+n}Kx^i\ra V/aV$.
Then $\ker(a_{K[x]})=\ker (a_V) \bigoplus \bigoplus_{i=1}^s
K(v_i-u_i)$ where $\{ v_1, \ldots , v_s\}$ is a basis for the
kernel $\ker (\d )$ of $\d$ and $u_1, \ldots , u_s\in V$ be any
elements such that $a'\der^n*v_i = a*u_i$. $\Box $

\begin{corollary}\label{a25Sep10}
Let an element  $a\in \mI_1$ be such that the map $a_{K[x]}$ is
surjective. Then the element $a$ can be uniquely written as the
product $a=\der^na'$ where $n=\dim_K(\ker (a_{K[x]}))$,
$a'=\sum_{i\geq 0}a_i'\der^i+f$, $a_0'\neq 0$, all $a_i'\in K[H]$,
$f=\sum_{i,j\geq 0} \l_{i+n, j}e_{i+n, j}$, $\l_{i+n,j}\in K$.
\end{corollary}

{\it Proof}. The existence and uniqueness follow from Theorem
\ref{15Jun10} and the equalities $\alpha (H) \cdot \der^n = \der^n
\cdot \alpha (H+n)$, for all elements $\alpha (H) \in K[H]$, and
$\der^n\cdot e_{ij} = e_{i-n, j}$.  $\Box $

$\noindent $


{\bf Classification of elements $a\in \mI_1$ such that $a_{K[x]}$
is an injection}.

\begin{lemma}\label{a16Jun10}
Let $a=\sum_{i\geq 1} a_{-i} \der^i +a_0+f\in \mI_1$ where all
$a_j\in K[H]$, $a_0\neq 0$, and $f\in F$. Then $\ker_{K[x]}(a\cdot
) = \ker_{K[x]_{\leq m}}(a\cdot )$ and $\coker_{K[x]}(a\cdot )
\simeq \coker_{K[x]_{\leq m}}(a\cdot )$ where $m=m(a) = \max \{
\deg_F(a), s\}$ and $s=\max \{ j\in \N \, | \, j+1$ is a root of
the polynomial $a_0\in K[H]\}$ or $s:=-1$ if there is no such a
root $j+1$.
\end{lemma}

{\it Proof}. Notice that $K[x]_{\leq m}$ is a $K[a]$-submodule of
$K[x]$ such that the map $a|_{K[x]/K[x]_{\leq m}}$ is a bijection
by the choice of $m$. Then (\ref{VUa}) yields the result. $\Box $

$\noindent $

Let $a\in \mI_1$ be such that $n=\max \{ i>0\, | \, a_i\neq
0\}\neq 0$, see (\ref{acan}). By Proposition \ref{a12Jun10}.(1),
$n=-\ind_{K[x]}(a)>0$. Then the element $a$ is a {\em unique} sum
(Proposition \ref{a12Jun10}.(2b)) 
\begin{equation}\label{acan3}
a=(\sum_{i\geq 1}b_{-i}\der^i +b_0+f)\int^n
\end{equation}
where all $b_j\in K[H]$, $b_0\neq 0$, $f\in \sum_{i,j\geq 0}
Ke_{i,j+n}$.

\begin{theorem}\label{16Jun10}
Let $a\in \mI_1$. Then the map $a_{K[x]}$ is an injection iff $a=
(\sum_{i\geq 1}a_{-i}\der^i +a_0+f)\int^n$ where $n\geq 0$, all
$a_j\in K[H]$, $a_0\neq 0$, $f\in \sum_{i,j\geq 0}Ke_{i,j+n}$ (the
presentation for $a$ is unique), and
\begin{enumerate}
\item when $n=0$ the map  $a_{K[x]}$ is a bijection (see Theorem
\ref{12Jun10}.(1)); 
 and
\item
 when $n\geq 1$,  $\ker_{K[x]_{\leq m}}(a'\cdot ) \bigcap x^n
 K[x]=0$ where $a':= \sum_{i\geq 1} a_{-i}\der^i +a_0+f$ and $m =
 m(a')$ is as in Lemma \ref{a16Jun10}; in this case there is the
 short exact sequence of finite dimensional vector spaces
 $$ 0\ra \ker_{K[x]_{\leq m}}(a'\cdot ) \ra K[x]/(x^n)\ra \coker
 (a_{K[x]})\ra K[x]_{\leq m} /a'*K[x]_{\leq m}\ra 0.$$
\end{enumerate}
If the map $a'_{K[H]}$ is a bijection then  $n=-\ind_{K[x]}(a) =
\dim_K(\coker (a_{K[x]}))$.
\end{theorem}

{\it Proof}. Suppose that the map $a_{K[x]}$ is  injective. Then
$n:=-\ind_{K[x]}(a) = \dim_K(\coker (a_{K[x]}))\geq 0$ and, by
Proposition \ref{a12Jun10}.(1) and (\ref{acan3}), the element $a$
is the unique sum $a=(\sum_{i\geq 1}a_{-i}\der^i +a_0+f)\int^n =
a'\int^n$ where all $a_j\in K[H]$, $a_0\neq 0$ and $f\in
\sum_{i,j\in \N} Ke_{i,j+n}$. To finish the proof of the theorem
we have to show that for each element $a$ that admits such a
presentation, $a=a'\int^n$ for some $n\in \N$, the map $a_{K[x]}$
is an injection iff the conditions 1 and 2 hold. If $n=0$ then
this follows from Theorem \ref{12Jun10}.(1,3). Suppose that $n\geq
1$. Then the map $a_{K[x]}$ is an injection iff
$\ker_{K[x]}(a'\cdot )\cap \im_{K[x]}(\int^n \cdot )=0$ iff
$\ker_{K[x]_{\leq m}}(a'\cdot ) \cap (x^n) =0$, by Lemma
\ref{a16Jun10}. This completes the `iff' part of statement 2. It
remains to prove  existence of the long exact sequence in
statement 2. So, let $n\geq 1$ and the map $a_{K[x]}$ is
injective. Applying (\ref{ables}) for the product $a=a'\int^n$,
yields the long exact sequence of vector spaces
 $$ 0\ra \ker_{K[x]}(a'\cdot ) \ra \coker_{K[x]}(\int^n)\ra \coker
 (a_{K[x]})\ra \coker_{K[x]} (a'\cdot ) \ra 0.$$
Notice that $\ker_{K[x]} (a'\cdot )=\ker_{K[x]_{\leq m}} (a'\cdot
)$,   $\coker_{K[x]} (a'\cdot )\simeq \coker_{K[x]_{\leq m}} (a'\cdot )$
(Lemma \ref{a16Jun10}) and $\coker_{K[x]}(\int^n\cdot ) =
K[x]/(x^n)$. The proof is complete.  $\Box $

$\noindent $

In general, for a linear map acting in an infinite dimensional
space it is not easy to find its cokernel. The next several
results make this problem finite dimensional for
integro-differential operators.

\begin{corollary}\label{x15Jun10}
Let  an element $a\in \mI_1$ be such that the map $a_{K[x]}$ is an
injection with $n=\dim_K(\coker (a_{K[x]}))\geq 1$, and
$a=\int^na'+f$ be the  unique sum where $a'=\sum_{i\geq
1}a_{-i}\der^i +a_0$, all $a_j\in K[H]$, $a_0\neq 0$, and $f\in F$
(Theorem \ref{16Jun10}). Then the set $\{ 1,x,\ldots , x^{n-1}\}$
is a basis for $\coker (a_{K[x]})$ iff the map $(a'+\der^n
f)_{K[x]}$ is a bijection (see Theorem \ref{12Jun10}.(1) for the
classification of  bijections).
\end{corollary}

{\it Remark}. The existence and uniqueness of the presentation
$a=\int^na'+f$ follows from Theorem \ref{16Jun10} and
(\ref{Iidi}).

$\noindent $

{\it Proof}. $(\Rightarrow )$  Suppose that the set  $\{
1,x,\ldots , x^{n-1}\}=\ker_{K[x]}(\der^n \cdot )$ is a basis for
the cokernel of the map  $a_{K[x]}$. Then $K[x] =
\ker_{K[x]}(\der^n \cdot ) \bigoplus \im (a_{K[x]})$, and so the
map $(\der^n a)_{K[x]}= (a'+\der^n f)_{K[x]}$ is an injection,
hence it is a bijection, by Theorem \ref{12Jun10}.(3).

$(\Leftarrow )$ Suppose that the map $(\der^n a)_{K[x]}=
(a'+\der^n f)_{K[x]}$ is a bijection. Then $\ker_{K[x]}(\der^n
\cdot ) \cap \im (a_{K[x]})=0$, and so $\ker_{K[x]}(\der^n \cdot )
\bigoplus \im (a_{K[x]})=K[x]$ since $\dim_K(\coker (a_{K[x]}))=n$
and $ n=\dim_K(\ker_{K[x]}(\der^n\cdot ))$. Therefore, the set $\{
1,x,\ldots , x^{n-1}\}$ is a basis for $\coker (a_{K[x]})$. $\Box
$


\begin{proposition}\label{a17Jun10}
Let $V$ be a nonzero finite dimensional  subspace of $K[x]$ of
dimension $n$. Then there exists a unit $s\in (1+F)^*$ such that
$V=\ker (s\der^ns^{-1})_{K[x]}$ and $K[x]= V\bigoplus \im
(s\int^ns^{-1})_{K[x]}$, i.e.,  $V\simeq \coker
(s\int^ns^{-1})_{K[x]}$. In particular, $s\der^ns^{-1}=\der^n+g$
and $s\int^ns^{-1} = \int^n +f$ for some elements $g,f\in F$.
\end{proposition}

{\it Proof}. If $V=\langle  1,x,\ldots , x^{n-1}\rangle =\ker
(\der^n)_{K[x]}$ then take $s=1$ as $K[x]= \ker
(\der^n)_{K[x]}\bigoplus \im (\int^n)_{K[x]}$. In the general
case, fix a natural number $m\geq n$ and  subspaces
$U,V,W\subseteq K[x]_{\leq m}$ such that $K[x]_{\leq m} =
V\bigoplus U = \ker (\der^n)_{K[x]} \bigoplus W$. Since
$(1+F)^*=\GL_\infty (K)$, we can find an element $s$ of $(1+F)^*$
such that $s^{-1} (V) = \ker (\der^n)_{K[x]}$, $s^{-1} (U) = W$
and $s(u) = u$ for all elements $u\in K[x]_{>m}:=\bigoplus_{i>m}
Kx^i$. Then the element $s$ satisfies the conditions of the
proposition:
$$ K[x]=  s\, \ker(\der^n)_{K[x]}\bigoplus s\,  \im
(\int^n)_{K[x]}= \ker (s\der^n s^{-1})_{K[x]}\bigoplus \im
(s\int^ns^{-1}) = V\bigoplus \im (s\int^ns^{-1}).$$ In particular,
$s=1+h$ and $s^{-1}=1+h'$ for some elements $h,h'\in F$, and so
$s\der^n s^{-1} = \der^n +g$ and $s\int^n s^{-1} = \int^n +f$ for
some elements $g,f\in F$. $\Box $


\begin{corollary}\label{b17Jun10}
Let an element $a\in \mI_1$ be such that the map $a_{K[x]}$ is an
injection with $n=\dim_K(\coker (a_{K[x]}))\geq 1$; then $a=\int^n
a'+f$ is the unique sum where $a'=\sum_{i\geq 1} a_{-i}\der^i
+a_0$, all $a_j\in K[H]$, $a_0\neq 0$, and $f\in F$ (Theorem
\ref{16Jun10}). Let $g\in F$ be such that the map $(\der^n
+g)_{K[x]}$ is a surjection. Then $\ker_{K[x]} (\der^n +g) \simeq
\coker (a_{K[x]})$ (i.e.,  $K[x] = \ker_{K[x]}(\der^n +g) \bigoplus
\im (a_{K[x]})$) iff the map $((\der^n+g)a)_{K[x]}= (a'+h)_{K[x]}$
is a bijection (see Theorem \ref{12Jun10}.(1) for the
classification  of  bijections) where $h:=\der^n f+ga\in F$. 
\end{corollary}

{\it Proof}. $(\Rightarrow )$ Suppose that $K[x]=
\ker_{K[x]}(\der^n +g) \bigoplus \im (a_{K[x]})$. Then the map
 $((\der^n+g)a)_{K[x]}= (a'+h)_{K[x]}$ is a injection, hence it is
 a bijection, by Theorem \ref{12Jun10}.(3).

$(\Leftarrow )$ Suppose that the map $((\der^n+g)a)_{K[x]}=
(a'+h)_{K[x]}$ is a bijection. Then $\ker_{K[x]}(\der^n+g)\cap \im
(a_{K[x]})= 0$, and so $\ker_{K[x]}(\der^n+g)\bigoplus  \im
(a_{K[x]})=K[x]$ since $n=-\deg_{\der^{-1}}(\der^n +g+F) =
\ind_{K[x]}(\der^n +g) = \dim_K(\ker_{K[x]}(\der^n +g))$, by
Proposition \ref{a12Jun10}.(1),  and $n= \dim_K(\coker
(a_{K[x]}))$.
 $\Box $

$\noindent $

Corollary \ref{b17Jun10} is an effective tool in finding a basis
for the cokernel of an injection $a_{K[x]}$.

$\noindent $

{\it Example}. Let $a= \der +\int$. Then $a_{K[x]}$ is an
injection and the map $\der_{K[x]}$ is a surjection such that
$(\der a)_{K[x]} =(\der ^2+1)_{K[x]}\in \Aut_K(K[x])$ is a
bijection since $(1+\der^2)^{-1}_{K[x]}= \sum_{i\geq 0} (-1)^i
\der^{2i}_{K[x]}$. By Corollary \ref{b17Jun10},
$\coker_{K[x]}(\der +\int ) \simeq \ker_{K[x]}(\der )= K$. Notice
that $(1+\der^2)^{-1} \not\in \mI_1^*$ where $\mI_1^*$ is the
group of  units of the algebra $\mI_1^*$.

\begin{lemma}\label{a18Jun10}
\begin{enumerate}
\item For each element $a\in \mI_1\backslash F$, there exists an
idempotent $f\in F$ such that $\ker (a_{K[x]}) = \im ( f_{K[x]})$.
\item Let $a\in \mI_1\backslash F$. Then there exists an element
$g\in F$ such that $\im (a_{K[x]})= \ker ( g_{K[x]})$ iff there
exists a natural number $d\geq 0$ such that $x^{d+1}K[x]\subseteq
\im (a_{K[x]})$ and $$\codim_{K[x]_{\leq d}}(K[x]_{\leq d}\bigcap
\im (a_{K[x]}))=\dim_K(\coker (a_{K[x]})).$$ In this case, the
element $g$ can be chosen to be an idempotent.
\end{enumerate}
\end{lemma}

{\it Proof}. 1. Since $a\in \mI_1\backslash F$, the kernel of the
linear map $a_{K[x]}$ is a finite dimensional vector space
(Theorem \ref{B30May10}.(1)), and so $\ker (a_{K[x]})\subseteq
K[x]_{\leq m}$ for some natural number $m$. Then $K[x]_{\leq m}=
\ker (a_{K[x]}) \bigoplus V$ for some subspace $V$ of $K[x]_{\leq
m}$. If $f$ is the projection onto the direct summand $\ker
(a_{K[x]})$ of $K[x]_{\leq m}$ extended by zero on $K[x]_{>m}$
then $f^2=f$, $f\in F$ and $\ker (a_{K[x]}) = \im ( f_{K[x]})$.

2. $(\Rightarrow )$ Suppose that there exists an element $g\in F$
such that $\im (a_{K[x]})= \ker ( g_{K[x]})$. Let $d=\deg_F(g)$.
Then $K[x]= K[x]_{\leq d} \bigoplus (x^{d+1})$ and $ \ker
(g_{K[x]}) = \ker (g_{K[x]_{\leq d}}) \bigoplus (x^{d+1})$.
Therefore $(x^{d+1}) \subseteq \im ( a_{K[x]})$ and
$\codim_{K[x]_{\leq d}} (K[x]_{\leq d}\bigcap \im
(a_{K[x]}))=\codim_{K[x]_{\leq d}}(\ker (g_{K[x]_{\leq
d}}))=\codim_{K[x]}(\ker (g_{K[x]}))=  \dim_K( \coker (a_{K[x]}))$.

$(\Leftarrow )$ Suppose that $(x^{d+1}) \subseteq \im ( a_{K[x]})$
and $$\codim_{K[x]_{\leq d}} (K[x]_{\leq d}\bigcap \im
(a_{K[x]}))= \dim_K( \coker (a_{K[x]}))=:n.$$ Then $K[x]_{\leq d}
= V\bigoplus (K[x]_{\leq d} \bigcap \im ( a_{K[x]}))$ for some
subspace $V$ of $K[x]_{\leq d}$ with $\dim_K(V) = n$. Then $\im
(a_{K[x]}) = K[x]_{\leq d} \bigcap \im ( a_{K[x]})\bigoplus
(x^{d+1})$. It suffices to take $g$ which is the projection map
$K[x]_{\leq d} \ra V$ extended by zero on the ideal $(x^{d+1})$.
$\Box $

$\noindent $

{\em Non-Example}. The conditions of statement 2 are very
restrictive. In particular, not for every element $a\in
\mI_1\backslash F$ there exists an element $g$ such that $\im (
a_{K[x]}) = \ker (g_{K[x]})$, eg, $a= \der +\int$ since
$x^{2n}\not\in \im (\der +\int )$ for all $n\geq 0$.

{\it Proof}. We have to show that there is no polynomial solution
$u$ to the equation $(\der +\int ) *u=x^{[2n]}$. Note that
$$ \der +\int = \der \int (\der +\int ) = \der (1-e_{00}+\int^2) =
\der (1-e_{00}) (1+\int^2).$$ Then, $(1-e_{00})(1+\int^2)
*u=x^{[2n+1]}+C$ for some constant $C$ necessarily $C=0$ as $\im
(1-e_{00})_{K[x]}= (x)\ni x^{[2n+1]}$. We can write $u=\l +v$ for
some $\l \in K$ and $v\in (x)$. The linear maps $1-e_{00}$ and
$1+\int^2$acting in $K[x]$ respect the ideal $(x)$. By taking the
equality $(1-e_{00})(1+\int^2)*u=x^{[2n+1]}$ modulo $(x)$ yields
$(1-e_{00})\l =0$ and so $\l =0$. Since the $(1-e_{00})_{K[x]}$ is
the projection onto the ideal $(x)$ in the decomposition $K[x]=
K\bigoplus (x)$, $\l =0$ and $x^{[2n+1]}\in (x)$, we can drop
$1-e_{00}$ in the equation, i.e.,  $(1+\int^2) *u=x^{[2n+1]}$. The
only solution $u=\sum_{i\geq 0}(-1)^i\int^{2i}*x^{[2n+1]}\in
K[[x]]$ is obviously not a polynomial. $\Box$

\begin{proposition}\label{a19Jun10}
\begin{enumerate}
\item For each element $a\in \mI_1\backslash F$ with $n:=
\dim_K(\coker (a_{K[x]}))$, there exists an element $\der^n+f$ for
some $f\in F$ (resp. $s\in (1+F)^*$)   such that the map
$(\der^n+f)a_{K[x]}$ (resp. $s\der^n s^{-1}a_{K[x]}$) is a
surjection. In this case, $\ker((\der^n+f)a_{K[x]}) = \ker
(a_{K[x]})$ (resp. $\ker (s\der^ns^{-1}a_{K[x]})= \ker
(a_{K[x]}))$. \item For each element $a\in \mI_1\backslash F$ with
 $n:= \dim_K(\ker (a_{K[x]}))$, there exists an element $\int^n +g$
  for
some $g\in F$ (resp. $s\in (1+F)^*$) such that the map $a(\int^n
+g)_{K[x]}$ (resp. $as\int^n s^{-1}_{K[x]}$)  is an injection. In
this case, $\im (a(\int^n +g)_{K[x]})=\im ( a_{K[x]})$ (resp.
$\im (as\int^ns^{-1}_{K[x]})= \im (a_{K[x]}))$. \item For each
element $a\in \mI_1\backslash F$ with $m:= \dim_K(\ker
(a_{K[x]}))$ and $n :=\dim_K(\coker (a_{K[x]}))$,  there exist
elements $\der^n+f$ and $\int^m+g$ for some $f,g\in F$ (resp.
$s,t\in (1+F)^*$) such that the map $(\der^n+f) a(\int^m
+g)_{K[x]}$ (resp. $s\der^ns^{-1} at\int^mt^{-1}_{K[x]}$) is a
bijection.
\end{enumerate}
\end{proposition}

{\it Proof}. 1. Notice that $K[x]= V\bigoplus \im (a_{K[x]})$ for
some $n$-dimensional subspace $V$ of $K[x]$. By Proposition
\ref{a17Jun10}, there exists a unit $s\in (1+F)^*$ such that
$V=\ker (s\der^ns^{-1}_{K[x]})$ and $s\der^ns^{-1} = \der^n+f$ for
some $f\in F$. Then the map $(s\der^ns^{-1})a_{K[x]}$ is
surjective since $K[x]= s\der^ns^{-1} *K[x]= s\der^ns^{-1}
(\ker(s\der^ns^{-1}) \bigoplus \im (a_{K[x]}))= \im
((s\der^ns^{-1})a_{K[x]})$. Then, by Lemma \ref{D30May10},
\begin{eqnarray*}
\dim_K \ker ((\der^n+f)a_{K[x]})&=&\ind_{K[x]}((\der^n+f)a) =
\ind_{K[x]}(\der^n+f) +\ind_{K[x]}(a)\\
& = &\ind_{K[x]}(\der^n) +\ind_{K[x]}(a) = n+\dim_K \ker
(a_{K[x]})-n= \dim_K \ker (a_{K[x]}).
\end{eqnarray*}
Therefore, $\ker((\der^n+f)a_{K[x]}) = \ker(a_{K[x]})$ since $\ker
((\der^n+f)a_{K[x]}) \supseteq \ker (a_{K[x]})$.

2.  Since $a\not\in F$, the kernel $V$ of the linear map
$a_{K[x]}$ is finite dimensional (Theorem \ref{B30May10}). By
Proposition \ref{a17Jun10}, there exists a unit $s\in (1+F)^*$
such that $K[x]= V\bigoplus \im (s\int^n s^{-1})_{K[x]}$ and
$s\int^ns^{-1} = \int^n+g$ for some element $g\in F$. It follows
that the map $a(\int^n +g)_{K[x]}$ is an injection. Then, by Lemma
\ref{D30May10},
\begin{eqnarray*}
-\dim_K \coker (a(\int^n+g)_{K[x]})&=&\ind_{K[x]}(a(\int^n+g)) =
\ind_{K[x]}(a)+ \ind_{K[x]}(\int^n+g) \\
& = &\ind_{K[x]}(a)+\ind_{K[x]}(\int^n)  = n-\dim_K \coker
(a_{K[x]})-n\\
&=& -\dim_K \coker (a_{K[x]}).
\end{eqnarray*}
Therefore, the natural inclusion  $\im
(a(\int^n+g)_{K[x]})\subseteq  \im (a_{K[x]})$ is an equality.

3. By statement 1, there exists an element $\der^n+f$ for some
$f\in F$ (resp. $s\in (1+F)^*$) such that the map $a':= (\der^n
+f)a_{K[x]}$ (resp. $a':= s\der^ns^{-1}a_{K[x]}$) is a surjection
with $\ker(a'_{K[x]}) = \ker(a_{K[x]})$. Then, by statement 2, for
the element $a'$,  there exists an element
$\int^m+g$ for some $g\in F$ (resp. $t\in (1+F)^*$) such that  the
map $a'':= a' (\int^m+g)_{K[x]}$ (resp. $a'':=
a't\int^mt^{-1}_{K[x]}$) is an injection with $\im
(a''_{K[x]})=\im (a'_{K[x]})=K[x]$, i.e.,  the map $a''_{K[x]}$ is a
bijection. $\Box $

$\noindent $

{\it Example}. Let $ a= \der +\int$. We know already that
$a_{K[x]}$ is an injection with $\dim_K(\coker (a_{K[x]}))=1$.
Then $\der a=\der^2+1$ and $ (1+\der^2)_{K[x]}$ is a bijection.

$\noindent $

The next proposition is useful in proving that various Ext's and
Tor's are finite dimensional vector spaces or not.

\begin{proposition}\label{a6Jun10}
Let $a,b\in \mI_1$.
\begin{enumerate}
\item If $a,b\not\in F$ then the vector spaces
$\ker_{\ker_{\mI_1}(\cdot b)}(a\cdot )$,
$\coker_{\ker_{\mI_1}(\cdot b)}(a\cdot )$,
 $\ker_{\coker_{\mI_1}(\cdot b)}(a\cdot )$,

 $ \coker_{\coker_{\mI_1}(\cdot b)}(a\cdot
 )$ are finite dimensional.
\item If $a\not\in F$ and $b\in F$ then
\begin{enumerate}
\item $\ker_{\ker_F(\cdot b)}(a\cdot )=\ker_{\ker_{\mI_1}(\cdot
b)}(a\cdot )$ and $\dim_K ( \ker_{\ker_{\mI_1}(\cdot b)}(a\cdot
))=\begin{cases}
0& \text{if } \ker_{K[x]}(a\cdot )=0,\\
\infty & \text{otherwise}.\\
\end{cases}$
\item The sequence $0\ra \coker_{\ker_F(\cdot b)}(a\cdot )\ra
\coker_{\ker_{\mI_1}(\cdot b)}(a\cdot )\ra \coker_{B_1}(a\cdot )
\ra 0$ is exact and $\dim_K (\coker_{\ker_{\mI_1}(\cdot b)}(a\cdot
))=\begin{cases}
0& \text{if } a=\l\der^i+f:K[x]\ra K[x] \text{ is a surjection},\\
\infty & \text{otherwise},\\
\end{cases}$

for some $\l \in K^*$, $i\geq 0$ and $f\in F$. \item
$\ker_{\coker_{\mI_1}(\cdot b)}(a\cdot ) \simeq
\ker_{\coker_F(\cdot b)}(a\cdot )$ and
$\dim_K(\ker_{\coker_{\mI_1}(\cdot b)}(a\cdot ))= \begin{cases}
0& \text{if }a_{K[x]} \, \text{is injective} ,\\
\infty & \text{otherwise}.\\
\end{cases}$
\item The sequence $0\ra \coker_{\coker_F(\cdot b)}(a\cdot ) \ra
\coker_{\coker_{\mI_1}(\cdot b)}(a\cdot ) \ra B_1/ aB_1\ra 0$ is
exact and $\dim_K(\coker_{\coker_{\mI_1}(\cdot b)}(a\cdot ))=
\begin{cases}
0& \text{if }a=\l\der^i+f, a_{K[x]} \, \text{is surjective} ,\\
\infty & \text{otherwise},\\
\end{cases}$

for some $\l \in K^*$, $i\geq 0$ and $ f\in F$.
\end{enumerate}

\item If $a\in F$ and $b\not\in F$ then
\begin{enumerate}
\item $\ker_{\ker_{\mI_1}(\cdot b)}(a\cdot )=\ker_{\ker_F(\cdot
b)}(a\cdot )$ and $\dim_K ( \ker_{\ker_{\mI_1}(\cdot b)}(a\cdot
))=\begin{cases}
\infty & \text{if } \ker_{\mI_1}(\cdot b)\neq 0,\\
0 & \text{otherwise}.\\
\end{cases}$

$\ker_{\mI_1} (\cdot b) =0\Leftrightarrow (\cdot b)_{\mI_1}$ is an
injection $\Leftrightarrow$ $\cdot b:\mI_1/ \int \mI_1\ra \mI_1/
\int \mI_1$ is an injection $\Leftrightarrow$ $b^*\cdot :K[x]\ra
K[x]$ is an injection.
 \item  $\coker_{\ker_{\mI_1}(\cdot b)}(a\cdot )\simeq \coker_{\ker_F(\cdot b)}(a\cdot
 )$ and
$\dim_K (\coker_{\ker_{\mI_1}(\cdot b)}(a\cdot ))=\begin{cases}
\infty& \text{if } \ker_{\mI_1}(\cdot b)\neq 0,\\
0 & \text{otherwise}.\\
\end{cases}$
 \item
$\dim_K(\ker_{\coker_{\mI_1}(\cdot b)}(a\cdot ))= \begin{cases}
0& \text{if } b=\l \int^i+f, (b^*\cdot )_{K[x]} \, \text{is surjective},\\
\infty & \text{otherwise},\\
\end{cases}$

where $0\neq \l \in K$, $i\geq 0$ and $f\in F$.  \item
$\dim_K(\coker_{\coker_{\mI_1}(\cdot b)}(a\cdot ))=
\begin{cases}
0& \text{if } b=\l \int^i+f, (b^*\cdot )_{K[x]} \, \text{is surjective},\\
\infty & \text{otherwise},\\
\end{cases}$

where $0\neq \l \in K$, $i\geq 0$ and $f\in F$.
\end{enumerate}
\item If $a,b\in F$ then the vector spaces
$\ker_{\ker_{\mI_1}(\cdot b)}(a\cdot )$,
$\coker_{\ker_{\mI_1}(\cdot b)}(a\cdot )$,
 $\ker_{\coker_{\mI_1}(\cdot b)}(a\cdot )$,

 $ \coker_{\coker_{\mI_1}(\cdot b)}(a\cdot
 )$ are infinite dimensional.
\end{enumerate}
\end{proposition}

{\it Remark}. Earlier, necessary and sufficient conditions were
given for the map $a\cdot :K[x]\ra K[x]$ to be bijective (Theorem
\ref{12Jun10}), surjective (Theorem \ref{15Jun10}) or injective
(Theorem \ref{16Jun10}).

$\noindent $

{\it Proof}. 1. Since $b\not\in F$, the left $\mI_1$-modules
$\ker_{\mI_1}(\cdot b)$ and $\coker_{\mI_1}(\cdot b)$ have finite
length (Theorem \ref{A30May10}). Then, by Theorem \ref{B30May10},
all four vector spaces are finite dimensional.

2(a,b). Since $b\in F$, $\ker_{B_1}(\cdot b) = \ker_{B_1}(0)=B_1$.
Since $a\not\in F$, $\ker_{\ker_{B_1}(\cdot b)}(a\cdot ) =
\ker_{B_1}(a\cdot ) =0$. Therefore, the long exact sequence in
Lemma \ref{a7Jun10}.(2a) brakes down into two short exact
sequences $\ker_{\ker_F(\cdot b)}(a\cdot
)=\ker_{\ker_{\mI_1}(\cdot b)}(a\cdot )$ and $0\ra
\coker_{\ker_F(\cdot b)}(a\cdot )\ra \coker_{\ker_{\mI_1}(\cdot
b)}(a\cdot )\ra \coker_{B_1}(a\cdot ) \ra 0$. Since $b\in F$,
${}_{\mI_1}\ker_F(\cdot b) \simeq K[x]^{(\N )}$, a direct sum of
countably many copies of the left $\mI_1$-module $K[x]$, then
statement $(a)$ follows. Using the short exact sequence, we see
that the vector space $\coker_{\ker_{\mI_1}(\cdot b)}(a\cdot )$ is
finite dimensional iff so are the vector spaces
$\coker_{\ker_F(\cdot b)}(a\cdot )$ and $\coker_{B_1}(a\cdot )$.
Since $\ker_F(\cdot b) \simeq K[x]^{(\N )}$, the vector space
$\coker_{\ker_F(\cdot b)}(a\cdot )$ is finite dimensional iff it
is a zero space  iff the map $a_{K[x]}$ is surjective. Since $
a\not\in F$, the vector space $\coker_{\ker_{B_1}(\cdot b)}(a\cdot
) = \coker_{B_1} (a\cdot )$ is finite dimensional iff $a+F$ is a
unit of the algebra $B_1$ iff $a+F= \l\der^i\in B_1$ where $\l \in
K^*$ and $i\in \Z$ iff $\coker_{B_1} (a\cdot )=0$. By Theorem
\ref{15Jun10}, the vector space $\coker_{\ker_{\mI_1}(\cdot
b)}(a\cdot )$ is finite dimensional iff it is a zero space iff
$a=\l \der^i+f$ for some $\l \in K^*$, $i\geq 0$, and $f\in F$
such that the map $(\l \der^i+f)_{K[x]}$ is surjective.

2(c,d). Since $a\not\in F$ and $b\in F$, we have
$\coker_{B_1}(\cdot b) = B_1$ and $\ker_{\coker_{B_1}(\cdot
b)}(a\cdot  )=0$. Therefore, the long exact sequence in Lemma
\ref{a7Jun10}.(2b) collapses to $\ker_{\coker_F(\cdot b)}(a\cdot
)\simeq U:=\ker_{\coker_{\mI_1}(\cdot b)}(a\cdot )$ and the short
exact sequence
$$ 0\ra V_1:= \coker_{\coker_F(\cdot b)}(a\cdot ) \ra V:=
\coker_{\coker_{\mI_1}(\cdot b)}(a\cdot ) )\ra B_1/aB_1\ra 0.$$
Since ${}_{\mI_1}F_{\mI_1}\simeq K[x]\t (\mI_1 / \int \mI_1 )$ and
$b\in F$, ${}_{\mI_1}\coker_F(\cdot b) \simeq K[x]^{(\N )}$.
Therefore, $$\dim_K(U)= \begin{cases}
0& \text{if } a_{K[x]}\,  \text{is an injection},\\
\infty& \text{otherwise}.\\
\end{cases}$$
Notice that $\dim_K(V)=\dim_K(V_1)+\dim_K(B_1/
aB_1)$.
$$\dim_K(B_1/ aB_1) =\begin{cases}
0& \text{if } a=\l\der^i+f, \l\int^i+f,\\
\infty& \text{otherwise},\\
\end{cases}$$
where $0\neq \l \in K$, $i\geq 0$ and $f\in F$. $\dim_K(V_1)=
\begin{cases}
0& \text{if } a_{K[x]} \, \text{is surjective},\\
\infty& \text{otherwise}.\\
\end{cases}$

By Theorem \ref{15Jun10}, $\dim_K(V) =0$ iff $a=\l
\der^i+f$ and $a_{K[x]}$ is a surjection where $0\neq \l \in K$,
$i\geq 0$ and $f\in F$; otherwise $\dim_K(V) = \infty$.

3(a,b). Notice that $\ker_{\mI_1}(\cdot b) =0 \Leftrightarrow
(\cdot b)_{\mI_1}$ is an injection $\Leftrightarrow$ $(\cdot b)_F$
is an injection (since $b\not\in F$) $\Leftrightarrow$  $\cdot
b:\mI_1/ \int \mI_1\ra \mI_1/ \int \mI_1$ is an injection (since
$F_{\mI_1}\simeq (\mI_1/ \int \mI_1)^{(\N )})$ $\Leftrightarrow$
$b^*\cdot :K[x]\ra K[x]$ is an injection since
${}_{\mI_1}K[x]\simeq \mI_1/\mI_1 \der$ and the map $*:\mI_1/ \int
\mI_1 \ra \mI_1/ \mI_1\der $, $u+\int\mI_1\mapsto u^*+\mI_1\der$,
is a bijection such that $(cu)^* = u^*c^*$ for all elements $c\in
\mI_1$ and $u\in \mI_1/ \int\mI_1$.

Since $ b\not\in F$, we have $\ker_{B_1} (\cdot b) =0$, and so the
short exact sequence in Theorem \ref{a7Jun10}.(2a) yields
$\ker_{\ker_F(\cdot b)}(a\cdot )= \ker_{\ker_{\mI_1}(\cdot
b)}(a\cdot )$ and $\coker_{\ker_F(\cdot b)}(a\cdot )=
\coker_{\ker_{\mI_1}(\cdot b)}(a\cdot )$. Clearly, if
$\ker_{\mI_1}(\cdot b) =0$ then these vector spaces are equal to
zero. Suppose that $\ker_{\mI_1}(\cdot b)\neq 0$. Then
$\ker_{\mI_1}(\cdot b)=\ker_F(\cdot b)$ since $b\not\in F$, and
$\ker_{\mI_1}(\cdot b)\simeq K[x]^m$ for some $m\geq 1$, by
Theorem \ref{A30May10}.(1a). Since $a\in F$, the kernel and the
cokernel of the linear map $a\cdot : K[x]\ra K[x]$ are infinite
dimensional vector spaces, then so are the vector spaces
$\ker_{\ker_{\mI_1}(\cdot b)}(a\cdot )$ and
$\coker_{\ker_{\mI_1}(\cdot b)}(a\cdot )$.

3(c). Let $U:=\ker_{\coker_{\mI_1}(\cdot b)}(a\cdot )$. Suppose
that $\coker_F(\cdot b)\neq 0$. By Theorem \ref{B30May10}.(1a),
${}_{\mI_1}\coker_F(\cdot b)\simeq K[x]^m$ for some $m\geq 1$
(since $b\not\in F$), and so $\dim_K(\ker_{\coker_F(\cdot
b)}(a\cdot ))=\infty$ since $a\in F$. Since there is the inclusion
$\ker_{\coker_F(\cdot b)}(a\cdot )\subseteq U$ (Lemma
\ref{a7Jun10}.(2b)), we must have $\dim_K(U)=\infty$.

Suppose that $\coker_F(\cdot b)=0$. Then the long exact sequence
in Lemma \ref{a7Jun10}.(2b) yields as isomorphism $U\simeq
\ker_{\coker_{B_1}(\cdot b)}(a\cdot )=\coker_{B_1}(\cdot b)\simeq
B_1/B_1b$. Notice that $\coker_F(\cdot b) =0$ iff the map $\cdot
b: \mI_1/\int\mI_1\ra \mI_1/ \int\mI_1$ is surjective since
$F_{\mI_1}\simeq (\mI_1/\int \mI_1)^{(\N )}$ iff the map
$(b^*\cdot )_{K[x]}$ is surjective. $\dim_K(B_1/B_1b)$ is either
$0$ or $\infty$, and $\dim_K(B_1/ B_1b)=0$ iff $b=\l \int^i+f$,
$\l \der^i+f$ where $0\neq \l \in K$, $i\geq 0$ and $f\in F$. By
Theorem \ref{15Jun10}, the two conditions, $\coker_F(\cdot b)=0$
and $\dim_K(B_1/ B_1b)=0$, hold iff $b=\l \int^i+f$ such that
$(b^*\cdot )_{K[x]}$ is a surjection where $0\neq \l \in K$,
$i\geq 0$ and $f\in F$. These are precisely the conditions when
$U=0$, otherwise $\dim_K (U) = \infty$.

3(d). Let $V=\coker_{\coker_{\mI_1}(\cdot b)}(a\cdot )$. If
$\coker_{B_1} (\cdot b) \neq 0$, i.e.,  $\dim_K(B_1/ B_1b)=\infty$,
then the end of the long exact sequence in Lemma
\ref{a7Jun10}.(2b) yields the surjection
$V\ra\coker_{\coker_{B_1}(\cdot b)}(a\cdot )=\coker_{B_1}(\cdot b)
=B_1/ B_1b$ (since $a\in F$), and so $\dim_K(V)=\infty$.

If $\coker_{B_1}(\cdot b)=0$ then the long exact sequence in Lemma
\ref{a7Jun10}.(2b) yields an isomorphism of vector spaces
$\coker_{\coker_F(\cdot b)}(a\cdot )\simeq V$. Notice that the
condition $\coker_{B_1}(\cdot b) =0$ holds iff $b=\l\int^i+f$,
$\l\der^i+f$ where $0\neq \l \in K$, $i\geq 0$ and $f\in F$.

If $\coker_F(\cdot b)=0$ then $V=0$. The condition $\coker_F(\cdot
b)=0$ holds iff the map $\cdot b :\mI_1/ \int \mI_1\ra \mI_1/ \int
\mI_1$ is surjective since $F_{\mI_1} \simeq
(\mI_1/\int\mI_1)^{(\N )}$ iff the map $(b^*\cdot )_{K[x]}$ is
surjective. By Theorem \ref{15Jun10}, the two conditions,
$\coker_{B_1}(\cdot b)=0$ and $\coker_F(\cdot b)=0$, hold iff
$b=\l\int^i +f$ such that $(b^*\cdot )_{K[x]}$ is a surjection
where $0\neq \l \in K$, $i\geq 0$ and $f\in F$.

If $\coker_F(\cdot b)\neq 0$ then, by Theorem \ref{B30May10}.(1a),
${}_{\mI_1}\coker_F(\cdot b) \simeq K[x]^m$ for some $m\geq 1$
(since $b\not\in F$) and so $\dim_K(V)=\dim_K(\coker
(a\cdot)_{K[x]^m})=\infty$ since $a\in F$. This proves statement
3(d).

4. Since ${}_{\mI_1}F_{\mI_1}\simeq K[x]\t (\mI_1/\int\mI_1)$ and
$a,b\in F$, the vector spaces at the beginning and the end of each
left exact sequence  in Lemma \ref{a7Jun10}.(2a,b)  are infinite
dimensional, and so are the four vector spaces in statement 4.
$\Box $


\section{Classification of one-sided invertible elements of $\mI_1$}\label{COSI}

In this section, a classification of elements of the algebra
$\mI_1$ that admit a {\em one-sided} inverse is given (Corollary
\ref{c25Sep10}), an explicit description of all one-sided inverses
is found (Theorem \ref{A3Oct10}). It is proved that the monoid
$\CL (\mI_1)$ (respectively, $\CR (\mI_1)$)  of all elements  of
$\mI_1$ that admit a left (respectively, right) inverse is
generated by the group $\mI_1^*$ of units of the algebra $\mI_1$
and the element $\int$ (respectively, $\der $), Theorem
\ref{3Oct10}.

The algebra $K+F=K+M_\infty (K)$ has the obvious {\em determinant
map} $\det : K+F\ra K$, $a\mapsto \det (a)$. The element $a\in
K+F$ is the unique finite sum  $a=\l +\sum \l_{ij}e_{ij}=\l +\sum
\l_{ij}\frac{i!}{j!}E_{ij}$ where $\l , \l_{ij}\in K$. Then
\begin{equation}\label{detaKF}
\det (a) :=\begin{cases}
0& \text{if }\l =0\\
\l \det (a')& \text{if }\l \neq 0,\\
\end{cases}
=\begin{cases}
0& \text{if }\l =0\\
\l \det (a'')&
\text{if }\l \neq 0,\\
\end{cases}
\end{equation}
where $a':=\sum_{i\in \N} e_{ii}+\sum \l^{-1}\l_{ij}e_{ij}$, $a'':=\sum_{i\in \N} E_{ii}+\sum
\l^{-1}\l_{ij}\frac{i!}{j!}E_{ij}$ (notice that $a'=a''=\l^{-1}a$).  In the first equality the usual determinant is taken w.r.t.
the matrix units $\{ e_{ij}\}_{i,j\in \N}$ but in the second
equality -  w.r.t. the matrix units $\{ E_{ij}\}_{i,j\in \N}$.
Both determinants are equal since for the element $a\in K+F$,
$$[a]_e= S^{-1}[a]_ES,$$ where $[a]_e$ and $[a]_E$ are the matrix
forms of the element $a$ w.r.t. the matrix units $\{ e_{ij}\}$ and
$\{ E_{ij}\}$ respectively, and $S:=\sum_{i\in \N}
i!e_{ii}=\sum_{i\in \N} i!E_{ii}$ is the infinite diagonal
invertible matrix.

\begin{theorem}\label{25Sep10}
\begin{enumerate}
\item Let an element $a\in \mI_1$ be such that the map $a_{K[x]}$
is an injection,  $n:= \dim_K(\coker (a_{K[x]}))$, and $b\in \mI_1$.
\begin{enumerate}
\item Then  
 $ba\in \mI_1^*$ iff $a=a'\int^n$
where $a'\in K^*+F$ and  $b=\der^nb'$ where $b'\in K^* +F$ with
$\det (ba) \neq 0$ (notice that $ba \in K+F$).  \item If elements
$a=a'\int^n$, where $n\in \N$, and  $a'\in K^*+F$ are such that the
map $a_{K[x]}$ is injective (necessarily, $n=\dim_K(\coker
(a_{K[x]}))$, by Proposition \ref{a12Jun10}.(1)) then there is at least one
element $b$ such that in the  statement (a).
\end{enumerate}
\item Let an element $b\in \mI_1$ be such that the map $b_{K[x]}$
is a surjection,  $n:= \dim_K(\ker (b_{K[x]}))$, and $a\in \mI_1$.
\begin{enumerate}
\item Then 
 $ba\in \mI_1^*$ iff $b=\der^n b'$
where $b'\in K^*+F$ and  $a=a'\int^n$ where $a'\in K^* +F$ with
$\det (ba) \neq 0$.
 \item If elements $b=\der^n b'$, where $n\in \N$, and  $b'\in K^*+F$ are such that the map $b_{K[x]}$ is a surjection
 (necessarily, $n=\dim_K(\ker (b_{K[x]}))$, by Proposition
 \ref{a12Jun10}.(1))  then there is at
least one element $a$ such as in  statement (a).
\end{enumerate}
\end{enumerate}
\end{theorem}

{\it Proof}. 1(a). $(\Rightarrow )$ Suppose that there exists an
element $b\in \mI_1$ such that  $ba\in \mI_1^*$. Then necessarily
the map $b_{K[x]}$ is surjective.  Then taking the inclusion
$ba\in \mI_1^*$  modulo the ideal $F$ yields the inclusion
$\overline{b}\overline{a}\in K^*$ where $\overline{b}=b+F$ and
$\overline{a}=a+F$ since $\mI_1^* = K^* (1+F)^*$. Bearing in mind
that  $a_{K[x]}$ is an injection and $b_{K[x]}$ is a surjection,
we must have $a=a'\int^n$ where $a'\in K^*+F$ (Theorem
\ref{16Jun10}) and $b=\der^n b'$ where $b'\in K^*+F$ (Corollary
\ref{a25Sep10}). Since $ba\in \mI_1^* =K^*(1+F)^* = \{ c\in K+F\,
| \, \det (c)\neq 0\}$, we must have $\det (ba) \neq 0$.

$(\Leftarrow )$ Suppose that $a$ and $b$ satisfy the conditions
after `iff'. Then $\det (ba ) \neq 0$ implies that $ba \in
\mI_1^*$. 

1(b). By Theorem \ref{a19Jun10}.(1), for the element $a$ there
exists an element $b=\der^n +f$ with $ f\in F$ such that the map
$(ba)_{K[x]}$ is a surjection and $\ker (ba_{K[x]})= \ker
(a_{K[x]})=0$, i.e.,  $ba \in \Aut_K(K[x])\cap (K+F) = (K+F)^* =
\mI_1^*$. 

2(a). $(\Rightarrow )$ Suppose that there is an element $a\in
\mI_1$ such that  $ba\in \mI_1^*$. Then necessarily the map
$a_{K[x]}$ is injective.  Then taking the inclusion modulo the
ideal $F$ yields the inclusion $\overline{b}\overline{a}\in K^*$
where $\overline{b}=b+F$ and $\overline{a}=a+F$ since $\mI_1^* =
K^* (1+F)^*$. Bearing in mind that  $a_{K[x]}$ is an injection and
$b_{K[x]}$ is a surjection, we must have $a=a'\int^n$ where $a'\in
K^*+F$ (Theorem \ref{16Jun10}) and $b=\der^n b'$ where $b'\in
K^*+F$ (Corollary \ref{a25Sep10}). Since $ba\in \mI_1^*
=K^*(1+F)^* = \{ c\in K+F\, | \, \det (c)\neq 0\}$, we must have
$\det (ba) \neq 0$.

$(\Leftarrow )$ Suppose that $a$ and $b$ satisfy the conditions
after `iff'. Then $\det (ba ) \neq 0$ implies that $ba \in
\mI_1^*$. 

2(b). By Theorem \ref{a19Jun10}.(2), for the element $b$ there
is an element $a=\int^n +g$ with $ g\in F$ such that the map
$(ba)_{K[x]}$ is an injection and $\im (ba_{K[x]})= \im
(b_{K[x]})=K[x]$, i.e.,  $ba \in \Aut_K(K[x])\cap (K+F) = (K+F)^* =
\mI_1^*$. $\Box $  

$\noindent $

Let $\CL (\mI_1) :=\{ a\in \mI_1\, | \, ba=1$ for some $b\in
\mI_1\}$ and $\CR (\mI_1) :=\{ b\in \mI_1\, | \, ba=1$ for some
$a\in \mI_1\}$, i.e.,  $\CL (\mI_1 )$ and $\CR (\mI_1)$ are the sets of
all the left and right invertible elements of the algebra $\mI_1$
respectively. The sets $\CL (\mI_1 )$ and $\CR (\mI_1)$ are
monoids and the group $\mI_1^*$ of invertible elements of the
algebra $\mI_1$ is  also the group of invertible elements of the
monoids $\CL (\mI_1)$ and $\CR (\mI_1)$, $\CL (\mI_1)\cap \CR
(\mI_1)=\mI_1^*$. For an element $u\in \mI_1$, let $\linv (u) :=\{
v\in \mI_1\, | \, vu=1\}$ and $\rinv (u) :=\{ v\in \mI_1\, | \,
uv=1\}$, the {\em sets of left and right inverses} for the element
$u$. The next theorem describes all the left and right inverses of
elements in $\mI_1$.

\begin{corollary}\label{c25Sep10}

\begin{enumerate}
\item An element $a\in \mI_1$ admits a left inverse iff $a=a'
\int^n$ for some natural number $n\geq 0$ and an element $a'\in K^*+F$
such that  $a_{K[x]}$ is an injection (necessarily,
$n=\dim_K(\coker(a_{K[x]})))$. In this case, $\linv (a) = \{ b=
\der^n b'\, | \, b'\in K^*+F, \int^n \der^n b'a'\int^n
\der^n=\int^n \der^n \}$. $\CL (\mI_1) = \{ a\in (K^*+F)\int^n\, |
\, n\in \N, a_{K[x]}$ is an injection$\}$.
 \item An element $b\in \mI_1$ admits a right inverse iff $b\in
\der^n b'$ for some natural number $n\geq 0$ and an element $b'\in K^*+F$ 
such that  $b_{K[x]}$ is a surjection (necessarily,
$n=\dim_K(\ker(b_{K[x]})))$. In this case, $\rinv (b) = \{
a'\int^n \, | \, a'\in K^*+F, \int^n \der^n b'a'\int^n
\der^n=\int^n \der^n \}$. $\CR (\mI_1) =\{ b\in \der^n (K^*+F)\, |
\, n\in \N, b_{K[x]}$ is a surjection$\}$.
\end{enumerate}
\end{corollary}

{\it Proof}. 1. Suppose that $ba=1$ for some element $b\in \mI_1$.
Then $a_{K[x]}$ is an injection and $b_{K[x]}$ is a surjection
since $\mI_1\subseteq \End_K(K[x])$. By Theorem \ref{25Sep10}.(1),
$a=a'\int^n$ and  $b=\der^n b'$ for some elements $a',b'\in K^*+F$
such that $\der^nb'a'\int^n=1$, or, equivalently, $\int^n \der^n
b'a'\int^n \der^n=\int^n \der^n$. The second equality is obtained
from the first by applying $\int^n (\cdot ) \der^n$, and the first
equality is obtained from the second by applying $\der^n (\cdot
)\int^n$. The last equality of statement 1 follows from Theorem
\ref{25Sep10}.(1).

2. Suppose that $ba=1$ for some element $a\in \mI_1$. Then
$a_{K[x]}$ is an injection and $b_{K[x]}$ is a surjection. By
Theorem \ref{25Sep10}.(2), $a=a'\int^n$ and  $b= \der^n b'$  for
some elements $a', b'\in K^*+F$ such that $\der^n b'a'\int^n=1$,
or, equivalently, $\int^n \der^n b'a'\int^n \der^n=\int^n \der^n$.
The last equality of statement 2 follows from Theorem
\ref{25Sep10}.(2).
 $\Box $

$\noindent $

 Since $\mI_1\subseteq
\End_K(K[x])$ we identify each element of the algebra $\mI_1$ with
its matrix with respect to the basis $\{ x^i/i!\, | \, i\in \N\}$
of the vector space $K[x]$.  {\em The equality $\int^n \der^n
b'a'\int^n \der^n=\int^n \der^n$ holds iff $b'a'=\begin{pmatrix}
A & B \\ C & 1_n \\
\end{pmatrix}$ for some matrices $A\in M_n(K)$, $B$, $C$ and}
$1_n:= e_{n,n}+e_{n+1, n+1}+\cdots = 1-e_{00}-\cdots e_{n-1,
n-1}$.

{\em Proof}. $(\Rightarrow )$ Trivial.

$(\Leftarrow )$ For an arbitrary choice of the matrices $A$, $B$
and $C$,
$$ \der^n \begin{pmatrix}
A & B \\ C & 1_n \\
\end{pmatrix} \int^n = \der^n 1_n\int^n = \der^n (1-e_{00}-\cdots
- e_{n-1,n-1}) \int^n = 1.\;\;\; \Box $$
The set $\CL (\mI_1)$ is
the disjoint union 
\begin{equation}\label{LIn}
\CL (\mI_1) =\coprod_{n\in \N} \CL (\mI_1)_n
\end{equation}
where $ \CL (\mI_1)_n:= \{ a\in (K^*+F)\int^n\, | \, a_{K[x]}$ is
injective$\}=\{ a'\int^n\, | \, a'\in K^*+F, \ker_{K[x]}(a')\cap
(x^n)=0\} = \{ a\in \CL (\mI_1) \, | \, \dim_K(\coker
(a_{K[x]}))=n\} = \{ a\in \CL(\mI_1) \, | \, -\ind_{K[x]}(a) =n\}$
and $\CL (\mI_1)_0 = \mI_1^*$, by Theorem \ref{12Jun10}.(3).
Similarly, the set $\CR (\mI_1)$ is the disjoint union
\begin{equation}\label{RIn} \CR (\mI_1) =\coprod_{n\in \N} \CR
(\mI_1)_n
\end{equation}
where $ \CR (\mI_1)_n:= \{ b\in \der^n (K^*+F)\, | \, b_{K[x]}$ is
surjective$\}= \{ \der^n b'\, | \, b'\in K^*+F,
\im_{K[x]}(b')+\sum_{i=0}^{n-1}Kx^i = K[x]\}= \{ b\in \CR
(\mI_1)\, | \, \dim_K(\ker (b_{K[x]}))=n\}= \{ b\in \CR (\mI_1) \,
| \, \ind_{K[x]}(b) = n\}$ and $\CR (\mI_1)_0= \mI_1^*$, by
Theorem \ref{12Jun10}.(3). Using the additivity of the index map
$\ind_{K[x]}$ we see that $\CL(\mI_1)$ and $\CR (\mI_1)$ are
$\N$-graded monoids, that is $\CL (\mI_1)_n\CL (\mI_1)_m\subseteq
\CL (\mI_1)_{n+m}$ and $\CR (\mI_1)_n\CR (\mI_1)_m\subseteq \CR
(\mI_1)_{n+m}$ for all $n,m\in \N$. In particular, $\mI_1^* \CL
(\mI_1)_n \mI_1^* = \CL (\mI_1)_n$ and $\mI_1^* \CR (\mI_1)_n
\mI_1^* = \CR (\mI_1)_n$ for all $n\in \N$.  Since $ab=1$ iff
$b^*a^* =1$, we have the equalities 
\begin{equation}\label{LRn}
\CL (\mI_1)^* =\CR (\mI_1), \;\; \CR (\mI_1)^* =\CL (\mI_1), \;\;
\CL (\mI_1)_n^* =\CR (\mI_1)_n, \;\;  \CR (\mI_1)_n^* =\CL
(\mI_1)_n,
\end{equation}
for all $n\in \N$. In particular, for all elements $a\in \CL
(\mI_1) \bigcup \CR (\mI_1)$, 
\begin{equation}\label{indaa*}
\ind_{K[x]}(a^*) =-\ind_{K[x]}(a).
\end{equation}
In general, the equality (\ref{indaa*}) does not hold.

$\noindent $

{\em Example}.  $\ind_{K[x]}(1+\der ) =0$ but $\ind_{K[x]}((1+\der
)^* ) =\ind_{K[x]}(1+\int )=-1$.

$\noindent $

 Clearly, $$\CL (\mI_1)_n \cap \CR (\mI_1)_m
= \begin{cases}
\mI_1^*& \text{if }n=m=0,\\
\emptyset& \text{otherwise}.\\
\end{cases}$$

\begin{theorem}\label{3Oct10}
\begin{enumerate}
\item The monoid $\CL (\mI_1) $ is generated by the group
$\mI_1^*$ and the element $\int$. \item For all $n\in \N$, $\CL
(\mI_1)_n = \mI_1^* \int^n$. For all $n,m\in \N$, $\CL
(\mI_1)_n\CL (\mI_1)_m=\CL (\mI_1)_{n+m}$.   \item The monoid $\CR
(\mI_1) $ is generated by the group $\mI_1^*$ and the element
$\der$. \item For all $n\in \N$, $\CR (\mI_1)_n = \der^n \mI_1^*
$. For all $n,m\in \N$, $\CR (\mI_1)_n\CR (\mI_1)_m=\CR
(\mI_1)_{n+m}$.
\end{enumerate}
\end{theorem}

{\it Proof}. 1. Statement 1 follows from statement 2 and
(\ref{LIn}).

2.  It is obvious that $\CL (\mI_1) \supseteq \mI_1^*\int^n$. Let
$a\in \CL (\mI_1)_n$. To prove that the opposite inclusion holds
we have to show that $a=u\int^n$ for some $u\in \mI_1^*$. Clearly,
$a=a'\int^n$ where $a'\in K^*+F$ and $\ker_{K[x]}(a')\cap
(x^n)=0$. The last condition means that the columns $c_n, 
c_{n+1}, \ldots $ of the $\N \times \N$ matrix $a'\in K^*+F=
K^*+M_\infty (K)$ are linearly independent (where $c_i$, $i\in \N$
are the columns of the matrix $a'$), and so  there exists an
element $f=\sum_{j=0}^{n-1}\sum_{i\in \N} \l_{ij}e_{ij}\in F$,
$\l_{ij}\in K$, such that all the columns of the matrix $a'+f\in
K^*+F= K^*+M_\infty (K)$ are linearly independent, i.e.,  $(a'+f)\in
\mI_1^*$. The equality $f\int^n=0$ implies that $(a'+f)\int^n =
a'\int^n =a$, i.e.,  $\CL (\mI_1)_n = \mI_1^* \int^n$. Then $\CL
(\mI_1)_n \CL (\mI_1)_m= \mI_1^* \int^n \mI_1^* \int^m = \mI_1^*
\int^{n+m}= \CL (\mI_1)_{n+m}$ for all $n,m\in \N$.

3 and 4. Statements 3 and 4 are obtained from statements 1 and 2
by applying the involution $*$ of the algebra $\mI_1$, see
(\ref{LRn}). $\Box $

$\noindent $

{\em Remark}. For all $n\geq 1$, $\CL (\mI_1)_n = \mI_1^* \int^n
\varsupsetneqq \int^n \mI_1^*$ and $\CR (\mI_1)_n = \der^n \mI_1^*
\varsupsetneqq \mI_1^* \der^n$, Corollary \ref{b3Oct10}.
\begin{corollary}\label{a3Oct10}
\begin{enumerate}
\item The decomposition $\CL (\mI_1) = \bigsqcup_{n\in \N} \CL
(\mI_1)_n$ is the orbit decomposition of the action of the group
$\mI_1^*$  on $\CL (\mI_1)$ by left multiplication, $\CL (\mI_1)_n
= \mI_n^* \int^n$ and the stabilizer of $\int^n$ is equal to ${\rm
st}_{\mI_1^*} (\int^n) = \begin{pmatrix} \GL_n(K) & 0\\ * & 1
 \end{pmatrix}\subseteq (1+F)^*$. In particular, there are only
 countably many orbits and the action of $\mI_1^*$ is not free.
\item The decomposition $\CR(\mI_1) = \bigsqcup_{n\in \N} \CR
(\mI_1)_n$ is the orbit decomposition of the action of the group
$\mI_1^*$  on $\CR (\mI_1)$ by right multiplication, $\CR (\mI_1)_n
= \der^n \mI_n^* $ and the stabilizer of $\der^n$ is equal to
${\rm st}_{\mI_1^*} (\der^n) = \begin{pmatrix} \GL_n(K) & *\\ 0 &
1
 \end{pmatrix}\subseteq (1+F)^*$. In particular, there are only
 countably many orbits and the action of $\mI_1^*$ is not free.
\end{enumerate}
\end{corollary}


The next theorem describes in explicit terms all the left and right inverses for
elements of $\mI_1$.

\begin{theorem}\label{A3Oct10}
\begin{enumerate}
\item Let $a\in \CL (\mI_1)$, i.e.,  $a=a'\int^n$ where $a'\in
\mI_1^*$ and $n\in \N$ (Theorem \ref{3Oct10}.(2)). Then $\linv (a)
= (\der^n+E_{\N , 0}+\cdots + E_{\N , n-1})a'^{-1}$ where $E_{\N ,
i}:= \sum_{j\in \N} Ke_{ji}= \sum_{j\in \N} KE_{ji}$. \item Let
$\in \CR (\mI_1)$, i.e.,  $b=\der^n b'$ where $b'\in \mI_1^*$ and
$n\in \N$ (Theorem \ref{3Oct10}.(4)). Then $\rinv (a) =
b'^{-1}(\int^n+E_{0, \N }+\cdots + E_{n-1, \N })$ where $E_{i, \N
}:= \sum_{j\in \N} Ke_{ij}= \sum_{j\in \N} KE_{ij}$.
\end{enumerate}
\end{theorem}

{\it Proof}. 1. It is obvious that $\linv (a) =
\der^na'^{-1}+\ker_{\mI_1} (\cdot a)$.  Since $\ker_{\mI_1}(\cdot
\int^n) =E_{\N , 0}+\cdots + E_{\N , n-1}=:\CE$ and $a'\in
\mI_1^*$, $\ker_{\mI_1}(\cdot a'\int^n) = \CE a'^{-1}$. Therefore,
$\linv (a) = (\der^n +\CE ) a'^{-1}$.

2. Statement 2 follows from statement 1 by applying the involution
$*$ and using the equalities $\CL (\mI_1)_n^* =\CR (\mI_1)_n$ and
$\rinv (a^*) = \linv (a)^*$. $\Box $

$\noindent $

Recall that  $\mI_1^* = K^* (1+F)^* = K^*\times (1+F)^*\simeq K^*
\times \GL_\infty (K)$. We use the matrix units $\{
e_{ij}\}_{i,j\in \N}$ to define the isomorphism $(e_{ij}\lra
E_{ij})$. Define the group monomorphism $\kappa :\mI_1^* \ra
\mI_1^*$ by the rule: for all $\l \in K^*$ and $u\in (1+F)^*
\simeq \GL_\infty (K)$, 
\begin{equation}\label{kk}
\kappa (\l ) = \l, \;\; \kappa (u) = \begin{pmatrix}
 1& 0\\ 0  & u
 \end{pmatrix} .
\end{equation}
By the very definition, $\kappa (K^*) = K^*$ and $\kappa ((1+F)^*)
\subseteq (1+F)^*$. Moreover, $\kappa (u) = e_{00}+\int u \der $
and, by induction, 
\begin{equation}\label{kk1}
\kappa^n(u) = e_{00}+e_{11}+\cdots +e_{n-1, n-1} +\int^n u \der^n=
\begin{pmatrix} E_n &0 \\0  & u
 \end{pmatrix}
\end{equation}
where $E_n=e_{00}+e_{11}+\cdots +e_{n-1,n-1}$ is the $n\times n$
identity matrix. Clearly, 
\begin{equation}\label{kk11}
\kappa (u^*) = \kappa (u)^*
\end{equation}
and 
\begin{equation}\label{kk4}
\det (\kappa (u)) = \det (u)
\end{equation}
where the determinant map $\det :(1+F)^*\ra (1+F)^*$ is taken with
respect to the basis $\{ e_{ij}\}$ (or $\{ E_{ij}\}$; both
determinants coincide, see (\ref{detaKF})). For all $u\in (1+F)^*$
and all $n\geq 1$, 
\begin{equation}\label{kk2}
\int^n u = \kappa^n (u) \int^n \;\; {\rm and}\;\; u\der^n = \der^n
\kappa^n (u).
\end{equation}
By (\ref{kk1}), (\ref{kk11}) and $(\mI_1^*)^*=\mI_1^*$,  it
suffices to prove that $\int u =\kappa (u) \int$: $ \int u = \int
u \der \int = (e_{00}+\int u \der ) \int = \kappa (u) \int $. Let
$u,v\in (1+F)^*$. Then 
\begin{equation}\label{kk3}
\kappa^n (u) \int^n = \kappa (v)\int^n \; \Rightarrow \; u=v.
\end{equation}
{\it Proof}. The equality can be written as $\int^n u = \int^n v
$.  Then multiplying this equality by the element $\der^n$ on the
left yields the result. $\Box$

\begin{corollary}\label{b3Oct10}

\begin{enumerate}
\item  For all $n\geq 1$, $\CL (\mI_1)_n = \mI_1^* \int^n
\varsupsetneqq \int^n \mI_1^*$.  The right action of the group
$\mI_1^*$ on the set $\CL (\mI_1)$ is free (i.e.,  the stabilizer of
each point is the identity group). \item   For all $n\geq 1$,
$\CR (\mI_1)_n = \der^n \mI_1^* \varsupsetneqq \mI_1^* \der^n$.
The left action of the group $\mI_1^*$ on the set $\CR (\mI_1)$ is
free.
\end{enumerate}
\end{corollary}

{\it Proof}. 1. By (\ref{kk1}) and (\ref{kk2}), $\CL (\mI_1)_n =
\mI_1^* \int^n \varsupsetneqq \int^n \mI_1^*$. The freeness of the
right action of the group $\mI_1^*$  follows from (\ref{LIn}) and
the facts that $\CL (\mI_1)_n = \mI_1^* \int^n$, $\CL (\mI_1)_n
\mI_1^* \subseteq \CL (\mI_1)_n$: if
$u\int^n v = u\int^n$ where $u,v\in \mI_1^*$ then $\int^n v =
\int^n$, and so $v=1$.

2. Applying the involution $*$ to statement 1 yields statement 2.
 $\Box $

$\noindent $


\section{The algebras $\tmI1$ and $\tmJ1$}\label{TAI1J1}

The aim of  this section is to introduce  and study the algebras
$\tmI1$ and $\tmJ1$ that  have remarkable properties (Theorem
\ref{14Jun10} and Corollary \ref{a14Jun10}). Name just a few: both
algebras are obtained from $\mI_1$ by inverting certain elements;
they contain the only proper ideal, $\CC(\tmI1) \triangleleft
\tmI1$ and $\CC (\tmJ1 )\triangleleft\tmJ1$;  the factor algebras
$\tmI1/\CC (\tmI1 )$ and $\tmJ1/\CC (\tmJ1 )$ are canonically
isomorphic to the skew field of fractions $\Frac (A_1)$ of the
Weyl algebra $A_1$ and its opposite skew field $\Frac (A_1)^{op}$
respectively. The algebras $\tmI1$ and $\tmJ1$ will turn out to be
the largest right and the largest left quotients rings of the
algebra $\mI_1$ respectively (Theorem \ref{18Jun10}).

Let $V$ be an infinite dimensional vector space. The set $\CV =
\CV (V)$ of all vector subspaces $U$ of $V$ of finite codimension,
i.e.,  $\codim (U) := \dim _K(V/U)<\infty$, is a filter. This means
that

(i) if $U\subseteq W$ are subspaces of $V$ and $U\in \CV$ then
$W\in \CV$;

(ii)  if $U_1, \ldots , U_n\in \CV$ then $\bigcap_{i=1}^n U_i \in
\CV$;

(iii) $\emptyset \not\in \CV$.

Let $\CF (V)$ be the set of all Fredholm linear maps/operators in
an infinite dimensional space in $V$. The set $\CF (V)$ is a
monoid with respect to composition of maps, $\CF (V)+\CC (V) = \CF
(V)$ and $\CF (V)\bigcap \CC (V) =\emptyset$. For each element
$f\in \CF (V)$ and for each vector space $U\in \CV (V)$, $f(U),
f^{-1}(U)\in \CV (V)$. Let $A$ be a subalgebra of $\End_K(V)$ and
$\CC (A) := A \bigcap \CC (V)$ be the ideal of compact operators
of the algebra $A$. Suppose that $A\backslash \CC (A) \subseteq
\CF (V)$. Example: $V=K[x]$, $A= \mI_1\subseteq \End_K(K[x])$
 and $\CC (\mI_1)=F$ (Theorem \ref{B30May10}, Corollary \ref{aB30May10}). Note that the
disjoint union $\CF (V)\bigsqcup \CC (V)$ is a semigroup but {\em
not} a ring but $\CF (\mI_1) \coprod \CC (\mI_1) = \mI_1$ {\em is}
a ring. A subalgebra $A$ of $\End_K(V)$ satisfies the condition
$A\backslash \CC (A)\subseteq \CF (V)$ iff $A\subseteq\CF
(V)\bigsqcup \CC (V)$ iff the Compact-Fredholm Alternative holds
for the left $A$-module $V$. For such a subalgebra $A$,
 we say that two elements $a$ and $b$ of  $A$
are {\em equivalent}, $a\sim b$, if $a|_V= b|_V$ for some $V\in
\CV (V)$. The relation $\sim$ is an equivalent relation on the
algebra $A$. Let $[a]$ be the equivalence class of the element
$a$. Then the set of equivalence classes $\bA = \{ [a]\, | \, a\in
A\}$ is a $K$-algebra where $[a]+[b]:=[a+b]$, $[a][b]:=[ab]$ and
$\l [a]= [\l a ]$ for all $a,b\in A$ and $\l \in K$  (the
inclusion $A\backslash \CC (A) \subseteq \CF (A)$ guarantees that
the multiplication is well defined). Since $[a]=0$ iff $a\in \CC
(A)$, there is a canonical algebra isomorphism $A/\CC (A) \ra
\bA$, $a+\CC (A) \mapsto [a]$. The inclusion $A\backslash \CC (A)
\subseteq \CF (A)$ implies that the algebra $\bA$ is a domain.
 Similarly, the set $\bCF (V) := \CF (V) / \CC (V) = \{ [a]=
 a+\CC (V) \, | \, a\in \CF (V)\}$ is a monoid, $[a][b]= [ab]$.
 For  vector spaces $U$ and $W$, $\Iso_K(U,W)$ is the set of all
 bijective linear maps from $U$ to $W$.  Via the map $[a]\mapsto [a]$,
 the monoid $\bCF (V)$ is canonically isomorphic to
 the {\em group} $\CG (V) = \{ [a]\, | \, a\in
 \Iso_K(U,W)$ for some $U,W\in \CV (V)\}$;  $[a]=[a']$ where $a'\in
 \Hom_K(U', W')$ and $U', W'\in \CV (V)$ iff $a|_{U''}= a'|_{U''}$
 for some $U''\in \CV (V) $ such that $U''\subseteq U\cap U'$; and for $b\in
 \Iso_K(U', W')$ and $b\in \CG (V)$, $[ab]= [ ab|_{b^{-1}(W'\cap
 U)}]$ and $[a]^{-1} = [ a^{-1} :W\ra U]$. It is convenient to
 identify the groups $\CF (V)$ and $\CG (V)$ via the group
 isomorphism $[a]\mapsto [a]$.


$\noindent $

 {\it Definition}. Let the algebra $A$ be as above, i.e.,
 $A\subseteq \End_K(V)$ and $A\backslash \CC (A) \subseteq \CF (V)$.
 Then $\ffrac_V(A) :=
 \ffrac_V(\bA )$ is the intersection of all subalgebras $B$ (if
 they exist) of the set $\CG (V)\bigcup \{ 0\}$ such that  $\bA,
 \bA^{-1} \subseteq B$ where $\bA^{-1} := \{ [a]^{-1} \, | \,
 0\neq [a]\in \bA\}$. So, $\ffrac_V(\bA )$ is the subalgebra (if
  it exits) in $\CG (V) \bigcup \{ 0\}$ generated by $\bA$ and
 $\bA^{-1}$.

\begin{proposition}\label{A14Jun10}
Let the algebra $A$  be as above, i.e.,  $A\subseteq \End_K(V)$ and
$A\backslash \CC (A) \subseteq \CF (A)$. Suppose that $\bA$ is a
left (resp. right) Goldie domain and $\Frac (\bA ) $ be its left
(rep. right) skew field of fractions. Then the map $\Frac (\bA
)\ra \ffrac_V (\bA ) $, $[s]^{-1}[a]\mapsto [s]^{-1}[a]$ (resp.
$[a][s]^{-1}\mapsto [a][s]^{-1}$) is an algebra isomorphism.
\end{proposition}

{\it Proof}. The statement follows from the left (resp. right) Ore
condition and the fact that each nonzero element of the algebra
$\bA$ is invertible in $\CG (V)$. In more detail, let $\bA$ be a
left (resp. right) Goldie domain. The left (resp. right) Ore
condition in the domain $\bA$ implies that each element of $\Frac
(\bA )$ and of ${\rm frac}_V(\bA )$ has the form $[s]^{-1}[a]$
(resp. $[a][s]^{-1}$) where $[a]\in \bA$ and $0\neq [s]\in \bA$.
By the universal property of the quotient ring $\Frac (\bA )$, the
map $[s]^{-1}[a]\mapsto [s]^{-1}[a]$ (resp. $[a][s]^{-1}\mapsto
[a][s]^{-1}$) is an algebra epimorphism with zero kernel, i.e.,  an
isomorphism. $\Box $


\begin{corollary}\label{aA14Jun10}
\begin{enumerate}
\item For all $\mI_1$-modules $M$ of finite length,
$\ffrac_M(\mI_1) \simeq  \Frac (B_1) = \Frac (A_1)$. \item  For
all $A_1$-modules $M$ of finite length,  $\ffrac_M(A_1) \simeq
\Frac (A_1)$.
\end{enumerate}
\end{corollary}

{\it Proof}. 1. This follows from Theorem \ref{B30May10} and
Proposition \ref{A14Jun10}.

2. This follows from \cite{McCon-Rob-JA-73} and Proposition
\ref{A14Jun10}. $\Box $

$\noindent $

 Recall that the algebra $B_1 = K[H][\der, \der^{-1}; \tau ]$
is the (left and right) localization of the Weyl algebra $A_1$ at
the powers of the element $\der$, $B_1= S_\der^{-1}A_1$,  and
$\Frac (B_1) = \Frac (A_1)$. Moreover, the multiplicatively closed
set $B_1^0:=\{ \sum_{i\geq 0} a_{-i} \der^i \in B_1\, | \, a_0\neq
0, {\rm all}\; a_j\in K[H]\}$ is a left and right Ore set in $B_1$
such that 
\begin{equation}\label{FracB}
\Frac (A_1) = \Frac (B_1) = {B_1^0}^{-1}B_1 = B_1{B_1^0}^{-1}.
\end{equation}
Notice that under the natural algebra epimorphism $\pi : \mI_1\ra
B_1= \mI_1/F$, $a\mapsto \oa := a+F$, 
\begin{equation}\label{1FracB}
\pi (\mI_1^0) = B_1^0,
\end{equation}
by Theorem \ref{12Jun10}.(1), where $\mI_1^0:= \mI_1\bigcap
\Aut_K(K[x])$, the multiplicative submonoid of the group
$\Aut_K(K[x])$.  Let $\tmI1$ be the subalgebra of $\End_K(K[x])$
generated by $\mI_1$ and ${\mI_1^0}^{-1}$, i.e.,  it is obtained
from the algebra $\mI_1$ by adding the inverse elements of all the
elements of the set $\mI_1^0$.

\begin{theorem}\label{14Jun10}
\begin{enumerate}
\item ${\mI_1^0}^{-1}F = F$ but $F\varsubsetneqq F{\mI_1^0}^{-1}$.
\item The multiplicatively closed set $\mI_1^0$ is a right but not
left Ore subset of  $\mI_1$ and $\mI_1{\mI_1^0}^{-1} \simeq
\tmI1$. \item $F{\mI_1^0}^{-1}$ is the only proper ideal of the
algebra $\tmI1$, $(F{\mI_1^0}^{-1})^2= F{\mI_1^0}^{-1}$. \item
$\tmI1/ F{\mI_1^0}^{-1}\simeq \Frac (A_1)$. \item $\CC (\tmI1 ) =
F{\mI_1^0}^{-1}$ (where $\CC (\tmI1 ) := \tmI1 \bigcap \CC
(K[x])$). \item $\ffrac_{K[x]}(\tmI1 ) \simeq \tmI1 / \CC (\tmI1
)\simeq \Frac (A_1)$. \item
\begin{enumerate}
\item There are only two (up to isomorphism) simple left
$\tmI1$-modules, $K[x]$ and $\Frac (A_1)= \tmI1 / \CC (\tmI1 )$.
 The first one is faithful and the second one is not. \item There
are only two (up to isomorphism) simple right  $\tmI1$-modules,
$(\mI_1 / \int \mI_1) {\mI_1^0}^{-1}\simeq K[\der ]
{\mI_1^0}^{-1}$ and $\Frac (A_1)$.
 The first one is
faithful and the second one is not.
\end{enumerate}
\item ${}_{\tmI1}F=\bigoplus_{i\in \N}E_{\N , i }$ is an infinite
direct sum of nonzero left ideals $E_{\N , i } =F e_{0i}\simeq
K[x]$, therefore the algebra $\tmI1$ is not  left Noetherian.
$F{\mI_1^0}^{-1} =(\bigoplus_{i\in \N} E_{i, \N }){\mI_1^0}^{-1}
\simeq \bigoplus_{i\in \N}E_{i, \N } {\mI_1^0}^{-1}$ is an
infinite direct sum of nonzero right ideals $E_{i,\N}
{\mI_1^0}^{-1}$ where $E_{i, \N } = e_{i0}F$, therefore the
algebra $\tmI1$ is not right Noetherian.
\end{enumerate}
\end{theorem}

{\it Proof}. 1. The equality ${\mI_1^0}^{-1}F=F$ follows from the
inversion formula (Theorem \ref{12Jun10}.(4)): $\begin{pmatrix} a
& b\\ 0 &c
 \end{pmatrix} \begin{pmatrix}
a' & 0 \\ 0 &0
 \end{pmatrix}=\begin{pmatrix}
aa' & 0\\  0&0
 \end{pmatrix}\in F$.
 The inclusion $F\subset
F{\mI_1^0}^{-1}$ is strict as $e_{00} (1-\der)^{-1} =
e_{00}\sum_{i\geq 0} \der^i = \sum_{i\geq 0} e_{0i}\in
F{\mI_1^0}^{-1}\backslash F$.

2. Notice that all elements of the set $\mI_1^0$ are (left and
right) regular, i.e.,  non-zero divisors in $\mI_1$. To prove that
$\mI_1^0$ is a right Ore set we have to show that for all elements
$a\in \mI_1$ and $s\in \mI_1^0$ there exist elements $b\in \mI_1$
and $ t\in \mI_1^0$ such that $at = sb$. By (\ref{FracB}) and
(\ref{1FracB}), $at=sb_1+f$ for some elements $t\in \mI_1^0$,
$b_1\in \mI_1$  and $f\in F$. By statement 1, $s_1^{-1} f\in F$,
and so $at=s(b_1+s^{-1}f)=sb$ where $b=b_1+s^{-1}f\in \mI_1$, as
required. Recall that the algebra $\tmI1$ is the subalgebra of the
endomorphism algebra $\End_K(K[x])$ generated by the algebra
$\mI_1$ and ${\mI_1^0}^{-1}$. Since $\mI_1^0$ is a right Ore set,
there is a natural algebra epimorphism $\mI_1 {\mI_1^0}^{-1}\ra
\tmI1$ such that its restriction to the subalgebra $\mI_1$ gives
the identity map $\mI_1\ra \mI_1\subseteq \tmI1$. Since the
subalgebra $\mI_1$ is an essential right $\mI_1$-module of
$\mI_1{\mI_1^0}^{-1}$, the epimorphism is necessarily an
isomorphism.

The set $\mI_0$ is not a left Ore set, otherwise we would have the
equality ${\mI_1^0}^{-1}F = F{\mI_1^0}^{-1}$ which would have
contradicted to statement 1 (the equality follows easily from the
fact that $F$ is an ideal of $\mI_1$ and the algebra $\mI_1/ F$ is
a domain. In more detail, we have to show that every fraction
$s^{-1} f$ (resp. $fs^{-1}$) where $s\in {\mI_1^0}$ and $f\in F$
can be written as $gt^{-1}$ (resp. $t^{-1} g$) for some $t\in
{\mI_1^0}$ and $g\in F$. Notice that $s^{-1} f = gt^{-1}$ (resp.
$fs^{-1} = t^{-1} g$) for some $ t\in {\mI_1^0}$ and $g\in \mI_1$,
i.e.,  $ft=sg\in F$ (resp. $tf = gs\in F$). Since $\mI_1/F$ is a
domain and $s\in {\mI_1^0}$ we conclude that $g\in F$ (resp. $g\in
F$)).

3. By statement 1,
$F{\mI_1^0}^{-1}=\mI_1{\mI_1^0}^{-1}F\mI_1{\mI_1^0}^{-1}$ is a
proper ideal of the algebra $\mI_1{\mI_1^0}^{-1}$  (since $0\neq
F\subseteq F{\mI_1^0}^{-1}$ and $F{\mI_1^0}^{-1}\neq
\mI_1{\mI_1^0}^{-1}$ as $F\bigcap \mI_1^0 = \emptyset$).
Therefore, by statement 1,
$(F{\mI_1^0}^{-1})^2=F{\mI_1^0}^{-1}F{\mI_1^0}^{-1}=F^2{\mI_1^0}^{-1}=F{\mI_1^0}^{-1}$
since $F=F^2$. Since $F$ is the only proper ideal of the algebra
$\mI_1$ and an essential $\mI_1$-bimodule \cite{algintdif}, the
ideal $F{\mI_1^0}^{-1}$ is the only proper ideal of the algebra
$\tmI1$.

4. $\tmI1 / F{\mI_1^0}^{-1}=
\mI_1{\mI_1^0}^{-1}/F{\mI_1^0}^{-1}\simeq
(\mI_1/F){\mI_1^0}^{-1}\stackrel{(\ref{1FracB})}{=}B_1{B_1^0}^{-1}=\Frac
(B_1)=\Frac (A_1)$.

5. Since $F{\mI_1^0}^{-1}\subseteq \CC (K[x])$, we have the
inclusion $F{\mI_1^0}^{-1}\subseteq \tmI1 \bigcap \CC (K[x])= \CC
 (\tmI1 )$ of proper ideals of the algebra $\tmI1$. By statement 3,
the inclusion is the equality.

6. By statements 4 and 5, the factor algebra $\tmI1 / \CC (\tmI1
)= \tmI1 / F{\mI_1^0}^{-1}\simeq \Frac (A_1)$ is a skew field.
Therefore, $\ffrac_{K[x]}(\tmI1 ) = \tmI1 / \CC (\tmI1 ) \simeq
\Frac (A_1)$.

7(a). Let $M$ be a simple left $\tmI1$-module. Then either
$F{\mI_1^0}^{-1}M=0$ or $F{\mI_1^0}^{-1}M=M$. In the first case,
$M\simeq \tmI1 / F{\mI_1^0}^{-1}\simeq \Frac (A_1)$. In the second
case, $FM\neq 0$. Since $F=\sum_{i\geq 0}Fe_{0i}$, must have
$e_{0i}m\neq 0$ for some nonzero element $m$ of $M$ and some
natural number $i$. Then ${}_{\tmI1}Fe_{0i}m\simeq K[x]$ since
${\mI_1^0}^{-1}F=F$, by statement 1.

7(b). Let $N$ be a simple right $\tmI1$-module. Then either
$NF{\mI_1^0}^{-1}=0$ or $NF{\mI_1^0}^{-1}=N$. In the first case,
$N\simeq \tmI1 / F{\mI_1^0}^{-1}\simeq \Frac (A_1)$. In the second
case, $NF\neq 0$. Since $F=\sum_{i\geq 0}e_{i0}F$, must have
$ne_{i0}\neq 0$ for some nonzero element $n$ of $N$ and some
natural number $i$. Then the right $\mI_1$-module $ne_{i0}F\simeq
\mI_1/ \int \mI_1\simeq K[\der ]$ is simple, hence so is its
localization $ne_{i0}F{\mI_1^0}^{-1}$ provided it is not equal to
zero which is the case as $0\neq n\in N_{\tmI1 }$.

8. Statement 8 is obvious. $\Box $

$\noindent $

Theorem \ref{14Jun10} shows  that the algebra $\tmI1$ is not
left-right symmetric. In particular, the involution $*$ of the
algebra $\mI_1$ {\em cannot} be extended to an involution of the
algebra $\tmI1$ since otherwise we would have the inclusion
$(\mI_1^0)^* \subseteq \Aut_K(K[x])$ which is not the case as
$1+\der\in \mI_1^0$ but $(1+\der )^* = 1+\int\not\in
\Aut_K(K[x])$, by Theorem \ref{12Jun10}.(1).

In fact, we {\em can} `extend' the involution $*$ of the algebra
$\mI_1$  in the following sense: we will construct an algebra
$\tmJ1$, that contains the algebra $\mI_1$, and two $K$-algebra
{\em anti-isomorphisms} 
\begin{equation}\label{IJab}
*:\tmI1 \ra \tmJ1, \;\; a\mapsto a^*, \;\; *:\tmJ1\ra\tmI1, \;\;
b\mapsto b^*,
\end{equation}
such that $a^{**}=a$, $b^{**} = b$, $\mI_1^* =\mI_1$ and the
restriction of the map $*$ in (\ref{IJab}) to the subalgebra
$\mI_1$ coincides with the involution $*$. Recall that a
$K$-algebra {\em anti-isomorphism} $\v : A\ra B$ is a bijective
linear map such that $\v (uv) = \v (v) \v (u)$ for all elements
$u,v\in A$; equivalently,  $\v :A\ra B^{op}$ is a $K$-algebra
isomorphism where $B^{op}$ is the opposite algebra to $B$. The
equalities $(\mI_1\der )^* = \int \mI_1$ and $(\int \mI_1)^* =
\mI_1\der$ yield the $K$-linear isomorphisms: 
\begin{equation}\label{KId1}
K[\int]\simeq \mI_1/\mI_1\der\stackrel{*}{\ra}\mI_1/ \int \mI_1
\simeq K[\der ] , \;\; p+\mI_1\der \mapsto p^*+\int \mI_1,
\end{equation}
\begin{equation}\label{KId2}
K[\der]\simeq \mI_1/\int \mI_1\stackrel{*}{\ra}\mI_1/  \mI_1\der
\simeq K[\int ] , \;\; q+\int \mI_1 \mapsto q^*+ \mI_1\der ,
\end{equation}
such that $(am)^* = m^*a^*$, $m^{**} = m$, $(na)^*=a^* n^*$ and
$n^{**} = n$ for all elements $a\in \mI_1$, $m\in K[\int]$ and
$n\in K[\der ]$. Notice that $K[\int ]$ is a faithful simple left
$\mI_1$-module isomorphic to the faithful simple $\mI_1$-module
$K[x]$, and $\mI_1/\int \mI_1$ is a faithful simple right
$\mI_1$-module. These are the only (up to isomorphism) faithful
simple left and right modules of the algebra $\mI_1$ respectively.

For the right $\mI_1$-module $K[\der ]$, let $\End_K(K[\der
])^{op}$ be its $K$-endomorphism algebra  where we write the
argument of a linear maps $\v$ of $\End_K(K[\der ])^{op}$ on the
{\em left}, i.e.,  $(\cdot ) \v : K[\der ] \ra K]\der ]$, $p\mapsto
(p)\v$. The upper script `op' indicates this fact, it also means
that the algebra $\End_K(K[\der ])^{op}$ is the {\em opposite}
algebra to the usual $K$-endomorphism algebra $\End_K(K[\der ])$
where we write the argument of a linear map on the {\em right}.
Since the right $\mI_1$-module $K[\der ]$ is faithful, there is
the algebra monomorphism $\mI_1\ra \End_K(K[\der ])^{op}$, $
\mapsto (\cdot a: p\mapsto pa)$. We identify the algebra $\mI_1$
with its isomorphic image, i.e.,  $\mI_1\subseteq \End_K(K[\der
])^{op}$. Let $\tmJ1$ be the subalgebra of $\End_K(K[\der] )^{op}$
generated by the algebra $\mI_1$ and the set $\mJ_1^0:=
\mI_1\bigcap \Aut_K(K[\der ])^{op}$ where $\Aut_K(K[\der ])^{op}$
is the group of units of the algebra $\End_K(K[\der ] )^{op}$. Using
(\ref{KId1}) and (\ref{KId2}), we see that 
\begin{equation}\label{KID3}
(\mI_1^0)^* = \mJ_1^0, \;\; (\mJ_1^0)^* = \mI_1^0.
\end{equation}
Therefore, the maps (\ref{KId1}) and (\ref{KId2}) yield the
$K$-algebra anti-isomorphisms (\ref{IJab}).

\begin{corollary}\label{a2Jul10}
Let $a\in \tmI1$ and $b\in \tmJ1$. Then
\begin{enumerate}
\item $\ker(a_{K[\int]}\cdot )^* = \ker(\cdot a_{K[\der]}^*)$ and
$\coker(a_{K[\int]}\cdot )^* = \coker(\cdot a_{K[\der]}^*)$. \item
$\ker(\cdot b_{K[\der]})^* = \ker(b_{K[\int]}^*\cdot )$ and
$\coker(\cdot b_{K[\der]})^* = \coker(b_{K[\int]}^*\cdot)$.
\end{enumerate}
\end{corollary}

{\it Proof}. Straightforward, see (\ref{IJab}). $\Box $

$\noindent $

The next corollary is a straightforward consequence of
(\ref{IJab}) and Theorem \ref{14Jun10}.

\begin{corollary}\label{a14Jun10}
Let $\mJ_1^0:= \mI_1\bigcap \Aut_K(K[\der])^{op}$ and $\tmJ1$ be
the subalgebra of $\End_K(K[\der ] )^{op}$ generated by the
algebra $\mI_1$  and ${\mJ_1^0}^{-1}$. Then
\begin{enumerate}
\item $F{\mJ_1^0}^{-1} = F$ but $F\varsubsetneqq {\mJ_1^0}^{-1}F$.
\item The multiplicatively closed set $\mJ_1^0$ is a left  but not
 right  Ore subset of  $\mI_1$ and ${\mJ_1^0}^{-1}\mI_1 \simeq \tmJ1$.
\item ${\mJ_1^0}^{-1}F$ is the only proper ideal of the algebra
$\tmJ1$, $({\mJ_1^0}^{-1}F)^2= {\mJ_1^0}^{-1}F$. \item $\tmJ1/
{\mJ_1^0}^{-1}F\simeq \Frac (A_1)^{op}$. \item $\CC (\tmJ1 ) =
{\mJ_1^0}^{-1}F$. \item $\ffrac_{K[\der]}(\tmJ1 ) \simeq \tmJ1 /
\CC (\tmJ1 )\simeq \Frac (A_1)^{op}$. \item
\begin{enumerate}
\item There are only two (up to isomorphism) simple right
$\tmJ1$-modules, $K[\der]$ and $\Frac (A_1)^{op}= \tmJ1 /
{\mJ_1^0}^{-1}F$. The first one is faithful and the second one is
not. \item There are only two (up to isomorphism) simple left
$\tmJ1$-modules ${\mJ_1^0}^{-1}(\mI_1 /  \mI_1\der ) \simeq
{\mJ_1^0}^{-1}K[x]$ and $\Frac (A_1)^{op}$.
 The first one is
faithful and the second one is not.
\end{enumerate}
\item $F_{\tmJ1}=\bigoplus_{i\in \N}E_{i,\N }$ is an infinite
direct sum of nonzero right ideals $E_{i,\N } = e_{i0}F\simeq
K[\der ]_{\tmJ1}$, therefore the algebra $\tmJ1$ is not  right
Noetherian. ${\mJ_1^0}^{-1} F={\mJ_1^0}^{-1}(\bigoplus_{i\in \N}
E_{\N , i}) \simeq \bigoplus_{i\in \N} {\mJ_1^0}^{-1}E_{\N , i}$
is an infinite direct sum of nonzero left ideals
${\mJ_1^0}^{-1}E_{\N, i}$ where $E_{\N , i} = Fe_{0i}\simeq K[x]$,
therefore the algebra $\tmJ1$ is not left Noetherian.
\end{enumerate}
\end{corollary}


\section{The largest left and right quotient rings of the algebra $\mI_1$}\label{LLRQR}

In this section, it is proved that neither left nor right quotient
ring for the algebra $\mI_1$ exists, and the largest left and the
largest right quotient rings of $\mI_1$ are found, they are not
$\mI_1$-isomorphic but $\mI_1$-anti-isomorphic (Theorem
\ref{18Jun10}). The sets of right regular, left regular and
regular elements of the algebra $\mI_1$ are described (Lemma
\ref{b18Jun10}.(1), Corollary \ref{c18Jun10}.(1) and Corollary
\ref{b2Jul10}).

Let $R$ be a ring. An element $r$ of a ring $R$ is {\em right
regular} if $rs=0$ implies $s=0$ for $s\in R$. Similarly, {\em
left regular} is defined, and {\em regular} means both right and
left regular (and hence not a zero divisor). We denote by $\CC_R$,
$\CC_R'$ and ${}'\CC_R$ the sets of regular, right regular and
left regular elements of $R$ respectively. All these sets are
monoids. A multiplicatively closed subset $S$ of $R$ is said to be
a {\em right Ore set} if it satisfies the {\em right Ore
condition} if, for each $r\in R$ and $s\in S$, $rS\bigcap sR\neq
\emptyset$. Let $S$ be a (non-empty) multiplicatively closed
subset of $R$, and let $\ass (S) :=\{ r\in R\, | \, rs=0$ for some
$s\in S\}$. Then a {\em right quotient ring} of $R$ with respect
to $S$ (a {\em right localization} of $R$ at $S$) is a ring $Q$
together with a homomorphism $\v :R\ra Q$ such that

(i) for all $s\in S$, $\v (s)$ is a unit of $Q$,

(ii) for all $q\in Q$, $q=\v (r) \v (s)^{-1}$ for some $r\in R$,
$s\in S$, and

(iii) $\ker (\v ) = \ass (S)$.

If such a ring $Q$ exists the ring $Q$ is unique up to isomorphism, usually it is
denoted by $RS^{-1}$. For a right Ore set $S$, the set $\ass (S)$
is an ideal of the ring $R$.  Recall that $RS^{-1}$ exists iff $S$
is a right Ore set and the set $\bS=\{ s+\ass (S)\in R/\ass (S)\,
| \, s\in S\} $ consists of regular elements (\cite{MR}, 2.1.12).
Similarly, a {\em left Ore set}, the {\em left Ore condition} and
the {\em left quotient ring} $S^{-1}R$ are defined. If both rings 
$S^{-1}R$ and $RS^{-1}$ exist then they are isomorphic (\cite{MR},
2.1.4.(ii)). The right quotient ring of $R$ with respect to the
set $\CC_R$ of all regular elements is called the {\em right
quotient ring} of $R$. If it exists, it is denoted by $\Frac_r(R)$ or
$Q_r(R)$. Similarly, the {\em left quotient ring},
$\Frac_l(R)=Q_l(R)$, is defined. If both left and right quotient
rings of $R$ exist then they are isomorphic and we write simply
$\Frac (R)$ or $Q(R)$ in this case. We will see that neither
$\Frac_l(\mI_1)$ nor $\Frac_r(\mI_1)$ exists (Theorem
\ref{18Jun10}.(1)). Therefore, we introduce the following new
concepts.

$\noindent $

{\it Definition}. For a ring $R$, a maximal with respect to
inclusion right Ore set $S$ of regular elements of $R$ (i.e.,
$S\subseteq \CC_R$) is called a {\em maximal regular right Ore
set} in $R$, and the quotient ring $RS^{-1}$ is called a {\em
maximal right quotient ring} of $R$. If a maximal right Ore set
$S$ is unique we say that $S$ is the {\em largest regular right
Ore set} in $R$, and the quotient ring $RS^{-1}$ is called the
{\em largest right quotient ring} of $R$ denoted $\Frac_r (R)$. In
\cite{larglquot}, it is proved that $\Frac_r (R)$ exists for an
arbitrary ring.

$\noindent $

Notice that if $\CC_R$ is a right Ore set then $\CC_R$ is the
largest regular right Ore set and the right quotient ring
$\Frac_r(R)$ is the largest right quotient ring of $R$. That is
why we use the same notation $\Frac_r(R)$ for the largest right
quotient ring of $R$. Similarly, a {\em maximal regular left Ore
set}, the {\em largest regular left Ore set} and the {\em largest
left quotient ring} of $R$, $\Frac_l(R)$, are defined.  The two
natural questions below have negative solutions as the case of the
algebra $R=\mI_1$ demonstrates (Theorem \ref{18Jun10}).

{\it Question 1}. Is the largest regular right Ore set of $R$ also
the largest regular left Ore of $R$?

{\it Question 2}. Are the rings $\Frac_l(R)$ and $\Frac_r(R)$
$R$-isomorphic, i.e.,  there is a ring isomorphism $\v :
\Frac_l(R)\ra \Frac_r(R)$ such that $ \v (r) = r$ for all elements
$r\in R$?

$\noindent $ It is obvious that if $S$ is the largest regular left
Ore set of $R$ which is also a right Ore set and if $S'$ is the
largest regular right Ore set of $R$ which is a left Ore set then
$S=S'$ and $\Frac_l(R) \simeq \Frac_r(R)$. In this case,  we
simply write $\Frac (R)$.

The next proposition gives answers to two questions:

{\it Question 3}. {\em What is the group $\tmI1^*$ of units of the
algebra $\tmI1$?}

{\it Question 4}.  {\em What is the image of the group $\tmI1^*$
under the natural algebra epimorphism $\pi : \tmI1\ra \tmI1/ \CC
(\tmI1 )=\Frac (B_1)$?}

Notice that the obvious concept of the degree $ \deg_\der$ on
$B_1$ can be extended to $\Frac (B_1)$ by the rule $\deg_\der
(s^{-1} a) = \deg_\der (a)-\deg_\der (s)$.
\begin{proposition}\label{A21Jun10}
\begin{enumerate}
\item $\tmI1^* = \mI_1^0{\mI_1^0}^{-1}:=\{ ts^{-1}\, | \, t,s\in
\mI_1^0\}$. \item  $\pi ( \tmI1^*) = B_1^0{B_1^0}^{-1}\neq \Frac
(B_1)^*$. \item $B_1^0{B_1^0}^{-1}$ is a normal subgroup of the
group $\Frac (B_1)^*$ of units of the skew field $\Frac (B_1)$.
\item $\Frac (B_1)^*=\bigsqcup_{i\in \Z}B_1^0{B_1^0}^{-1}\der^i$,
a disjoint union.\item  The group $B_1^0{B_1^0}^{-1}$ is the
kernel of the group epimorphism $\deg_\der :\Frac (B_1)^*=\Frac
(A_1)^*\ra \Z$.
\end{enumerate}
\end{proposition}

{\it Proof}. 1. Let $u\in \tmI1^*$. Then $u=ts^{-1}$ for some
elements $t\in \mI_1$ and $s\in \mI_1^0$ (Theorem
\ref{14Jun10}.(2)). Clearly, $t\in \mI_1\bigcap \tmI1^*\subseteq
\mI_1\bigcap \Aut_K(K[x])=\mI_1^0$.

2. The equality follows from statement 1 and (\ref{1FracB}). The
inequality follows from statement 4.

4. Since $B_1=K[H][\der, \der^{-1}; \tau ]$, statement 4 is
obvious.

3. Since $B_1^0{B_1^0}^{-1}\der^i= \der^iB_1^0{B_1^0}^{-1}$,
statement 3 follows from statement 4.

5. Since $\deg_\der (B_1^0{B_1^0}^{-1})=0$ and $\deg_\der (\der^i)
= i$ for all $i\in \Z$, statement 5 follows from statement 4.
$\Box $

$\noindent $

Let $S$ be a right Ore set of a ring $R$, $\bR :=R/\ass (A)$ and
$\bS = \{ s+\ass (S)\, | \, s\in S\}$. The right Ore set $S$ such
that the elements of $\bS$ are regular in $\bR$ (and so $RS^{-1}$
exists) is called a {\em right denominator set}. Similarly, a {\em
left denominator set} is defined.

\begin{corollary}\label{aA21Jun10}
\begin{enumerate}
\item The multiplicatively closed set $S_\der:=\{ \der^i\, | \,
i\in \N\}$ is a left denominator set in $\tmI1$ with $\ass (S_\der
) = F{\mI_1^0}^{-1}$ and $S_\der^{-1} \tmI1 \simeq \Frac (B_1)$.
\item The multiplicatively closed set $S_{\int} := \{ \int^i\, |
\, i\in \N\}$ is a right denominator set in $\tmJ1$ with $\ass
(S_{\int} ) = {\mJ_1^0}^{-1}F$ and $ \tmJ1 S_{\int}^{-1}\simeq
\Frac (B_1)^{op}$.
\end{enumerate}
\end{corollary}

{\it Proof}. 1. Since $F=\bigcup_{i\geq 1} \ker_F(\der^i\cdot )$
and $F{\mI_1^0}^{-1}$ is the only proper ideal of the algebra
$\tmI1$ (Theorem \ref{14Jun10}.(3)), $\ass (S_{\der }) =
F{\mI_1^0}^{-1}$ where $\ass (S_\der ) = \{ a\in \mI_1 \, | \,
\der^ia=0$ for some $i\geq 1\}$. Since $\tmI1 / \ass
(S_{\der})=\tmI1 / F{\mI_1^0}^{-1}\simeq \Frac (B_1)$ (Theorem
\ref{14Jun10}.(5,6)), we have the isomorphism $S_\der^{-1}
\tmI1\simeq \Frac (B_1)$, by Proposition \ref{A21Jun10}.(4).

2. Statement 2 is obtained from statement 1 by using the
anti-isomorphism $*$, see  (\ref{IJab}). $\Box $

$\noindent $

 Lemma
\ref{b18Jun10}.(1), Corollary \ref{c18Jun10}.(1) and Corollary
\ref{b2Jul10} describe the sets of  right regular, left regular
and regular elements of the algebra $\mI_1$ respectively.

\begin{lemma}\label{b18Jun10}
Let $a\in \mI_1$.
\begin{enumerate}
\item The following statements are equivalent.
\begin{enumerate}
\item The map $a_{\mI_1}\cdot$ is an injection. \item  The map
$a_F\cdot$ is an injection. \item  The map $a_{K[x]}\cdot$ is an
injection (see Theorem \ref{16Jun10}).
\end{enumerate}
\item The following statements are equivalent.
\begin{enumerate}
\item The map $a_{\mI_1}\cdot$ is a surjection. \item $a=\l
\der^i+f$ for some $\l \in K^*$, $i\geq 0$ and $f\in F$ such that
the map $a_{K[x]}\cdot$ is surjective (see Theorem \ref{15Jun10}).
\end{enumerate}
\end{enumerate}
\end{lemma}

{\it Proof}. 1. If $a\in F$ then none for the three statements is
true. If $a\not\in F$ then the three statements hold since
$\ker_{\mI_1}(a\cdot ) =\ker_F(a\cdot )$ (as $B_1=\mI_1/ F$ is a
domain) and ${}_{\mI_1}F\simeq K[x]^{(\N )}$.

2. Applying the Snake Lemma to the commutative diagram of right
$\mI_1$-modules $$
\xymatrix{0\ar[r] & F\ar[r]\ar[d]^{a\cdot}  & \mI_1 \ar[r]\ar[d]^{ a\cdot} & B_1 \ar[r]\ar[d]^{a\cdot} & 0 \\
0\ar[r] & F\ar[r]  & \mI_1 \ar[r] & B_1\ar[r] & 0 }
$$
yields the long exact sequence 
\begin{equation}\label{res1}
0\ra \ker_F(a\cdot)\ra \ker_{\mI_1} (a\cdot)\ra \ker_{B_1}
(a\cdot)\ra  \coker_F(a\cdot)\ra \coker_{\mI_1} (a\cdot)\ra
\coker_{B_1} (a\cdot)\ra 0.
\end{equation}
$(a)\Rightarrow (b)$ If the map $a_{\mI_1}\cdot$ is a surjection
then so is the map $a_{B_1}\cdot$, i.e.,  $\overline{a}:=a+F\in
B_1^*$, and so $\overline{a}=\l \der^i$ for some $\l\in K^*$ and
$i\in \Z$, and then, by (\ref{res1}), $\coker_F(a\cdot
)=\coker_{\mI_1}(a\cdot )=0$. Since ${}_{\mI_1}F\simeq K[x]^{(\N
)}$, the map $a_{K[x]}$ is a surjection. By Theorem \ref{15Jun10},
$i\geq 0$.

$(a)\Leftarrow (b)$ If statement (b) holds then $a+F\in B_1^*$ and
$\coker_F(a\cdot ) =0$ (since ${}_{\mI_1}F\simeq K[x]^{(\N )}$),
then, by (\ref{res1}), $0=\coker_F(a\cdot ) =
\coker_{\mI_1}(a\cdot )$. $\Box $

$\noindent $

It is obvious that $a_{\mI_1}\cdot$ is an injection (resp. a
surjection) iff $\cdot a_{\mI_1}^*$ is an injection (resp. a
surjection). Using the involution $*$ defined in (\ref{IJab}),
$F^* = F$, ${}_{\mI_1}K[x]\simeq \mI_1/ \mI_1 \der \simeq K[\int
]$ and $(\mI_1/ \mI_1 \der )^* =\mI_1/\int \mI_1 \simeq K[\der
]_{\mI_1}$, we obtain  similar results but for right
multiplication.

\begin{lemma}\label{c18Jun10}
Let $a\in \mI_1$.
\begin{enumerate}
\item The following statements are equivalent.
\begin{enumerate}
\item The map $\cdot a_{\mI_1}$ is an injection. \item  The map
$\cdot a_F$ is an injection. \item  The map $\cdot a_{K[\der]}$ is
an injection.
\end{enumerate}
\item The following statements are equivalent.
\begin{enumerate}
\item The map $\cdot a_{\mI_1}$ is a surjection. \item $a=\l
\int^i+f$ for some $\l \in K^*$, $i\geq 0$ and $f\in F$ such that
the map $\cdot a_{K[\der ]}$ is surjective (see Theorem
\ref{16Jun10}).
\end{enumerate}
\end{enumerate}
\end{lemma}

\begin{corollary}\label{b2Jul10}
Let $a\in \mI_1$. Then $a\in \CC_{\mI_1}$ iff the maps $a_{K[x]}$
and $a_{K[x]}^*$ are injections (Theorem \ref{16Jun10}).
\end{corollary}

The next lemma is a useful criterion of left regularity.

\begin{lemma}\label{d18Jun10}
Let $a\in \mI_1\backslash F$. Then the map $\cdot a_{\mI_1}$ is an
injection iff for all natural numbers $n\geq 0$,
$(e_{00}+e_{11}+\cdots +e_{nn})\im (a_{K[x]})=K[x]_{\leq n}$. For
example, the $\cdot (\der +\int )_{\mI_1}$ is an injection.
\end{lemma}

{\it Proof}. The map $\cdot a_{\mI_1}$ is not an injection iff
there exists a nonzero element $f\in F$ such that $fa=0$ since
$a\not\in F$ and $\mI_1/F$ is a domain iff $(e_{00}+e_{11}+\cdots
+e_{nn})\im (a_{K[x]})\neq K[x]_{\leq n}$ for some natural number
$n$ (e.g. $n=\deg_F(f)$ since $f=f(e_{00}+e_{11}+\cdots
+e_{nn})$):  if $V:=(e_{00}+e_{11}+\cdots +e_{nn})\im
(a_{K[x]})\neq K[x]_{\leq n}$ then $K[x]_{\leq n}= V\bigoplus V'$
for some nonzero subspace $V'$ of $K[x]_{\leq n}$ and take $f$ to
be the projection onto $V'$).

Let $a=\der +\int$ then since $a* x^{[0]}=x^{[1]}$, $a*x^{[i]}=
x^{[i-1]}+x^{[i+1]}$ for $i\geq 1$, we have $(e_{00}+e_{11}+\cdots
+e_{nn})\im (a_{K[x]})=K[x]_{\leq n}$ for all $n\geq 0$. Then, the
map $\cdot a_{\mI_1}$ is an injection, by Lemma \ref{d18Jun10}.
 $\Box $

$\noindent $

Let $R$, $S$ and $T$ be rings such that $R\subseteq S$ and
$R\subseteq T$. If there exists a ring isomorphism (respectively,
an anti-isomorphism) $\v :S\ra T$ such that $\v (r) =r$ for all
$r\in R$  (respectively, $\v (R) = R$) we say that the rings $S$
and $T$ are $R$-{\em isomorphic} (respectively, $R$-{\em
anti-isomorphic}).

\begin{theorem}\label{18Jun10}
\begin{enumerate}
\item The set $\CC_{\mI_1}$ of regular elements of the algebra
$\mI_1$ satisfies neither left nor right Ore condition. Therefore,
the left and the right  quotient rings of $\mI_1$,
$\CC_{\mI_1}^{-1}\mI_1$ and $\mI_1\CC_{\mI_1}^{-1}$,  do not
exist. \item The set $\mI_1^0$ is the largest regular right Ore
set in $\mI_1$, and so the largest right quotient ring of
fractions $\Frac_r(\mI_1):=\mI_1{\mI_1^0}^{-1}=\tmI1$ of $\mI_1$
exists.
\item The set $\mJ_1^0=(\mI_1^0)^*$ is the largest regular left
Ore set in $\mI_1$, and so the largest left  quotient ring of
fractions $\Frac_l(\mI_1):={\mJ_1^0}^{-1}\mI_1 =\tmJ1$  of $\mI_1$
exists.
\item The rings $\Frac_r(\mI_1)$ and $\Frac_l(\mI_1)$ are not
$\mI_1$-isomorphic but are
 $\mI_1$-anti-isomorphic (see (\ref{IJab})). In particular, the
 largest  regular right Ore set in $\mI_1$ is not a left
  Ore set and the largest  regular left  Ore set in $\mI_1$ is not
 a right Ore set.
\end{enumerate}
\end{theorem}

{\it Proof}. 1. By Lemma \ref{d18Jun10}, the element $a:= \der
+\int \in \mI_1$ is left regular, hence it is regular since
$a^*=a$.  To prove statement 1 it suffices to  show that
$\mI_1e_{00}\bigcap \mI_1 a=0$ and $e_{00}\mI_1\bigcap a\mI_1=0$.
In fact, it suffices to show that the first equality holds since
then the second equality can be obtained from the first by
applying the involution $*$. Suppose that the first equality
fails, then we can find a nonzero element, say $u$, in the
intersection, we seek a contradiction.  Since
$\mI_1e_{00}=\bigoplus_{i\in \N} Ke_{i0}$, the element $u$ is the
unique sum $\sum_{i\in I}\l_ie_{i0}$ where all $\l_i\in K^*$ and
$I$ is a non-empty finite subset of natural numbers. Clearly,
$u=ga$ for some element $g\in \mI_1$. Since $u\in F$ and $a\not\in
F$, we must have $g\in F$. Since $F_{\mI_1} =\bigoplus_{i\in \N}
E_{i,\N }$ is the direct sum of right submodules $E_{i,\N }=
\bigoplus_{j\in \N} Ke_{ij}\simeq K[\der ]_{\mI_1}$ and $e_{i0}\in
E_{i,\N }$ for all $i\in I$, we see that
$e_{i0}=\l_i^{-1}e_{ii}ga$, i.e.,  the element $1\in K[\der ]$
belongs to $\im (\cdot a_{K[\der ]})$ which is obviously
impossible since $\der^i*a=
\begin{cases}
\der& \text{if }i=0,\\
\der^{i+1}+\der^{i-1}& \text{if }i>0.\\
\end{cases}
$

2. Suppose that $S$ is a regular right Ore set in $\mI_1$. We have
to show that $S\subseteq \mI_1^0$. Let $u\in S$, we have to show
that $u\in \mI_1^0$. Recall that $\tmI1= \mI_1{\mI_1^0}^{-1}$,
$F{\mI_1^0}^{-1}$ is the only proper ideal of the algebra
 $\tmI1$ and $\tmI1 / F{\mI_1^0}^{-1}\simeq \Frac (B_1)$
 (Theorem \ref{14Jun10}). Let $\pi
 : \tmI1 \ra \tmI1 / F{\mI_1^0}^{-1}$, $a\mapsto \overline{a}:=
 a+F{\mI_1^0}^{-1}$. None of the elements of the ideal
 $F{\mI_1^0}^{-1}$ is regular. Therefore, $u\not\in
 F{\mI_1^0}^{-1}$. Then $\pi (u) \in \Frac (B_1)^*$, and so
 $u=\der^ist^{-1}+ft_1^{-1}$, $\int^jst^{-1}+ft_1^{-1}$ for some
 $i\geq 0$, $j>0$,  $s,t,t_1\in \mI_1^0$ and $f\in F$ (Proposition \ref{A21Jun10}).

 {\em Case 1}: $u=\der^ist^{-1}+ft_1^{-1}$. Taking a common right
 denominator we may assume that $t=t_1$ (changing $s$ and $f$ if
 necessary), i.e.,  $u=(\der^is+f)t^{-1}$. Since $\der^i F=F$, we
 may assume that $u=\der^i (s+g)t^{-1}$ for some element $g\in
 F$ such that $\der^i g=f$. Notice that $u\in \End_K(K[x])$, and
 (using Lemma \ref{D30May10}).
 $$ \ind_{K[x]}(u)=\ind_{K[x]}(\der^i)+
 \ind_{K[x]}(s+g)+\ind_{K[x]}(t^{-1})=
 i+\ind_{K[x]}(s)=i\geq 0.$$
The element $u$ is regular and ${}_{\mI_1}F\simeq K[x]^{(\N )}$,
hence $i=0$, i.e.,  $u=(s+g)t^{-1}$ and $\ker_{K[x]}(u)=0$. The two
conditions $\ind_{K[x]}(u)=0$ and $\ker_{K[x]}(u)=0$ imply that
 the map $u_{K[x]}$ is a bijection, and so $u\in \mI_1\cap
 \Aut_K(K[x])=\mI_1^0$.

{\em Case 2}: $u=\int^jst^{-1}+ft_1^{-1}$, $j>0$. By the same
reason as in  Case 1, we may assume that $t=t_1$, i.e.,
$u=(\int^js+f)t^{-1}$. Then (using Lemma \ref{D30May10})
$$ \ind_{K[x]}(u)=\ind_{K[x]}(\int^js+f)+\ind_{K[x]}(t^{-1})=
\ind_{K[x]}(\int^js)=\ind_{K[x]}(\int^j)+\ind_{K[x]}(s)=-j<0.$$
Therefore, we can fix an element, say $v\in K[x]$, such that
$v\not\in \im_{K[x]}(u\cdot )$. Fix an element, say $e\in F$, such
that $\im_{K[x]}(e\cdot ) = Kv$. Then $e\mI_1\cap u\mI_1=0$ since,
otherwise, we get a contradiction: fix $0\neq w\in e\mI_1\cap
u\mI_1$, then $w=ea=ub$ for some elements $a,b\in \mI_1$, and so
$0\neq wK[x]\subseteq \im_{K[x]}(e\cdot ) \cap \im_{K[x]}(u\cdot )
= Kv\cap \im_{K[x]}(u\cdot ) =0$, a contradiction. The equality
$e\mI_1\cap u\mI_1=0$ implies $eS\cap u\mI_1=\emptyset$, this
contradicts to the fact that $S$ is a  right Ore set in $\mI_1$.
This means that $S\subseteq \mI_1^0$, as required.

3. Statement 3 follows from statement 2 using (\ref{IJab}).

4. Suppose that $\v : \mI_1{\mI_1^0}^{-1}\ra {\mJ_1^0}^{-1}\mI_1$
is an $\mI_1$-isomorphism, we seek a contradiction. By Proposition
\ref{a30Sep10}, $\mJ_1^0\subseteq (\mI_1{\mI_1^0}^{-1})^* =
\mI_1^0{\mI_1^0}^{-1}$, by Proposition \ref{A21Jun10}.(1).
Clearly, $1+\der \in \mI_1^0$, then $u:=1+\int = (1+\der )^* \in
(\mI_1^0)^* = \mJ_1^0$, by (\ref{KID3}). By Proposition
\ref{a12Jun10}.(1), $\ind_{K[x]}(1+\int ) = -1$. Since $u\in
\mI_1{\mI_1^0}^{-1}\subseteq \Aut_K(K[x])$, we must have
$\ind_{K[x]}(1+\int ) =0$, a contradiction. $\Box $


\begin{proposition}\label{a30Sep10}
Let $R$ be a ring, $T$ and $S$ be a left and right denominator set
of regular elements of the ring $R$  respectively (and so
$R\subseteq T^{-1}R$ and $R\subseteq RS^{-1}$). Then there is an
$R$-isomorphism of rings $\v : T^{-1}R\ra RS^{-1}$ (i.e.,  $\v (r) =
r$ for all $r\in R$) iff $S\subseteq (T^{-1}R)^*$ and $T\subseteq
(RS^{-1})^*$ where $(T^{-1}R)^*$ and $(RS^{-1})^*$ are the groups
of units of the rings $T^{-1}R$ and $RS^{-1}$ respectively.
\end{proposition}

{\it Proof}. $(\Rightarrow )$ Trivial.

$(\Leftarrow )$ The subring $\CT$ of $RS^{-1}$ generated by $R$
and the set $T^{-1} = \{ t^{-1}\, | \, t\in T\}$ is canonically
isomorphic to the ring $T^{-1}R$ since $T$ is a left denominator
set of regular elements of the ring $R$. Similarly, the subring
$\CC$ of $T^{-1}R$ generated by $R$ and the set $S^{-1} =
\{s^{-1}\, | \, s\in S\}$ is canonically isomorphic to the ring
$RS^{-1}$. Therefore, we have the ring $R$-monomorphisms $\v :
T^{-1}R\ra RS^{-1}$ and $\psi : RS^{-1}\ra T^{-1}R$. Since $\psi
\v (T^{-1}R)=T^{-1}R$ and $ \v \psi (RS^{-1}) = RS^{-1}$, the maps
$\v $ and $\psi$ are surjective, and so $\v $ is an
$R$-isomorphism. $\Box $


\begin{corollary}\label{21Aug10}
 The the inclusion $A_1\ra
\mI_1$ cannot be lifted neither to a ring monomorphism $\Frac
(A_1)\ra \Frac_l (\mI_1)$ nor to $\Frac (A_1)\ra \Frac_r(\mI_1)$.
\end{corollary}

{\it Proof}. The element $\der $ of the Weyl algebra $A_1$ is
invertible in $\Frac (A_1)$ but is not in $\Frac_l(\mI_1)$ and
$\Frac_r(\mI_1)$ since $\der e_{00}=0$. $\Box $


$${\bf Acknowledgements}$$

 The author would like to thank the referee for comments.

\small{

Department of Pure Mathematics

University of Sheffield

Hicks Building

Sheffield S3 7RH

UK

email: v.bavula@sheffield.ac.uk }


\begin{thebibliography}{99}
\bibitem{Bav-cdixintdif} V. V. Bavula, An analogue of the Conjecture of Dixmier is true for the algebra
 of polynomial integro-differential operators. (Arxiv:math.RA: 1011.3009).
\end{thebibliography}
\end{document}